\documentclass[11pt,reqno,a4paper]{amsart}

\usepackage[T1]{fontenc}
\usepackage{lmodern}
\usepackage[letterpaper,left=1.6in,right=1.6in,
  top=1.25in,bottom=1.25in,headheight=14pt,headsep=20pt]{geometry}
\usepackage[protrusion=true]{microtype}

\usepackage{mathtools,amssymb,amsmath,amsthm}
\usepackage{thmtools}
\usepackage{xfrac}
\usepackage{mathrsfs,bm}
\usepackage{enumitem}

\usepackage[dvipsnames]{xcolor}
\usepackage[
  pagebackref,
  colorlinks=true,
  linkcolor=MidnightBlue,
  citecolor=Black,
  urlcolor=MidnightBlue,
  hypertexnames=false
]{hyperref}
\usepackage{doi}
\usepackage[capitalise,nameinlink]{cleveref}
\usepackage{tikz,pgfplots,float}

\mathtoolsset{showonlyrefs=true}
\usetikzlibrary{arrows.meta,calc,decorations.pathreplacing,patterns}
\pgfplotsset{compat=1.18}

\crefformat{equation}{(#2#1#3)}
\Crefformat{equation}{(#2#1#3)}
\crefrangeformat{equation}{(#3#1#4)--(#5#2#6)}
\Crefrangeformat{equation}{(#3#1#4)--(#5#2#6)}
\crefmultiformat{equation}
  {(#2#1#3)}
  { and~(#2#1#3)}
  {, (#2#1#3)}
  { and~(#2#1#3)}
\Crefmultiformat{equation}
  {(#2#1#3)}
  { and~(#2#1#3)}
  {, (#2#1#3)}
  { and~(#2#1#3)}

\newtheorem{thm}{Theorem}[section]
\newtheorem*{thm*}{Theorem}
\newtheorem{lemma}[thm]{Lemma}
\newtheorem{proposition}[thm]{Proposition}
\newtheorem{corollary}[thm]{Corollary}

\theoremstyle{definition}
\newtheorem{definition}[thm]{Definition}
\newtheorem{rmk}[thm]{Remark}
\newtheorem{ipotesi}[thm]{Assumption}

\numberwithin{equation}{section}

\newcommand{\N}{\mathbb{N}}
\newcommand{\Z}{\mathbb{Z}}
\newcommand{\R}{\mathbb{R}}

\newcommand{\vol}{\mathrm{Vol}}

\newcommand{\p}{\mathbf{p}}

\newcommand\res{\mathop{\hbox{\vrule height 7pt width .3pt depth 0pt\vrule height .3pt width 5pt depth 0pt}}\nolimits}

\renewcommand{\flat}{\mathscr{F}}

\newcommand{\F}{\mathbf{F}}
\newcommand{\M}{\mathbf{M}}

\newcommand{\bF}{\mathbb{F}}

\newcommand{\V}{\mathbf{v}}

\newcommand{\bI}{\mathbf{I}}

\newcommand{\hau}{\mathcal{H}}

\newcommand{\eps}{\varepsilon}
\newcommand{\spt}{\mathrm{spt}}

\renewcommand{\epsilon}{\varepsilon}

\def\XXint#1#2#3{{\setbox0=\hbox{$#1{#2#3}{\int}$ }
  \vcenter{\hbox{$#2#3$ }}\kern-.6\wd0}}

\newcommand{\pa}{\partial}

\def\a#1{[\![#1]\!]} 

\newcommand{\de}{\partial}

\newcommand{\W}{\mathcal W}

\makeatletter
\renewcommand{\tocsection}[3]{%
  \indentlabel{\@ifnotempty{#2}{\bfseries\ignorespaces\makebox[\@ifempty{#1}{25pt}{75pt}][l]{#1 #2\quad}}}\bfseries#3}
\renewcommand{\tocsubsection}[3]{%
  \indentlabel{\@ifnotempty{#2}{\ignorespaces\makebox[30pt][l]{#1 #2\quad}}}#3}
\newcommand\@dotsep{4}
\def\@tocline#1#2#3#4#5#6#7{\relax
  \ifnum #1>\c@tocdepth \else
    \par \addpenalty\@secpenalty\addvspace{#2}%
    \begingroup \hyphenpenalty\@M
    \@ifempty{#4}{\@tempdima\csname r@tocindent\number#1\endcsname\relax}{\@tempdima#4\relax}%
    \parindent\z@ \leftskip#3\relax \advance\leftskip\@tempdima\relax
    \rightskip\@pnumwidth plus1em \parfillskip-\@pnumwidth
    #5\leavevmode\hskip-\@tempdima{#6}\nobreak
    \leaders\hbox{$\m@th\mkern \@dotsep mu\hbox{.}\mkern \@dotsep mu$}\hfill
    \nobreak
    \hbox to\@pnumwidth{\@tocpagenum{\ifnum#1=1\bfseries\fi#7}}\par
    \nobreak
    \endgroup
  \fi}
\def\l@section{\@tocline{1}{0pt}{1pc}{25pt}{}}
\def\l@subsection{\@tocline{2}{0pt}{2pc}{30pt}{}}
\makeatother

\begin{document}

\title[Anisotropic Weyl law]{The anisotropic Weyl law and Grassmannian equidistribution}
\author[A. De Rosa]{Antonio De Rosa}
\address{Department of Decision Sciences and BIDSA, Bocconi University, Milano, Italy}
\email{antonio.derosa@unibocconi.it}

\author[A. Halavati]{Aria Halavati}
\address{Department of Decision Sciences and BIDSA, Bocconi University, Milano, Italy}
\email{aria.halavati@unibocconi.it}


\begin{abstract}
We establish a Weyl law for the codimension-one anisotropic volume spectrum on compact Riemannian manifolds. We then use this law to obtain equidistribution in the full Grassmannian bundle while keeping the ambient metric fixed. More precisely, on every closed connected smooth Riemannian manifold of dimension three through six, we construct smooth, embedded hypersurfaces stationary for anisotropies smoothly converging to area. The cumulative sum of their area measures, when lifted to the Grassmannian bundle, converges to the uniform measure.
\end{abstract}

\maketitle

\section{Introduction}\label{sec:introduction}

The existence and distribution of closed minimal hypersurfaces are central
questions in geometric analysis. The topology of the space of cycles,
identified by Almgren~\cite{Almgren-homotopy}, provides the basis for a
variational approach to these questions. Combined with the min-max
construction of Pitts~\cite{Pitts1981} and the regularity theory of
Schoen-Simon~\cite{SS1981}, it produces smooth closed embedded minimal
hypersurfaces in every closed Riemannian manifold of dimension between
three and seven. Multiparameter families of cycles encode further
topological information and give rise to a sequence of variational
invariants, the volume spectrum. For a bumpy metric on a closed manifold of dimension three through seven, Zhou~\cite{Zhou} proved that every width in the volume spectrum admits a realization by a disjoint union of smooth, closed, embedded, two-sided minimal hypersurfaces, each with multiplicity one. Such metrics form a residual set in the space of smooth
Riemannian metrics by White's bumpy metrics
theorem~\cite{White1991,WhiteBumpy2017}.

The asymptotic behavior of this spectrum connects min-max theory with
the global geometry of the ambient manifold. The estimates of
Gromov~\cite{Gromov1988Dimension} and Guth~\cite{Guth2009Minimax} show
that the $p$-width of a fixed compact $d$-dimensional manifold grows at
the rate $p^{1/d}$. Liokumovich-Marques-Neves~\cite{Weyl-LMN}
identified the leading coefficient by proving the Weyl law
\begin{equation}\label{eq:intro-isotropic-weyl}
 \lim_{p\to\infty}p^{-1/d}\omega_p^{d-1}(M,1)
 =a_d\operatorname{Vol}_g(M)^{(d-1)/d},
\end{equation}
where $a_d>0$ depends only on the dimension. Thus the large-index
behavior of a nonlinear variational spectrum recovers the volume of
the manifold, in analogy with the classical Weyl law for eigenvalues.

These ideas have led to substantial progress on the abundance of
minimal hypersurfaces. In particular Yau~\cite[Problem~88]{Yau1982} conjectured that every
closed Riemannian three-manifold contains infinitely
many smooth, closed, immersed minimal surfaces. Marques-Neves~\cite{MN2017PositiveRicci}
proved the existence of infinitely many such hypersurfaces under
positive Ricci curvature, and Song~\cite{Song2023} established the
corresponding result for arbitrary closed metrics in dimensions three
through seven (see also \cite{Li2023,ChodoshMantoulidis2020,ChodoshMantoulidis2023}). The Weyl law also makes it possible to address their
distribution. Irie-Marques-Neves~\cite{IMN2018} proved that, for a
$C^\infty$-generic metric in these dimensions, the union of the closed
smooth embedded minimal hypersurfaces is dense. Subsequently,
Marques-Neves-Song~\cite{MNS} constructed a sequence whose cumulative
area measures equidistribute with respect to normalized Riemannian
volume. Their theorem also controls tangent-plane averages tested
against traces of symmetric two-tensors. Equidistribution coexists with concentration phenomena. Song-Zhou~\cite{SongZhou2021} proved that, for a $C^\infty$-generic metric in dimensions three through seven, every connected, closed, embedded, two-sided, stable minimal hypersurface is approached by a sequence of closed embedded minimal hypersurfaces whose areas diverge and whose normalized varifolds converge to its normalized varifold. 

The natural next question concerns the simultaneous distribution of
positions and tangent planes. Write $G_{d-1}(M)$ for the bundle of
unoriented hyperplanes, and let $\sigma_x$ be the invariant
probability measure on $G_{d-1}(T_xM)$. Full Grassmannian
equidistribution asks for convergence against every continuous
function $h$ on $G_{d-1}(M)$, with limiting probability measure $\Lambda_g$ satisfying:
\begin{equation}\label{eq:intro-grassmann-measure}
 \int_{G_{d-1}(M)} h\,d\Lambda_g
 :=\frac{1}{\operatorname{Vol}_g(M)}
   \int_M\int_{G_{d-1}(T_xM)}h(x,T)\,d\sigma_x(T)\,dV_g(x).
\end{equation}
The question of extending equidistribution to the Grassmannian bundle
is explicitly raised by Marques-Neves-Song in
\cite[Introduction, final paragraph]{MNS}. It asks for substantially more than
the tensor-trace conclusion: in codimension one,
$\operatorname{tr}_T k=\operatorname{tr}_g k-k(\nu,\nu)$ for a unit
normal $\nu$, so those tests see only constant and quadratic angular
components. A general continuous function on the hyperplane bundle
can contain angular components of arbitrarily high degree.

In this paper, we develop an anisotropic approach to this problem.
We keep the Riemannian metric fixed and allow the surface energy to
depend on the tangent plane. More precisely, for a positive function
$F:G_{d-1}(TM)\to(0,\infty)$, we consider
\[
 \mathcal A_F(\Sigma)
 =\int_\Sigma F(x,T_x\Sigma)\,dA_g.
\]
The area functional corresponds to $F\equiv1$. The integrand $F$ can be also seen as the
positively one-homogeneous integrand defined on the tangent bundle
$\varphi_x(\xi)=|\xi|_gF(x,\xi^\perp)$ for $x \in M$ and $\xi\ne0 \in T_xM$, with
$\varphi_x(0)=0$. Convexity of this extension is the natural
assumption. Varying $F$ permits arbitrary smooth angular tests, while
retaining the background metric and the measure
\eqref{eq:intro-grassmann-measure}.

Our first main result is a Weyl law for the anisotropic volume
spectrum. We use
$\mathbb Z_2$ coefficients, relative cycles when the ambient manifold
has boundary, and sweepouts continuous in the strong mass norm. First we define the Weyl constant for frozen anisotropies and prove that the limit exists:
\begin{thm}
\label{thm:main-frozen-weyl}
Let $H$ be an autonomous integrand satisfying
\ref{H:comparability}-\ref{H:ellipticity}.
The limit
\begin{equation}\label{eq:frozen-weyl-constant}
 a(H):=\lim_{p\to\infty}
 p^{-1/d}\omega_p^{d-1}((0,1)^d,H)=\sup_{q\geq1}
 q^{-1/d}\omega_q^{d-1}((0,1)^d,H)
\end{equation}
exists and belongs to $(0,\infty)$.
We call $a(H)$ the frozen Weyl constant associated with $H$. Moreover, for every bounded connected Lipschitz domain
$\Omega\subset\R^d$, we have
\begin{equation}\label{eq:frozen-domain-weyl-law}
 \lim_{p\to\infty}
 p^{-1/d}\omega_p^{d-1}(\Omega,H)=a(H)|\Omega|^{(d-1)/d}.
\end{equation}
\end{thm}
Rotation invariance of $a$ allows us to compute $a(F(x,\cdot))$ in any orthonormal frame at $x$. The resulting coefficient is intrinsic. The following theorem extends the frozen Weyl law to variable integrands on manifolds. 

\begin{thm}
\label{thm:main-weyl}
Let $(M^d,g)$ be compact, connected, and smooth, with possibly empty
smooth boundary, and let $d\geq2$. Assume that $F$ satisfies
\ref{H:comparability}-\ref{H:ellipticity}. Then
\begin{equation}\label{eq:intro-anisotropic-weyl}
\lim_{p\to\infty}p^{-1/d}\omega_p^{d-1}(M,F)
 =\mathcal W(F)
 :=\left(\int_M a(F(x,\cdot))^{d/(d-1)}\,dV_g(x)\right)^{(d-1)/d}.
\end{equation}
\end{thm}
In the isotropic
case, \eqref{eq:intro-anisotropic-weyl} reduces to \eqref{eq:intro-isotropic-weyl}. Both Theorem \ref{thm:main-frozen-weyl} and Theorem \ref{thm:main-weyl} are proved in Part~\ref{part:weyl}.

The second main result is the following Theorem \ref{thm:main-equidistribution}, which concerns the equidistribution of stationary hypersurfaces. Part~\ref{part:equidistribution} proves this conclusion under the assumption of the smooth realization of any width, which in turn is established in Theorem \ref{thm:main-realization} in Part~\ref{part:realization}.

\begin{thm}\label{thm:main-equidistribution}
Let $(M^d,g)$ be a closed connected smooth Riemannian
manifold, with $3\leq d\leq6$, and let $\Lambda_g$
be the probability measure defined in
\eqref{eq:intro-grassmann-measure}.
There exist smooth positive strictly elliptic anisotropies $F_j\rightarrow1$ in $C^\infty(G_{d-1}(TM))$,
and smooth, closed, embedded, connected
$F_j$-stationary hypersurfaces $\Sigma_j$ such that, for every $h\in C^0(G_{d-1}(TM))$,
\begin{equation}\label{eq:intro-equidistribution}
\begin{gathered}
\lim_{q\to\infty}\frac{\displaystyle\sum_{j=1}^q\int_{\Sigma_j}h(x,T_x\Sigma_j)\,dA_g}{\displaystyle\sum_{j=1}^q\operatorname{area}_g(\Sigma_j)}
=\int_{G_{d-1}(TM)}h\,d\Lambda_g.
\end{gathered}
\end{equation}
Moreover, for every $j$, there is an integer $p_j\geq1$ such that
$\Sigma_j$ is a nondegenerate component of a smooth
varifold realizing $\omega_{p_j}^{d-1}(M,F_j)$.
\end{thm}

Thus we obtain full Grassmannian equidistribution through hypersurfaces stationary for energies converging smoothly to area. This result also demonstrates how the theory of anisotropic minimal surfaces can be applied to questions arising in the classical area setting. It is also an anisotropic approximation statement in the direction of the question raised by Marques-Neves-Song in \cite[Introduction, final paragraph]{MNS}. The hypersurfaces $\Sigma_j$ are stationary for their respective integrands $F_j$; they are not asserted to be minimal for area. On the other hand, the conclusion holds for every fixed metric $g$, without a genericity assumption on the metric.

Part~\ref{part:realization} proves a key tool to achieve Theorem \ref{thm:main-equidistribution}: the smooth realization of any width.

\begin{thm}
\label{thm:main-realization}
For every fixed closed smooth $(M^d,g)$, $3 \leq d\leq 6$, and every fixed integer
$k_0\geq6$, there exists $\delta_0>0$ such that every smooth
$F:G_{d-1}(TM)\to(0,\infty)$ with
$\|F-1\|_{C^{k_0}}<\delta_0$ and every $p\geq1$ admit an integral
varifold
\[
 V=\sum_{i=1}^s m_i|\Sigma_i|,\qquad
 \int F\,dV=\omega_p^{d-1}(M,F),
\]
where $m_i\in\mathbb N$ and the $\Sigma_i$ are smooth closed
embedded connected pairwise disjoint $F$-stationary hypersurfaces. $\delta_0$ is indeed independent of $p$.
\end{thm}

The smooth realization requires developing a regularity theory adapted to
anisotropic min-max hypersurfaces. This subject presents additional
difficulties. The regularity theory for minimizing hypersurfaces of
Schoen-Simon-Almgren~\cite{SSA1977} gives weaker bounds on the
singular set than the corresponding theory for area, and the
usual monotonicity formula for the density ratio is unavailable for a general
anisotropic stationary varifold \cite{Allard1974}. The lack of monotonicity is a serious obstruction in the anisotropic minimal surface theory \cite{DeR-survey, DeRosa2026Lectures}. De Philippis-De Rosa~\cite{DPDR2024} developed a min-max construction of anisotropic minimal and CMC surfaces in $3$-manifolds. In the isotropic case, the min-max theory for CMC hypersurfaces was previously developed by
Zhou-Zhu~\cite{ZhouZhu2019} in closed Riemannian manifolds of dimension three through seven. De Philippis-De Rosa-Li~\cite{DPDRLi2024} subsequently established existence and regularity results in arbitrary
dimension. A complementary
approach through anisotropic Allen-Cahn energies was developed by
De Rosa-Pigati~\cite{DeRosaPigati2025}, producing possibly singular
stationary integral varifolds. The multiplicity of these varifolds has been analyzed in the recent work \cite{ChenLan2026}. More recently,
De Rosa-Halavati-Wang~\cite{DRHW2026} established anisotropic
minimization in isotopy classes and a Simon-Smith construction with
genus control, together with isolated-singularity removal under
pinching assumptions. Their removable-singularity theorem also
applies to the $3$-dimensional multiparameter Almgren-Pitts
construction.

In Theorem~\ref{thm:main-realization} we advance the aforementioned theory, by proving that on a fixed closed $d$-manifold, $3\leq d \leq 6$, for every anisotropy in a neighborhood of area independent of the index $p$, every width is realized by a smooth embedded stationary varifold. To this aim, we refine the regularity theory for minimizing hypersurfaces in a neighborhood of area, and we establish curvature
estimates, and isolated-singularity removal for stable hypersurfaces.

\subsection{Organization and main ideas of the proof}

 The paper is organized in three parts.
Part~\ref{part:weyl} determines the anisotropic Weyl law, that is, the asymptotic behavior of the
anisotropic widths. Part~\ref{part:equidistribution} converts this
spectral information into equidistribution of positions and tangent
planes in the Grassmannian bundle, assuming that the widths are attained by finite unions of
smooth stationary hypersurfaces with positive integer multiplicities.
Part~\ref{part:realization} establishes this smooth realization
property in dimension three through six.
We explain below the proof strategy of these three parts.

\medskip
\noindent\textbf{Part 1: The anisotropic Weyl law.}
The first task is to set up a class of sweepouts that can be localized
and assembled without losing control of either their topology or their
anisotropic energy. We work with relative $\mathbb Z_2$ cycles and
sweepouts continuous in the strong mass topology. Proposition \ref{prop:mass-approximation} is an approximation
result that converts the flat families arising
in geometric constructions into admissible strong mass families, with arbitrarily small increase of the maximal anisotropic energy.

First we establish the Weyl law for frozen anisotropies,
using the self-similarity argument of \cite{Weyl-LMN}.
The Lusternik-Schnirelmann inequality in
Proposition~\ref{prop:LS-anisotropic} bounds the width
of a domain from below in terms of the widths of disjoint
subdomains. Tiling and packing then give the frozen
Weyl law of Theorem~\ref{thm:main-frozen-weyl}.
The constant $a(H)$ is independent of the choice of
orthonormal coordinates by
Corollary~\ref{cor:frozen-weyl-rotation-invariance}.

For a non-autonomous integrand, we freeze the coefficients
on small coordinate cubes and apply the variable
Lusternik-Schnirelmann inequality of
Proposition~\ref{prop:LS-variable-anisotropy}.
Taking the limit gives the lower bound in
Corollary~\ref{cor:variable-weyl-lower}.
In the lower bound, the asymptotically optimal allocation of the index is governed by
\[
 \frac{p_i}{p}\simeq
\frac{\vol_g(C_i)\,a(H_i)^{d/(d-1)}}{\sum_j\vol_g(C_j)\,a(H_j)^{d/(d-1)}}.
\]
This gives the integral and the exponent appearing
in $\mathcal W(F)$.

For the upper bound, the self-similar argument no longer
suffices, so we use a variational construction.
First we cover the manifold by small coordinate cubes,
freeze the anisotropy at their centers, and place normalized
copies of the cubes in a long rectangle, recording the
corresponding volume factors; see
Figure~\ref{fig:manifold-cubes-to-rectangle}.
The goal is to combine local sweepouts while controlling
both their energy and the resulting sweepout index.

For each cube, choose an almost optimal $p$-sweepout with
parameter space $X_i$ and pass to its filling double cover
$\widetilde X_i$.
We take the product of these covers with the diagonal
involution and, for a large integer $K$, form the join
\[
 (\widetilde X_1\times\cdots\times\widetilde X_N)^{*K}.
\]
The join parameters determine a slab decomposition of the
long rectangle; see Figure~\ref{fig:join-subdivision}.
For each nonempty intersection of a slab with a cube,
we transform the assigned filling using the cube-to-slab
map of Proposition~\ref{prop:cube-to-slab-map}.
This map first applies a homothety and then cuts and
translations. The homothety gives the required scaling
of the anisotropic area, while the rearrangements preserve
the energy within each piece and add a controlled interface
cost. An interactive visualization is available
in \cite{HalavatiCubeToSlab}.

The resulting family descends to the quotient by the
involution. Lemma~\ref{lem:join-construction-is-higher-sweepout}
shows that it is a $K(p+1)-1$ sweepout.
Lemma~\ref{lem:join-construction-area-calculation}
bounds its maximal energy by an allocation problem
for the numbers of slabs meeting each cube.
Taking first $p\to\infty$, then $K\to\infty$, and finally
refining the coordinate cover gives the matching upper
bound in Theorem~\ref{thm:variable-weyl-upper}.
Together with the lower bound, this proves
Theorem~\ref{thm:main-weyl}.

\medskip
\noindent\textbf{Part 2: Grassmannian equidistribution.}
Here we prove
Theorem~\ref{thm:fixed-metric-anisotropic-equidistribution}
under Assumption~\ref{hyp:near-area-realization}, which
ensures that every width is realized by a finite collection
of smooth stationary hypersurfaces.
The goal is to use the Weyl law to construct a sequence
of such hypersurfaces whose normalized cumulative area
measures converge to $\Lambda_g$.
For this we keep the metric fixed and vary the anisotropic
integrand. Indeed, differentiating the energy with respect
to a perturbation $h(x,T)$ gives its integral over the
hypersurface, so these variations test both its position
and its tangent plane.

First we identify the corresponding derivatives of the
Weyl functional $\mathcal W$ near the area integrand.
The difficulty is that the frozen coefficient $a(H)$
is not explicitly known, and we do not assume that it
is differentiable at $H=1$.
Instead, we average derivatives over finite-dimensional
families of perturbations.
Lemma~\ref{lem:grassmann-parseval} constructs these families
so that a radial average in the parameter space induces
a rotation-invariant measure in each fiber.
This avoids choosing a global orthonormal frame of the
tangent bundle.
Rotation invariance of $a$ then identifies the angular
average, while homogeneity determines its normalization.
In Proposition~\ref{prop:global-averaged-weyl-derivative},
we show that, as the perturbations shrink to zero,
the averaged derivative in a direction $h$ converges to
\[
 \mathcal W(1)\int_{G_{d-1}(TM)}h\,d\Lambda_g.
\]

To express derivatives of the widths in terms of stationary
hypersurfaces, we next adapt White's deformation argument
to variations of the integrand in
Proposition~\ref{prop:anisotropic-white-structure}.
After an arbitrarily small smooth perturbation of a
finite-parameter family, every closed embedded stationary
hypersurface is nondegenerate at almost every parameter.
The inverse function theorem and the countable atlas
of the solution space then give countably many local
$C^1$ branches of nondegenerate connected stationary
hypersurfaces.

A varifold realizing a width may contain several connected
components with positive integer multiplicities.
We therefore take all finite collections of these branches
with such multiplicities.
Their energies form a countable family of $C^1$ functions,
and the smooth realization assumption ensures that each
width agrees with one of these functions at almost every
parameter.
Lemma~\ref{lem:countable-smooth-selection} shows that
their differentials also agree almost everywhere.
Consequently,
Proposition~\ref{prop:anisotropic-countable-derivative}
gives a realizing varifold $V_{p,t}$ such that
\[
 D_t\omega_p^{d-1}(M,F(t))[v]
 =
 \int_{G_{d-1}(TM)}
 D_tF(t)[v](x,T)\,dV_{p,t}(x,T)
 \qquad\text{for every }v.
\]
When differentiating the energy of the realizing collection,
the contribution from the motion of the hypersurfaces
vanishes by stationarity.
The same varifold therefore represents the derivative
in every parameter direction: there is no need of a compactness
argument.

To combine this formula with the Weyl law, we must pass
from convergence of the normalized widths to convergence
of their derivatives.
The Weyl law and comparison of integrands give locally
uniform convergence in each finite-dimensional family,
but this does not imply pointwise convergence of derivatives.
Instead, we integrate against a smooth parameter density
and use integration by parts, as in
Lemma~\ref{lem:averaged-spectral-derivatives}.
We first keep the perturbation scale fixed and let
$p\to\infty$, and then let the scale tend to zero.
The generic perturbation needed for the derivative formula
is controlled by
Lemma~\ref{lem:generic-perturbation-averaging}:
choosing it sufficiently small in $C^1$ relative to the
perturbation scale makes its effect on the averaged
derivatives arbitrarily small.
The estimate is uniform in $p$, so this perturbation
can be fixed before taking the limit in the widths.

Combining these estimates, and approximating the parameter
averages by finite convex combinations,
Proposition~\ref{prop:geometric-grassmann-collections}
produces finite weighted collections of smooth stationary
hypersurfaces whose normalized weighted sums of area measures
approximate $\Lambda_g$ against any prescribed finite family of
test functions.
We then take increasingly many functions from a dense
family in $C(G_{d-1}(TM))$, while making the corresponding
anisotropies converge to area in increasingly strong
$C^k$ norms.

Finally, in Lemma~\ref{lem:grassmann-block-sequence}, we arrange these collections into a single sequence of hypersurfaces. We approximate the weights by rational numbers and replace them by repetitions of the corresponding hypersurfaces.
By choosing the number of repetitions appropriately, the normalized sums of the area measures converge to $\Lambda_g$ as the number of terms tends to infinity, while the total area diverges. Each collection contributes only finitely many terms, so the associated anisotropies still converge to area in $C^\infty$.
This proves Theorem~\ref{thm:fixed-metric-anisotropic-equidistribution}.

\medskip
\noindent\textbf{Part 3: Smooth realization of the anisotropic widths.}
Part~\ref{part:realization} proves smooth realization Theorem~\ref{thm:main-realization} for finite-parameter Almgren-Pitts constructions in ambient dimensions $3\leq d\leq6$: We show that for anisotropies close to the area functional in strong enough topology, the min-max $p$-widths for arbitrary $p$ are achieved by a collection of smooth hypersurfaces (possibly with multiplicity). Indeed the main difficulty is carrying out the regularity theory without a monotonicity formula.

In Section~\ref{subsec:pull-tight} we perform the pull-tight procedure to min-max sequences and show that from every min-max sequence one can extract another sequence that has a limiting min-max varifold which is stationary and locally almost-minimizing away from finitely many points. This step is classical.

In Section~\ref{subsec:curvature-estimate-regularity-replacements} we prove that almost minimizing varifolds admit smooth replacements that glue together smoothly. First in Proposition~\ref{prop:stable-curvature-compactness} we prove interior curvature estimate for stable $F$-stationary hypersurfaces for $F$ sufficiently close to the area in strong enough topology. The proof uses a point-picking argument and Bernstein theorems that are already established for anisotropies sufficiently close to the area (see Lin~\cite{Lin1990},
Chodosh-Li~\cite[Corollary~1.5]{ChodoshLi2023},
and Li-Xia~\cite[Corollary~3]{LiXia2026}). Then in Proposition~\ref{prop:boundary-multiplicity} we use a barrier trick from \cite{Hardt1977} and show boundary regularity if the domain is a thin lens, depending on the multiplicity of the boundary. This result is a precursor to the smooth gluing of consecutive replacements, for which we show the standard existence procedure in Proposition~\ref{prop:smooth-replacements}.

Now for two consecutive replacements on two overlapping annuli, similar to the smooth gluing argument of \cite{DPDRLi2024}, we cover the boundary point with a series of shrinking dyadic overlapping lens, centered at an arbitrary boundary point. Then we use the boundary regularity for these shapes and taking the limit, conclude with smooth gluing of the two original replacements. This is all done in Proposition~\ref{prop:smooth-gluing}.

At this point, as stated in Corollary~\ref{cor:regularity-away-finite-points}, we have shown that the min-max varifold is indeed a smooth hypersurface away from finitely many points.

To remove the singular point in the limit, we use the removable singularity ideas developed in \cite{DRHW2026} using a \textit{separation function}. The main obstacle is to show that the singular points are indeed of finite density or more precisely of zero capacity. Once this is shown we can extend the stability inequality over the singular point and, using a slight modification of the near-area anisotropic calculations of \cite{Winklmann2005}, we can upgrade the stability inequality's exponent to $d-1$ and conclude smoothness using a blow-up method (as in \cite[Appendix A]{DPDRLi2024}).

The idea behind establishing an upper density bound at a singular point is to study the separation between nearby sheets. This function is an approximate supersolution of the anisotropic Jacobi equation. Stability away from the singularity gives a nonnegative quadratic form for the Jacobi operator. For anisotropies sufficiently close to area, the ambient distance $|x|$ from the singular point is a strict subsolution. Comparing the separation function with a small multiple of $|x|$ and applying stability with cutoffs around the singularity shows that the separation cannot vanish faster than linearly. Together with the curvature estimate, this implies the required density bound. This is done in Theorem~\ref{thm:local-removable-singularity}, which, combined with Corollary~\ref{cor:regularity-away-finite-points}, implies the smooth realization theorem~\ref{thm:main-realization} for anisotropies close to area.

\newtheorem*{ack*}{Acknowledgement}
\begin{ack*}
    The authors were funded by the European Union: the European Research Council (ERC), through StG ``ANGEVA'', project number: 101076411. Views and opinions expressed are however those of the authors only and do not necessarily reflect those of the European Union or the European Research Council. Neither the European Union nor the granting authority can be held responsible for them.
\end{ack*}

\medskip
\noindent\textbf{AI Statement.} The initial ideas underlying the join and slab constructions in sections~\ref{subsec:cube-to-slab}-\ref{subsec:long-rectangle-join} emerged from a discussion between the authors and GPT-6 Astra. In particular, the model suggested the original constructions, which the authors subsequently examined, verified, refined, and reorganized into their present form. The authors then came up with the cube-to-slab map and the long rectangle decomposition. In particular the current argument of sections~\ref{subsec:cube-to-slab}-\ref{subsec:long-rectangle-join} is human made. We also used GPT-6 Astra to improve the exposition of the article, including editing help, generating figures and reference checking. In particular figures were drawn by hand and fed into Astra to generate a TikZ file. The authors take full responsibility for the mathematical content, references, figures, and presentation of the manuscript.

\subsection{Notation}

The ambient dimension is denoted by $d$, and all cycles
and hypersurfaces have dimension $d-1$.
We write
\[
 G:=G_{d-1}(TM).
\]
For $U\subset M$, we use
\[
 G_{d-1}(U):=\{(x,T)\in G:x\in U\};
\]
in particular, $G_{d-1}(M)=G$.
We write $dV_g$ and $dA_g$ for Riemannian volume and
hypersurface area.

Unless otherwise stated, currents and cohomology groups
have coefficients in $\mathbb Z_2$.
For a current $T$, its ordinary mass, mass measure,
and associated varifold are denoted by
$\M(T)$, $\|T\|$, and $|T|$, respectively.
We use $\flat$ for flat distance.
Relative cycles, their canonical representatives,
and the sweepout class are defined in
Section~\ref{subsec:sweepouts-widths}.

The symbol $\F$ denotes anisotropic energy, whereas
$\bF$ denotes the varifold metric and the corresponding
metric on cycles used in Part~\ref{part:realization}.
The normal representation of an integrand and the
anisotropic perimeter are defined below.
The quantities $a$, $\mathcal W$, and $\Lambda_g$
retain the meanings given in the introduction.

\part{The anisotropic Weyl law}\label{part:weyl}
\section{Introduction and conventions}\label{sec:preliminaries}
Throughout, $d\geq2$ is the ambient dimension, and all cycles have
dimension $d-1$. We use the notation
\[
 \widetilde\omega_p^{d-1}(M,F)=p^{-1/d}\omega_p^{d-1}(M,F).
\]
In this article we only consider the case of codimension one.

\subsection{Anisotropic integrands}
Let $(M^{d},g)$ be compact, smooth, and connected, with
smooth boundary (possibly empty). Let $G_{d-1}(M)$ be the bundle of unoriented tangent
$(d-1)$-planes:
$$
\begin{aligned}
G_{d-1}(M) &= \bigl\{(x,T):\, x\in M,\,\\
&\qquad T\subset T_xM\text{ is a } (d-1)\text{-dimensional linear subspace}\bigr\},
\end{aligned}
$$
and let $F:G_{d-1}(TM)\to(0,\infty)$ be an anisotropic integrand. For a rectifiable varifold
$V$ and an open set $U$ write
\[
 \F(V,U)=\int_{G_{d-1}(U)}F(x,T)\,dV(x,T).
\]
The normal representation of $F$ is
\begin{equation}\label{eq:normal-extension}
 \varphi_x(\xi)=|\xi|_gF(x,\xi^\perp)\quad(\xi\ne0),
 \qquad \varphi_x(0)=0.
\end{equation}
It is even and positively one-homogeneous. We assume the following conditions on the anisotropy throughout the proof of this part:
\begin{enumerate}[label=(H\arabic*),ref=H\arabic*]
\item\label{H:comparability}(Comparability) There is $0<\lambda\le1$ such that
 $\lambda\le F(x,T)\le\lambda^{-1}$.
\item\label{H:regularity}(Regularity) $F$ is continuous on $G_{d-1}(M)$.
\item\label{H:ellipticity}(Convexity) For each $x$, $\varphi_x$ is convex.
\end{enumerate}
In particular these assumptions hold under the smooth uniformly
elliptic hypotheses used for smooth integrands: $F\in C^3$, a bounded
$C^3$ norm, and
\[
 D^2_{\xi\xi}\varphi_x(\xi)[\zeta,\zeta]
 \ge\lambda\frac{|\zeta|_g^2}{|\xi|_g}
 \quad\text{if }\xi\ne0,\quad \langle\zeta,\xi\rangle_g=0.
\]
 Convexity follows from this inequality and
homogeneity. We will frequently also denote an autonomous integrand on the Euclidean space with $H$.

For a finite-perimeter set $E$ in $U$, define
\[
 P_F(E;U):=\int_{\partial^*E\cap U}F(x,\nu_E^\perp)\,d\hau^{d-1},
 \qquad
 P(E;U):=\hau^{d-1}(\partial^*E\cap U).
\]
We write $\mathcal C(U)$ for the finite-perimeter sets
in $U$, identified modulo null sets and equipped with
the $L^1$ topology.
Since $\varphi_x$ is even, the anisotropic perimeter
is independent of the sign chosen for the normal.

\subsection{Sweepouts and widths}\label{subsec:sweepouts-widths} 
We follow the conventions of \cite{Weyl-LMN}, restricting throughout to codimension one and coefficients in $\Z_2$. For the relative-cycle framework and the corresponding free-boundary min-max theory for the area functional, see also Li-Zhou~\cite{LiZhou2021}. We briefly recall the notation.

For $m\in\N$, we set $I^m=[0,1]^m$. For each integer $j\geq0$,
we define $I(1,j)$ to be the cubical complex with vertices
$\{i3^{-j}: \, 0\leq i\leq3^j\}$, and one-dimensional cells $[i3^{-j},(i+1)3^{-j}]$, with $0\leq i<3^j$.
We define
\[
 I(m,j)=I(1,j)^{\otimes m}.
\]
A cell $\alpha=\alpha_1\otimes\cdots\otimes\alpha_m$ has dimension
$\sum_{i=1}^m\dim\alpha_i$, and we identify it with its support
$\alpha_1\times\cdots\times\alpha_m\subset I^m$.
If $X$ is a cubical subcomplex of $I(m,j_0)$, we define $X(j)$,
for $j\geq j_0$, to consist of the cells of $I(m,j)$ contained
in $X$. We denote its collection of $q$-cells by $X(j)_q$.
We identify each cubical complex with its underlying topological
space.

We fix an isometric embedding of $M$ into some Euclidean space
$\R^Q$. For $0\leq k\leq d$, we denote by
$\mathcal R_k(M;\Z_2)$ the space of $k$-dimensional rectifiable
flat chains with coefficients in $\Z_2$ supported in $M$.
For $k\geq1$, we set
\[
 \bI_k(M;\Z_2)
 :=
 \bigl\{
 T\in\mathcal R_k(M;\Z_2):
 \de T\in\mathcal R_{k-1}(M;\Z_2)
 \bigr\}.
\]
We use the analogous notation for chains supported in $\de M$.
For a rectifiable chain $T$, we denote its mass measure by
$\|T\|$, its mass by $\M(T)$, and its support by $\spt(T)$. In some parts of the paper we also use the notation $\M_F$ to mean the anisotropic mass of a current with respect to $F$. For $S,T\in\mathcal R_{d-1}(M;\Z_2)$, we define the flat distance
by
\[
 \flat(S,T)
 :=
 \inf\left\{
 \begin{aligned}
 &\M(R)+\M(Q):
 T-S=R+\de Q,\\
 &R\in\mathcal R_{d-1}(M;\Z_2),\,
 Q\in\mathcal R_d(M;\Z_2)
 \end{aligned}
 \right\}.
\]

We define
\[
 \mathcal Z_{d-1}(M,\de M;\Z_2)
 :=
 \bigl\{
 T\in\bI_{d-1}(M;\Z_2):
 \spt(\de T)\subset\de M
 \bigr\}.
\]
We identify two elements $S,T$ of this space when
\[
 S-T\in\bI_{d-1}(\de M;\Z_2),
\]
and we denote the resulting quotient by
$\mathcal Z_{d-1,\mathrm{rel}}(M,\de M;\Z_2)$.
We write
\[
 P:\mathcal Z_{d-1}(M,\de M;\Z_2)
 \longrightarrow
 \mathcal Z_{d-1,\mathrm{rel}}(M,\de M;\Z_2)
\]
for the quotient map. When $\de M=\emptyset$, we recover the
usual space $\mathcal Z_{d-1}(M;\Z_2)$ of cycles.

We define the mass and flat distance of relative cycles by
\begin{align*}
 \M(P(T))
 &:=
 \inf_{R\in\bI_{d-1}(\de M;\Z_2)}\M(T+R),\\
 \flat(P(S),P(T))
 &:=
 \inf_{R\in\bI_{d-1}(\de M;\Z_2)}\flat(S+R,T).
\end{align*}
We equip the relative cycle space either with this flat distance
or with the strong mass distance $(\sigma,\tau)\longmapsto\M(\sigma-\tau)$.
For a relative class $\tau=P(T)$, we define its mass measure by
\[
 \|\tau\|(A):=\|T\|(A\setminus\de M)
 \qquad\text{for every Borel set }A\subset M.
\]
This definition does not depend on the representative, and
\cite[Proposition~2.4]{Weyl-LMN} gives
$\M(\tau)=\|\tau\|(M)$.
Thus we do not count mass supported on the boundary.

For $T\in\mathcal R_{d-1}(M;\Z_2)$, we denote by $|T|$ its
associated rectifiable varifold and, consistently with the
preceding subsection, set
\[
 \F(T):=\F(|T|,M).
\]
We define the anisotropic mass of a relative cycle by
\[
 \F(P(T))
 :=
 \inf_{R\in\bI_{d-1}(\de M;\Z_2)}\F(T+R).
\]
Equivalently, we integrate only over the interior:
\[
 \F(P(T))=\F(|T|,M\setminus\de M).
\]
Indeed, representatives agree away from $\de M$, and the
boundary approximation used in
\cite[Proposition~2.4]{Weyl-LMN}, together with
(\ref{H:comparability}), removes the boundary contribution with
arbitrarily small anisotropic cost.
In particular, for relative cycles $\sigma,\tau$, we have
\begin{equation}\label{eq:relative-energy-comparison}
 \lambda\M(\tau)\leq\F(\tau)\leq\lambda^{-1}\M(\tau),
 \qquad
 |\F(\sigma)-\F(\tau)|
 \leq\lambda^{-1}\M(\sigma-\tau).
\end{equation}

We denote by
$\mathcal Z_{d-1,\mathrm{rel}}^0(M,\de M;\Z_2)$ the path
component of zero in the flat topology. Since $M$ is connected,
Almgren's isomorphism \cite{Almgren-homotopy} gives
\[
 \pi_1\bigl(
 \mathcal Z_{d-1,\mathrm{rel}}^0(M,\de M;\Z_2),0
 \bigr)
 \simeq H_d(M,\de M;\Z_2)
 \simeq\Z_2.
\]
We denote by $\bar\lambda$ the unique nonzero element of
$H^1\bigl(\mathcal Z_{d-1,\mathrm{rel}}^0(M,\de M;\Z_2);\Z_2\bigr)$,
where we compute cohomology using the flat topology.

Let $p\in\N$ and let $X$ be a finite cubical subcomplex of
some $I(m,j)$. We call a flat-continuous map
\[
 \Phi:X\longrightarrow
 \mathcal Z_{d-1,\mathrm{rel}}^0(M,\de M;\Z_2)
\]
a \emph{$p$-sweepout} if
\begin{equation}\label{eq:p-sweepout}
 \Phi^*(\bar\lambda^p)
 =
 \bigl(\Phi^*\bar\lambda\bigr)^p
 \ne0
 \qquad\text{in }H^p(X;\Z_2).
\end{equation}
This condition depends only on the flat homotopy class of $\Phi$.

Throughout this paper, we use sweepouts continuous in the
strong mass topology. We denote by $\mathcal P_p^{d-1}(M)$
the collection of all such $p$-sweepouts, allowing both the
parameter complex $X$ and its ambient dimension to vary.
Every map in this collection has no concentration of mass
by \cite[Lemma~2.7]{Weyl-LMN}; explicitly,
\begin{equation}\label{eq:no-mass-concentration}
 \lim_{r\downarrow0}
 \sup_{\substack{x\in X\\y\in M}}
 \|\Phi(x)\|\bigl(B_r^g(y)\bigr)=0,
\end{equation}
where $B_r^g(y)$ denotes the ball in $M$ for the Riemannian
distance. By \ref{H:comparability}, this condition also excludes
concentration of anisotropic mass.

We define the anisotropic $p$-width by
\begin{equation}\label{eq:anisotropic-p-width}
 \omega_p^{d-1}(M,F)
 :=
 \inf_{\Phi\in\mathcal P_p^{d-1}(M)}
 \sup_{x\in\operatorname{dmn}(\Phi)}
 \F(\Phi(x)).
\end{equation}
We retain the normalization
\[
 \widetilde\omega_p^{d-1}(M,F)
 :=
 p^{-1/d}\omega_p^{d-1}(M,F).
\]

\subsection{Approximation results}\label{subsec:approximation}

We follow
\cite[Section~2.9]{Weyl-LMN} and
explain the modification needed to control the anisotropic energy.

In our later constructions, we will obtain families by restricting
and assembling local sweepouts. These operations naturally produce
maps continuous in the flat topology. The next proposition allows
us to work with such families, provided we verify that they have
no concentration of mass. We can then recover strong mass
continuity on the same parameter space, preserve the sweepout
class, and increase the maximal anisotropic energy by an
arbitrarily small amount.

\begin{proposition}
\label{prop:mass-approximation}
Let $F$ satisfy
(\ref{H:comparability}-\ref{H:ellipticity}), let $X$ be a finite
cubical subcomplex of some $I(m,j_0)$, and let
\[
 \Phi:X\longrightarrow
 \mathcal Z_{d-1,\mathrm{rel}}^0(M,\de M;\Z_2)
\]
be continuous in the flat topology and have no concentration
of mass, in the sense of \eqref{eq:no-mass-concentration}.
For every $\eps>0$, there exists a map
\[
 \Phi_\eps:X\longrightarrow
 \mathcal Z_{d-1,\mathrm{rel}}^0(M,\de M;\Z_2)
\]
continuous in the strong mass topology and homotopic to $\Phi$
in the flat topology, such that
\begin{align}
 \sup_{x\in X}\flat(\Phi_\eps(x),\Phi(x))
 &<\eps,
 \label{eq:approximation-flat}\\
 \sup_{x\in X}\F(\Phi_\eps(x))
 &\leq \sup_{x\in X}\F(\Phi(x))+\eps.
 \label{eq:approximation-energy}
\end{align}
In particular, if $\Phi$ is a flat-continuous $p$-sweepout,
then $\Phi_\eps\in\mathcal P_p^{d-1}(M)$.
\end{proposition}

\begin{proof}
We follow \cite[Theorems~2.10-2.11 and Proposition~2.12]{Weyl-LMN}.
In the discretization argument of Appendix~A, local lower
semicontinuity of the anisotropic perimeter
\cite{Reshetnyak1968} gives the energy comparison for the
replaced portions, while comparability controls the added
slicing chains by their ordinary mass. Thus the same
construction, with smaller error tolerances, preserves the
anisotropic energy bound and the mass-fineness estimates.
The Almgren extension preserves this bound by
\eqref{eq:relative-energy-comparison}. Uniform flat approximation
and flat homotopy follow from the cited results, and homotopy
invariance preserves the $p$-sweepout condition.
\end{proof}

We next record the no-concentration property of
mass-continuous families, following
\cite[Lemma~2.7]{Weyl-LMN}. This shows that our admissible
sweepouts automatically satisfy the additional hypothesis
in Proposition~\ref{prop:mass-approximation}.

\begin{proposition}
\label{prop:mass-continuity-no-concentration}
Let $X$ be a finite cubical complex and let
\[
 \Phi:X\longrightarrow
 \mathcal Z_{d-1,\mathrm{rel}}^0(M,\de M;\Z_2)
\]
be continuous in the strong mass topology. Then
\[
 \lim_{r\downarrow0}
 \sup_{\substack{x\in X\\y\in M}}
 \|\Phi(x)\|(B_r^g(y))=0.
\]
In particular, $\Phi$ has no concentration of mass or
anisotropic mass.
\end{proposition}

\begin{proof}
We fix $\eps>0$. By compactness of $X$ and mass continuity,
we choose $x_1,\ldots,x_N\in X$ such that, for every $x\in X$,
there exists $i$ with
\[
 \M(\Phi(x)-\Phi(x_i))<\frac{\eps}{2}.
\]
Since $d-1\geq1$, each rectifiable mass measure
$\|\Phi(x_i)\|$ has no atoms. Compactness of $M$ then gives
\[
 \lim_{r\downarrow0}
 \sup_{y\in M}\|\Phi(x_i)\|(B_r^g(y))=0.
\]
Indeed, a sequence of shrinking balls with masses bounded
away from zero would, after taking a convergent subsequence
of their centers, produce an atom.

Since we have only finitely many measures, we choose $r>0$
such that
\[
 \sup_{y\in M}\|\Phi(x_i)\|(B_r^g(y))
 <\frac{\eps}{2}
 \qquad\text{for }i=1,\ldots,N.
\]
For $x\in X$, we choose $i$ as above and obtain, for every
$y\in M$,
\[
 \|\Phi(x)\|(B_r^g(y))
 \leq
 \|\Phi(x_i)\|(B_r^g(y))
 +\M(\Phi(x)-\Phi(x_i))
 <\eps.
\]
\emergencystretch=2em\relax
This proves the assertion for mass. The corresponding
statement for anisotropic mass follows from $F\leq\lambda^{-1}$.
\end{proof}

Combining the two propositions, we can compute
$\omega_p^{d-1}(M,F)$ using flat-continuous $p$-sweepouts
with no concentration of mass. Indeed, every admissible
mass-continuous sweepout belongs to this larger class, and
Proposition~\ref{prop:mass-approximation} gives the reverse
inequality between the two infima by letting $\eps\downarrow0$.

\section{The frozen anisotropy Weyl law}
The frozen Weyl law follows similarly to the packing and self-similar strategy of \cite{Weyl-LMN}. Here we outline the main adaptations to anisotropies.
\subsection{Preliminaries}\label{subsec:frozen-preliminaries}

We first compare the anisotropic widths with the ordinary volume
spectrum. We then establish the anisotropic version of the
restriction estimates in
\cite[Section~2.14 and Lemma~2.15]{Weyl-LMN}, following the same
construction and explaining the modification needed for the
energy estimate.

\begin{lemma}\label{lem:width-boundedness}
Let $(M^d,g)$ be as in the preceding section, and let $F$
satisfy \ref{H:comparability}. There exists
$C_{\lambda,d,M}>1$ such that
\[
 C_{\lambda,d,M}^{-1}
 \leq \widetilde\omega_p^{d-1}(M,F)
 \leq C_{\lambda,d,M}
 \qquad\text{for every }p\geq1.
\]
The same conclusion holds when $M$ is a bounded connected
Lipschitz domain in $\R^d$, with widths computed relative
to its boundary.
\end{lemma}

\begin{proof}
For every admissible sweepout $\Phi$, we have
\[
 \lambda\sup_{x\in\operatorname{dmn}(\Phi)}\M(\Phi(x))
 \leq
 \sup_{x\in\operatorname{dmn}(\Phi)}\F(\Phi(x))
 \leq
 \lambda^{-1}
 \sup_{x\in\operatorname{dmn}(\Phi)}\M(\Phi(x)).
\]
Taking the infimum over the same collection of sweepouts gives
\begin{equation}\label{eq:width-comparison}
 \lambda\omega_p^{d-1}(M,1)
 \leq\omega_p^{d-1}(M,F)
 \leq\lambda^{-1}\omega_p^{d-1}(M,1).
\end{equation}
We multiply by $p^{-1/d}$ and apply the Gromov-Guth bounds
\cite{Gromov1988Dimension,Guth2009Minimax}; see also
\cite[Theorem~2.8]{Weyl-LMN}.
\end{proof}

For the remainder of this subsection, we fix an autonomous
integrand $H$ satisfying
(\ref{H:comparability}-\ref{H:ellipticity}), and we use $\F$
for its associated energy. We consider bounded connected
Lipschitz domains in $\R^d$ and compute their widths relative
to their boundaries.

Let $R\subset\Omega$ be two such domains.
For a relative cycle $\tau=P(T)$ in $\Omega$, we define
\[
 \F(\tau\res R):=\F(|T|,R).
\]
This quantity does not depend on the representative.
We integrate over the open domain $R$, excluding energy
supported on $\de R$.
The chain $T\res R$ need not have finite boundary mass,
so it need not itself define an integral relative cycle
in $R$. We address this issue using the restriction
construction of \cite[Lemma~2.15]{Weyl-LMN}.

\begin{lemma}
\label{lem:localizing-slicing-lemma}
Let $R\subset\Omega$ be bounded connected Lipschitz domains
in $\R^d$, let $p\geq1$, and let
$\Phi\in\mathcal P_p^{d-1}(\Omega)$ with
$X=\operatorname{dmn}(\Phi)$.
Then the following statements hold.
\begin{itemize}[leftmargin=0.6cm]
\item[(i)]
For every $\eps>0$, there exists a map
$\widehat\Phi:X\longrightarrow
 \mathcal Z_{d-1,\mathrm{rel}}^0(R,\de R;\Z_2)$
continuous in the strong mass topology such that
\begin{equation}\label{eq:localization-energy}
 \F(\widehat\Phi(x))
 \leq(1+\eps)^{d-1}\F(\Phi(x)\res R)
 \qquad\text{for every }x\in X,
\end{equation}
and
\begin{equation}\label{eq:localization-class}
 \widehat\Phi^*\bar\lambda=\Phi^*\bar\lambda
 \qquad\text{in }H^1(X;\Z_2).
\end{equation}
In \eqref{eq:localization-class}, we use $\bar\lambda$ for
the distinguished class of the corresponding target
cycle space.

\item[(ii)]
For every $\eta>0$ and every integer $1\leq q\leq p$,
the set
\[
 \bigl\{
 x\in X:
 \F(\Phi(x)\res R)
 \leq\omega_q^{d-1}(R,H)-\eta
 \bigr\}
\]
is contained in an open set $U\subset X$ such that,
writing $\iota:U\hookrightarrow X$ for the inclusion,
\[
 \bigl(\iota^*\Phi^*\bar\lambda\bigr)^q=0
 \qquad\text{in }H^q(U;\Z_2).
\]
\end{itemize}
\end{lemma}

\begin{proof}
We follow \cite[Lemma~2.15]{Weyl-LMN}. Its inward
deformations $f_t$ converge to the identity in $C^1$.
By continuity and positivity of $H$, for all sufficiently
small $t>0$ we have
\begin{equation}\label{eq:restriction-distortion}
 \F((f_t^{-1})_\#S)
 \leq(1+\eps)^{d-1}\F(S),
 \qquad \spt S\subset f_t(R).
\end{equation}
The slicing and mass-continuity estimates are unchanged,
while lower semicontinuity of the anisotropic energy gives
\eqref{eq:localization-energy} in the limiting construction.
The same Almgren-isomorphism argument gives
\eqref{eq:localization-class}, proving \textup{(i)}.

For \textup{(ii)}, the indicated sublevel set is compact
by mass continuity. If it is empty, take $U=\emptyset$.
Otherwise, apply \textup{(i)} with $\eps$ sufficiently small
and choose a finite cubical neighborhood $Y$ of this set
such that
\[
 \sup_{x\in Y}\F(\widehat\Phi(x))
 <\omega_q^{d-1}(R,H).
\]
The definition of the width implies that
$((\widehat\Phi|_Y)^*\bar\lambda)^q=0$.
Taking $U=\operatorname{int}_X Y$ and using
\eqref{eq:localization-class} proves the assertion.
\end{proof}

\subsection{The Lusternik-Schnirelmann inequality}
\label{subsec:LS-inequality}

\leavevmode\newline
We now prove the anisotropic version of
\cite[Section~3.1]{Weyl-LMN}.
The restriction estimates allow us to compare a global width
with the sum of local widths on disjoint subdomains.
The cup-product argument ensures that one member of the
global sweepout simultaneously satisfies all the required
local energy bounds.

We write $|\Omega|$ for Euclidean volume and call two domains
similar if one is obtained from the other by a translation
and a positive dilation. For convenience, we set
\[
 \omega_0^{d-1}(\Omega,H)=\widetilde\omega_0^{d-1}(\Omega,H):=0.
\]

\begin{proposition}
\label{prop:LS-anisotropic}
Let $H$ be an autonomous integrand satisfying
(\ref{H:comparability}-\ref{H:ellipticity}).
Let $\Omega_0,\ldots,\Omega_N$ be bounded connected
Lipschitz domains in $\R^d$ with unit volume.
For each $i=1,\ldots,N$, let $\Omega_i^*\subset\Omega_0$
be similar to $\Omega_i$, and assume that the domains
$\Omega_i^*$ have pairwise disjoint interiors.
We set
\[
 v_i:=|\Omega_i^*|,
 \qquad
 V:=\min_{1\leq i\leq N}v_i,
 \qquad
 p_i:=\lfloor pv_i\rfloor.
\]
There exists $c=c(\lambda,d,\Omega_0)>0$ such that
\begin{equation}\label{eq:normalized-LS-frozen}
 \widetilde\omega_p^{d-1}(\Omega_0,H)
 \geq
 \sum_{i=1}^N
 v_i\,\widetilde\omega_{p_i}^{d-1}(\Omega_i,H)
 -\frac{c}{pV}
 \qquad\text{for every }p\geq1.
\end{equation}
\end{proposition}

\begin{proof}
We follow the proof of \cite[Section~3.1]{Weyl-LMN},
using Lemma~\ref{lem:localizing-slicing-lemma} in place
of their Lemma~2.15. Additivity of the anisotropic energy
over disjoint interiors and the same relative cup-product
argument give
\begin{equation}\label{eq:unrenormalized-LS-frozen}
 \omega_p^{d-1}(\Omega_0,H)
 \geq \sum_{i=1}^N\omega_{q_i}^{d-1}(\Omega_i^*,H),
 \qquad q_i\geq0,\quad \sum_iq_i\leq p.
\end{equation}
Translation invariance and scaling by $r^{d-1}$ give
the same normalization as in the cited proof.

Taking $q_i=p_i$, we have $p_i/(pv_i)\geq1/2$ whenever
$p_i\geq1$. Thus the preceding inequality and
Lemma~\ref{lem:width-boundedness} imply
\[
 \sum_{p_i\geq1}v_i
 \widetilde\omega_{p_i}^{d-1}(\Omega_i,H)
 \leq 2^{1/d}\widetilde\omega_p^{d-1}(\Omega_0,H)
 \leq 2^{1/d}C_{\lambda,d,\Omega_0}.
\]
The same integer-rounding estimate as in \cite{Weyl-LMN}
therefore bounds the error by
$2^{1/d}C_{\lambda,d,\Omega_0}/(pV)$.
With the convention at index zero, this proves the assertion.
\end{proof}

\subsection{The frozen constant}
\label{subsec:frozen-constant}

We first establish convergence on the unit cube $C=(0,1)^d$
and then pass to bounded connected Lipschitz domains using
the packing arguments of \cite{Weyl-LMN}. We also record
the orthogonal invariance and Lipschitz continuity of the
resulting coefficient.

\begin{proof}[Proof of Theorem~\ref{thm:main-frozen-weyl}]
Fix $q\geq1$ and, for $p\geq q$, set
$m=\lfloor(p/q)^{1/d}\rfloor$.
Partitioning $C$ into $m^d$ equal cubes, the unnormalized
Lusternik-Schnirelmann inequality
\eqref{eq:unrenormalized-LS-frozen} and scaling give
\[
 \widetilde\omega_p^{d-1}(C,H)
 \geq
 \left(\frac{m^dq}{p}\right)^{1/d}
 \widetilde\omega_q^{d-1}(C,H).
\]
Since $m^dq/p\to1$, letting $p\to\infty$ and then taking
the supremum over $q$ yields
\[
 \liminf_{p\to\infty}\widetilde\omega_p^{d-1}(C,H)
 \geq \sup_{q\geq1}\widetilde\omega_q^{d-1}(C,H)
 \geq \limsup_{p\to\infty}\widetilde\omega_p^{d-1}(C,H).
\]
This proves \eqref{eq:frozen-weyl-constant}; positivity
and finiteness follow from Lemma~\ref{lem:width-boundedness}.

For unit-volume Lipschitz domains, the packing arguments
of \cite[Lemmas~3.4-3.6]{Weyl-LMN} apply using
Proposition~\ref{prop:LS-anisotropic}. These constructions
require only translations and positive dilations, which
preserve the autonomous energy and multiply it by
$r^{d-1}$, respectively. They therefore give
\begin{equation}\label{eq:frozen-domain-liminf}
 \liminf_{p\to\infty}\widetilde\omega_p^{d-1}(\Omega,H)
 \geq a(H)
 \geq
 \limsup_{p\to\infty}\widetilde\omega_p^{d-1}(\Omega,H).
\end{equation}
Finally, dilation of a general domain to unit volume gives
\eqref{eq:frozen-domain-weyl-law}.
\end{proof}

The domain independence of the frozen coefficient implies
its invariance under orthogonal changes of coordinates.

\begin{corollary}
\label{cor:frozen-weyl-rotation-invariance}
Let $H$ satisfy
\ref{H:comparability}-\ref{H:ellipticity}.
For every $Q\in O(d)$, we have
\[
 a(H\circ Q)=a(H),
 \qquad (H\circ Q)(T):=H(QT).
\]
\end{corollary}

\begin{proof}
Orthogonal pushforward preserves the sweepout class and gives
\[
 \widetilde\omega_p^{d-1}(C,H\circ Q)
 =
 \widetilde\omega_p^{d-1}(Q(C),H).
\]
Since $|Q(C)|=1$, letting $p\to\infty$ and applying
Theorem~\ref{thm:main-frozen-weyl} proves the assertion.
\end{proof}

\begin{lemma}
\label{lem:frozen-weyl-lipschitz}
Let $H$ and $K$ satisfy
\ref{H:comparability}-\ref{H:ellipticity} with the same
comparability constant $\lambda$. Then
\[
 |a(H)-a(K)|
 \leq a(1)\lambda^{-2}\|H-K\|_\infty.
\]
\end{lemma}

\begin{proof}
Since $H,K\geq\lambda$, we have
\[
 K\leq(1+\|H-K\|_\infty/\lambda)H,
 \qquad
 H\leq(1+\|H-K\|_\infty/\lambda)K.
\]
Comparing the energies of the same sweepouts and passing
to the Weyl limits gives the corresponding inequalities
for $a(H)$ and $a(K)$. Since
$a(H),a(K)\leq\lambda^{-1}a(1)$, the conclusion follows.
\end{proof}

\begin{rmk}
    Consequently, for a continuous variable anisotropy $F$
satisfying \ref{H:comparability}-\ref{H:ellipticity},
the function
$x\longmapsto a(F(x,\cdot))$ is intrinsically defined and continuous.
Indeed, Corollary~\ref{cor:frozen-weyl-rotation-invariance}
makes its definition independent of the orthonormal frame,
and Lemma~\ref{lem:frozen-weyl-lipschitz} gives continuity
in every local orthonormal frame.
\end{rmk}

\section{The variable anisotropy Weyl law}
The main point of departure from the established Weyl law for the area functional in \cite{Weyl-LMN} is the case of variable anisotropies. Even on Euclidean domains, the self similar strategy breaks down and this is the reason we develop the explicit upper bound constructions in this section. The lower bound follows from the same techniques, using the Lusternik-Schnirelmann inequality. For the upper bound we perform a construction which then reduces the problem to an \textit{integer allocation}. This in turn selects the best constants and the integral Weyl law in the limit. 
\subsection{The lower bound}
Here we collect some preliminaries. First we show a readily available improvement of the Lusternik-Schnirelmann inequality for variable anisotropies. This will indeed provide the lower bound for the variable Weyl Law.

\begin{proposition}
\label{prop:LS-variable-anisotropy}
Let $(M^d,g)$ be compact, connected, and smooth, with possibly
empty smooth boundary, and let $F$ satisfy
\ref{H:comparability}-\ref{H:ellipticity}.
Set $C=(0,1)^d$. Let $C_1,\ldots,C_N$ be normal-coordinate
cubes with pairwise disjoint interiors and
$\overline{C_i}\subset\operatorname{int}M$.
Write $x_i$ and $r_i$ for their centers and side lengths,
and define the frozen integrands in the center of each cube $H_i=F(x_i,\cdot)$ in the orthonormal coordinates
defining $C_i$. For $p\geq1$, set:
\begin{itemize}
\item $a_i=a(H_i)$ and $v_i=\vol_g(C_i)$;
\item $S=\sum_{j=1}^N a_j^{d/(d-1)}v_j$;
\item $q_i=a_i^{d/(d-1)}v_i/S$ and
      $p_i=\lfloor pq_i\rfloor$;
\item $v_*=\min_i v_i$ and $r_*=\max_i r_i$.
\end{itemize}
There are $r_0>0$, $C_0=C_0(M,g,d,\lambda)$, and a function
$\eta_{F,g}(r)\to0$ as $r\downarrow0$ (depending on the modulus of continuity of $F$) such that, if
$r_*\leq r_0$, then
\begin{equation}\label{eq:normalized-LS-variable}
 \widetilde\omega_p^{d-1}(M,F)
 \geq
 \frac{S^{(d-1)/d}}{1+\eta_{F,g}(r_*)}
 \sum_{i=1}^N q_i
 \frac{\widetilde\omega_{p_i}^{d-1}(C,H_i)}{a_i}
 -\frac{C_0}{pv_*}.
\end{equation}
\end{proposition}

\begin{proof}
We first extend the restriction argument of
Lemma~\ref{lem:localizing-slicing-lemma} to an interior cube
$R=C_i$ and the variable integrand $F$.
For a relative cycle $\tau=P(T)$ in $M$, write
\[
 \F(\tau\res R):=\F(|T|,R).
\]
This does not depend on the representative and, by
\ref{H:comparability}, is continuous as a function of
$\tau$ in the strong mass topology.

The inward deformations used in the proof of
Lemma~\ref{lem:localizing-slicing-lemma} can be made in the
coordinate chart of $R$. They converge to the identity in
$C^1$. Positivity and uniform continuity of $F$ imply that
pushing a chain back from an inward deformation of $R$
increases its anisotropic energy by a factor tending to one.
The slicing and interpolation argument of that lemma therefore
gives, for every $\theta>0$, a strong mass-continuous
restriction $\widehat\Phi$ of a strong mass-continuous family
$\Phi:X\to\mathcal Z_{d-1,\mathrm{rel}}^0(M,\de M;\Z_2)$
such that
\[
 \F(\widehat\Phi(x))
 \leq(1+\theta)\F(\Phi(x)\res R)
 \quad\text{for every }x\in X,
 \qquad
 \widehat\Phi^*\bar\lambda=\Phi^*\bar\lambda.
\]
Here $\widehat\Phi$ takes values in the relative cycle space
of $R$. As in Lemma~\ref{lem:localizing-slicing-lemma}(ii),
it follows that, for $m\geq1$ and $\delta>0$, the set
\[
 \left\{x\in X:
 \F(\Phi(x)\res R)
 \leq\omega_m^{d-1}(R,F)-\delta\right\}
\]
is contained in an open set $U\subset X$ on which
$\bigl(\Phi^*\bar\lambda|_U\bigr)^m=0$.

We apply this observation to the cubes $C_i$.
Let $\Phi\in\mathcal P_p^{d-1}(M)$, put
$X=\operatorname{dmn}(\Phi)$, and write
$\alpha=\Phi^*\bar\lambda$. We omit indices with $p_i=0$;
if there are no remaining indices, the assertion is immediate.
Since $\sum_iq_i=1$, we have $\sum_i p_i\leq p$.
For $\delta>0$, the preceding restriction argument supplies
open sets $U_i\subset X$ containing
\[
 \left\{x\in X:
 \F(\Phi(x)\res C_i)
 \leq\omega_{p_i}^{d-1}(C_i,F)-\frac{\delta}{N}
 \right\}
\]
and satisfying $(\alpha|_{U_i})^{p_i}=0$.
If the $U_i$ covered $X$, their relative cup product would
give $\alpha^{\sum_i p_i}=0$, and hence $\alpha^p=0$.
This contradicts the $p$-sweepout condition.
Thus some $x\in X$ lies outside every $U_i$.
Disjointness of the cube interiors gives
\[
 \F(\Phi(x))
 \geq\sum_i\F(\Phi(x)\res C_i)
 \geq\sum_i\omega_{p_i}^{d-1}(C_i,F)-\delta.
\]
Taking the infimum over $\Phi$ and letting
$\delta\downarrow0$, we obtain
\begin{equation}\label{eq:LS-variable-unnormalized}
 \omega_p^{d-1}(M,F)
 \geq\sum_i\omega_{p_i}^{d-1}(C_i,F).
\end{equation}

Let $E_i:\R^d\to T_{x_i}M$ be the orthonormal
identification defining $C_i$, and put
\[
 \Psi_i(z)
 =\exp_{x_i}\bigl(r_iE_i(z-(1/2,\ldots,1/2))\bigr),
 \qquad z\in C.
\]
We compare each local width with its frozen counterpart.
For $z\in C$ and $T\in G_{d-1}(\R^d)$, the integrand
pulled back by $\Psi_i$, after removing the dilation factor,
is
\[
 H_{i,r_i}(z,T)
 :=
 r_i^{-(d-1)}J_{d-1}\Psi_i(z,T)\,
 F\bigl(\Psi_i(z),D\Psi_i(z)T\bigr).
\]
Normal-coordinate estimates and continuity of $F$ on the
compact Grassmannian bundle give a function
$\rho_{F,g}(r)\to0$ such that
\[
 \sup_{\substack{z\in C\\T\in G_{d-1}(\R^d)}}
 |H_{i,r_i}(z,T)-H_i(T)|
 \leq\rho_{F,g}(r_*)
 \qquad\text{for every }i.
\]
The function can be chosen uniformly over all admissible
normal-coordinate cubes of side at most $r_*$.
Since $H_i\geq\lambda$ and
$v_i=r_i^d(1+O_g(r_i^2))$, comparison under $\Psi_i$
gives, after decreasing $r_0$ if necessary,
\begin{equation}\label{eq:variable-local-freezing}
 \omega_m^{d-1}(C_i,F)
 \geq
 \frac{v_i^{(d-1)/d}}{1+\eta_{F,g}(r_*)}
 \omega_m^{d-1}(C,H_i)
 \qquad(m\geq1),
\end{equation}
where one may take
$\eta_{F,g}(r)=C_{M,g,\lambda}
\bigl(r^2+\rho_{F,g}(r)\bigr)$.
For $F\in C^1$, the normal-coordinate estimates instead
give $\rho_{F,g}(r)\leq C_{M,g,\|F\|_{C^1}}r$.

Set $\beta=(d-1)/d$ and, for $p_i\geq1$, write
$b_i=\widetilde\omega_{p_i}^{d-1}(C,H_i)$.
Equations~\eqref{eq:LS-variable-unnormalized} and
\eqref{eq:variable-local-freezing} imply
\[
 \widetilde\omega_p^{d-1}(M,F)
 \geq
 \frac{1}{1+\eta_{F,g}(r_*)}
 \sum_{p_i\geq1}
 v_i^\beta\left(\frac{p_i}{p}\right)^{1/d}b_i.
\]
For each index in this sum, let
$t_i=p_i/(pq_i)$. Then $1/2\leq t_i\leq1$ and
\[
 0\leq1-t_i^{1/d}
 \leq1-t_i
 <\frac{1}{pq_i}.
\]
The definitions of $S$ and $q_i$ give the identity
\[
 v_i^\beta q_i^{1/d}
 =S^\beta\frac{q_i}{a_i}.
\]
Consequently, replacing $t_i^{1/d}$ by $1$ in the
preceding lower bound incurs an error of at most
\[
 \frac{S^\beta}{p}
 \sum_{p_i\geq1}\frac{b_i}{a_i}.
\]
By Lemma~\ref{lem:width-boundedness}, the $b_i$ are
uniformly bounded in terms of $d$ and $\lambda$.
Moreover, $a_i\geq\lambda a(1)$,
$S\leq C_{d,\lambda}\vol_g(M)$, and
$Nv_*\leq\sum_i v_i\leq\vol_g(M)$.
The error is therefore at most $C/(pv_*)$.
Using $\widetilde\omega_0^{d-1}(C,H_i)=0$ to restore
the omitted indices proves
\eqref{eq:normalized-LS-variable}.
\end{proof}

The immediate corollary of Proposition \ref{prop:LS-variable-anisotropy} is the lower bound for the Weyl law:
\begin{corollary}
\label{cor:variable-weyl-lower}
Let $(M^d,g)$ and $F$ satisfy
\ref{H:comparability}-\ref{H:ellipticity}. Then
\[
 \liminf_{p\to\infty}\widetilde\omega_p^{d-1}(M,F)
 \geq
 \left(
 \int_M a(F(x,\cdot))^{d/(d-1)}\,dV_g(x)
 \right)^{(d-1)/d}.
\]
\end{corollary}
\begin{proof}
Fix a finite packing as in Proposition~\ref{prop:LS-variable-anisotropy}.
For this packing, $p_i\to\infty$ as $p\to\infty$.
The frozen Weyl law and \eqref{eq:normalized-LS-variable} give
$$
\liminf_{p\to\infty}\widetilde\omega_p^{d-1}(M,F)\geq
 \frac{S^{(d-1)/d}}{1+\eta_{F,g}(r_*)}
 \sum_{i=1}^N q_i
 \lim_{p\to\infty}
 \frac{\widetilde\omega_{p_i}^{d-1}(C,H_i)}{a_i} =\frac{S^{(d-1)/d}}{1+\eta_{F,g}(r_*)}.
$$
Now choose packings whose maximal side length tends to zero and
whose unions exhaust $M$ up to vanishing volume. By continuity
of $x\mapsto a(F(x,\cdot))$,
$$
 S=\sum_{i=1}^N\int_{C_i}a(F(x_i,\cdot))^{d/(d-1)}\,dV_g(x) \longrightarrow \int_M a(F(x,\cdot))^{d/(d-1)}\,dV_g(x).
$$
Finally, (again by continuity) $\eta_{F,g}(r_*)\to0$, which proves the claim.
\end{proof}

\subsection{An equivariant cube to slab map}\label{subsec:cube-to-slab}
We recall the notation $\mathcal C(M)$ for sets of finite perimeter and identify each such set with its associated top-dimensional mod-$2$ current.
In the sequel we will need to repeatedly use a map that takes a sweepout (of sets of finite perimeter) in a cube and with a controlled amount of cuts and translations, fits it into a thin slab with equal volume. Indeed we cannot use any rotation or direct anisotropic scaling, since these type of operations will change the anisotropic mass. We also need this map to be equivariant, meaning that it should commute with set complementing. The interested reader may skip the proof at a first reading and only inspect the figures~\ref{fig:cut-and-paste-1} to~\ref{fig:cut-and-paste-4}. An interactive visualization of the planar cube-to-slab map
and the long-rectangle construction of
Section~\ref{subsec:long-rectangle-join}
is available in~\cite{HalavatiCubeToSlab}. We describe this map in detail in the following proposition:
\begin{proposition}\label{prop:cube-to-slab-map}
There exists a family of maps
\begin{align}
    T_\nu:\mathcal{C}([0,1]^d)
    \longrightarrow\mathcal{C}([0,\nu]\times[0,1]^{d-1}),
    \qquad \nu\in[0,1],
\end{align}
with the following properties:
\begin{enumerate}[leftmargin=*]
    \item Let
    \[
    \iota_\nu:\mathcal{C}([0,\nu]\times[0,1]^{d-1})
    \longrightarrow\mathcal{C}([0,1]^d)
    \]
    be the map induced by set inclusion. Then
    $(\nu,S)\mapsto\iota_\nu(T_\nu(S))$ is jointly continuous
    in the usual topology on $[0,1]$ and the flat topology
    on fillings.
    \item The map is equivariant: for every
    $S\in\mathcal{C}([0,1]^d)$,
    \begin{align}
        T_\nu([0,1]^d\setminus S)
        =([0,\nu]\times[0,1]^{d-1})\setminus T_\nu(S).
    \end{align}
    \item For every fixed autonomous anisotropy $H$ satisfying
    \ref{H:comparability}-\ref{H:ellipticity} and every
    $S\in\mathcal{C}([0,1]^d)$,
    \begin{align}\label{eq:cub-to-slab-estimate}
        \M_H(\pa T_\nu(S))
        \leq\nu^{\frac{d-1}{d}}\M_H(\pa S)+C(d,\lambda),
    \end{align}
    where $C(d,\lambda)>0$ depends only on $d$ and $\lambda$.
    \item
The maps preserve no concentration of mass uniformly in $\nu$.
More precisely, if a family
$\mathscr E\subset\mathcal C([0,1]^d)$ satisfies
\[
 \lim_{r\downarrow0}
 \sup_{\substack{S\in\mathscr E\\y\in\R^d}}
 \|\partial S\|\bigl(B_r(y)\cap(0,1)^d\bigr)=0,
\]
then
\[
 \lim_{r\downarrow0}
 \sup_{\substack{\nu\in[0,1]\\S\in\mathscr E\\y\in\R^d}}
 \|\partial T_\nu(S)\|\bigl(B_r(y)\bigr)=0.
\]
\end{enumerate}
\end{proposition}
\begin{proof}
Given $\rho_1,\rho_2>0$ and $0<\theta<1$, we first show how
to rearrange the rectangle $[0,\rho_1]\times[0,\rho_2]$ into
$[0,\theta\rho_1]\times[0,\rho_2/\theta]$ in four steps.
We disregard the cutting lines and identify sets modulo null
sets throughout.
    
    For $(x,y)\in[0,\rho_1]\times[0,\rho_2]$, define (see Figure~\ref{fig:cut-and-paste-1})
\[
\tau_{\rho_1,\rho_2;\theta}^{(1)}(x,y)
=
\begin{cases}
(x,y),
& \rho_2x+(1-\theta)\rho_1y<\rho_1\rho_2,\\[2pt]
(x-\rho_1,y),
& \rho_2x+(1-\theta)\rho_1y>\rho_1\rho_2.
\end{cases}
\]
The cut joins $(\theta\rho_1,\rho_2)$ to $(\rho_1,0)$.
The upper-right triangle is translated by $(-\rho_1,0)$,
while the remaining piece is fixed. The resulting parallelogram
has vertices
\[
(0,0),\quad (\rho_1,0),\quad
(\theta\rho_1,\rho_2),\quad
(-(1-\theta)\rho_1,\rho_2).
\]
    

\begin{figure}[H]
\centering
\scalebox{0.7}{%
\begin{tikzpicture}[
  y=1.3cm,
  font=\small, line cap=round, line join=round,
  outline/.style={draw=cpInk,line width=.8pt},
  cut/.style={
    draw=cpInk,line width=.7pt,
    dash pattern=on 3pt off 2.5pt
  },
  brace/.style={
    draw=cpInk,line width=.5pt,decorate,
    decoration={brace,amplitude=4pt}
  },
  every node/.style={text=cpInk}
]
\definecolor{cpBlue}{HTML}{466C86}
\definecolor{cpOrange}{HTML}{C77C50}
\definecolor{cpInk}{HTML}{303A44}

\def\L{2.6}
\def\s{0.7}

\coordinate (O) at (0,0);
\coordinate (A) at (\L,0);
\coordinate (B) at (0,\L);
\coordinate (C) at (\L,\L);
\coordinate (P) at (\s,\L);

\fill[cpBlue!13] (O)--(A)--(P)--(B)--cycle;
\fill[cpOrange!18] (P)--(C)--(A)--cycle;
\path[pattern=north east lines,pattern color=cpOrange!65]
  (P)--(C)--(A)--cycle;

\draw[outline] (O)--(A)--(C)--(B)--cycle;
\draw[cut] (P)--(A);

\foreach \p in {O,A,B,C,P}
  \fill[cpInk] (\p) circle (1.3pt);
\node[below left=2pt] at (O) {$O$};

\draw[brace] ($(B)+(0,.18)$)--($(P)+(0,.18)$)
  node[midway,above=6pt] {$\theta\rho_1$};
\draw[brace] (-.22,0)--(-.22,\L)
  node[midway,left=6pt] {$\rho_2$};
\draw[brace] (\L,-.40)--(0,-.40)
  node[midway,below=6pt] {$\rho_1$};

\draw[-{Latex[length=2.6mm]},draw=cpInk,line width=.85pt]
  (3.15,1.3)--(4.65,1.3);

\begin{scope}[xshift=6.9cm]
  \coordinate (Q) at ({\s-\L},\L);
  \coordinate (R) at (0,\L);
  \coordinate (S) at (\s,\L);
  \coordinate (T) at (\L,0);
  \coordinate (U) at (0,0);

  \fill[cpOrange!18] (Q)--(R)--(U)--cycle;
  \path[pattern=north east lines,pattern color=cpOrange!65]
    (Q)--(R)--(U)--cycle;
  \fill[cpBlue!13] (U)--(T)--(S)--(R)--cycle;

  \draw[outline] (Q)--(S)--(T)--(U)--cycle;
  \draw[outline] (R)--(U);

  \foreach \p in {Q,R,S,T,U}
    \fill[cpInk] (\p) circle (1.3pt);

  \draw[brace] ($(Q)+(0,.18)$)--($(R)+(0,.18)$)
    node[midway,above=6pt] {$(1-\theta)\rho_1$};
  \draw[brace] ($(R)+(0,.18)$)--($(S)+(0,.18)$)
    node[midway,above=6pt] {$\theta\rho_1$};
  \draw[brace] (\L,-.40)--(0,-.40)
    node[midway,below=6pt] {$\rho_1$};
  \draw[draw=cpInk!65,line width=.5pt,dashed]
    (S)--({\L+.25},\L);

  \draw[brace] ({\L+.25},\L)--({\L+.25},0)
    node[midway,right=6pt] {$\rho_2$};
\end{scope}
\end{tikzpicture}%
}
\caption[The cut and translation of the map
$\tau^{(1)}_{\rho_1,\rho_2;\theta}$.]{%
The cut and translation of the map
$\tau^{(1)}_{\rho_1,\rho_2;\theta}$.\protect\footnotemark}
\label{fig:cut-and-paste-1}
\end{figure}
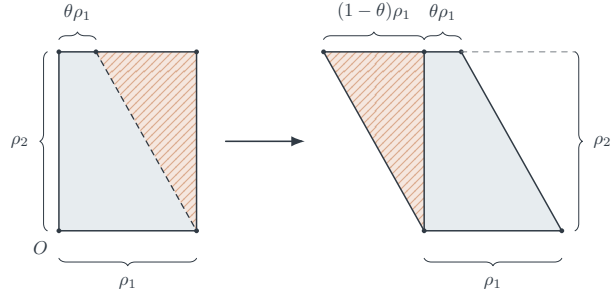
\footnotetext{This figure was converted from a hand-drawn
sketch to TikZ code using ChatGPT.}

On the parallelogram obtained in the first step, define (see Figure~\ref{fig:cut-and-paste-2})
\[
\tau_{\rho_1,\rho_2;\theta}^{(2)}(x,y)
=
\begin{cases}
(x,y),
& \rho_2x<\theta\rho_1y,\\[2pt]
(x-(1-\theta)\rho_1,\,y+\rho_2),
& \rho_2x>\theta\rho_1y.
\end{cases}
\]
The cut joins $(0,0)$ to $(\theta\rho_1,\rho_2)$.
The triangle on its right is translated by
$(-(1-\theta)\rho_1,\rho_2)$, while the other triangle is fixed.
The resulting parallelogram has vertices
$$
(0,0),\qquad (\theta\rho_1,\rho_2),\qquad ((2\theta-1)\rho_1,2\rho_2),\qquad
(-(1-\theta)\rho_1,\rho_2).
$$

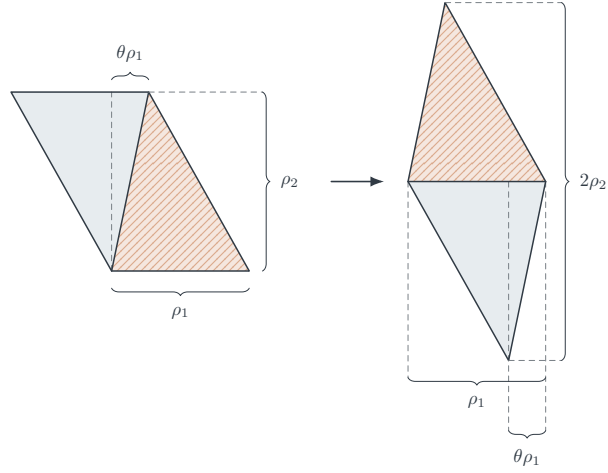
\begin{figure}[H]
\centering
\scalebox{0.7}{%
\begin{tikzpicture}[
  y=1.3cm,
  font=\small, line cap=round, line join=round,
  outline/.style={draw=cpInk,line width=.8pt},
  guide/.style={
    draw=cpInk!65,line width=.5pt,
    dash pattern=on 3pt off 2.5pt
  },
  brace/.style={
    draw=cpInk,line width=.5pt,decorate,
    decoration={brace,amplitude=4pt}
  },
  every node/.style={text=cpInk}
]
\definecolor{cpBlue}{HTML}{466C86}
\definecolor{cpOrange}{HTML}{C77C50}
\definecolor{cpInk}{HTML}{303A44}

\def\L{2.6}
\def\s{0.7}
\def\sideLabel{\nu^{1/d}}
\def\shortLabel{\nu}

\coordinate (A) at (0,\L);
\coordinate (B) at (\L,\L);
\coordinate (C) at ({2*\L-\s},0);
\coordinate (D) at ({\L-\s},0);
\coordinate (E) at ({\L-\s},\L);

\fill[cpBlue!13] (A)--(B)--(D)--cycle;
\fill[cpOrange!18] (B)--(C)--(D)--cycle;
\path[pattern=north east lines,pattern color=cpOrange!65]
  (B)--(C)--(D)--cycle;

\draw[guide] (D)--(E);
\draw[outline] (A)--(B)--(C)--(D)--cycle;
\draw[outline] (D)--(B);

\draw[brace] ($(E)+(0,.18)$)--($(B)+(0,.18)$)
  node[midway,above=6pt] {$\theta\rho_1$};
\draw[brace] ($(C)+(0,-.25)$)--($(D)+(0,-.25)$)
  node[midway,below=6pt] {$\rho_1$};
  \draw[guide] (B)--({2*\L-\s+.25},\L);
\draw[brace]
  ({2*\L-\s+.25},\L)--({2*\L-\s+.25},0)
  node[midway,right=6pt] {$\rho_2$};

\draw[-{Latex[length=2.6mm]},draw=cpInk,line width=.85pt]
  (6.05,1.3)--(7.05,1.3);

\begin{scope}[shift={(7.5,-1.3)}]
  \coordinate (P) at (0,\L);
  \coordinate (Q) at (\s,{2*\L});
  \coordinate (R) at (\L,\L);
  \coordinate (S) at ({\L-\s},0);

  \fill[cpBlue!13] (P)--(R)--(S)--cycle;
  \fill[cpOrange!18] (P)--(Q)--(R)--cycle;
  \path[pattern=north east lines,pattern color=cpOrange!65]
    (P)--(Q)--(R)--cycle;

  \draw[guide] (P)--(0,-.25);
  \draw[guide] (R)--(\L,-1.02);
  \draw[guide] ({\L-\s},\L)--({\L-\s},-1.02);

  \draw[outline] (P)--(Q)--(R)--(S)--cycle;
  \draw[outline] (P)--(R);

  \draw[brace] (\L,-.25)--(0,-.25)
    node[midway,below=6pt] {$\rho_1$};
  \draw[brace] (\L,-1.02)--({\L-\s},-1.02)
    node[midway,below=6pt] {$\theta\rho_1$};
    \draw[guide] (Q)--({\L+.3},{2*\L});
\draw[guide] (S)--({\L+.3},0);
  \draw[brace] ({\L+.3},{2*\L})--({\L+.3},0)
    node[midway,right=6pt] {$2\rho_2$};
\end{scope}
\end{tikzpicture}%
}
\caption[The cut and translation of the map
$\tau^{(2)}_{\rho_1,\rho_2;\theta}$.]{%
The cut and translation of the map
$\tau^{(2)}_{\rho_1,\rho_2;\theta}$.\protect\footnotemark}
\label{fig:cut-and-paste-2}
\end{figure}
\footnotetext{This figure was converted from a hand-drawn
sketch to TikZ code using ChatGPT.}

On the parallelogram obtained in the second step, define (see Figure~\ref{fig:cut-and-paste-3})
\[
\tau_{\rho_1,\rho_2;\theta}^{(3)}(x,y)
=
\left(
 x-\theta\rho_1
   \left\lfloor\frac{x}{\theta\rho_1}\right\rfloor,
 \quad
 y-\rho_2
   \left\lfloor\frac{x}{\theta\rho_1}\right\rfloor
\right).
\]
Away from the vertical cutting lines, the integer
\[
-\left\lfloor\frac{x}{\theta\rho_1}\right\rfloor
\]
counts the strips from right to left, starting at zero.
Thus the rightmost strip is fixed; each successive strip
is translated one additional distance $\theta\rho_1$ to the
right and $\rho_2$ upward. The possibly narrower leftmost
piece is included in the same formula.

The resulting parallelogram is
\[
\left\{
(x,y):
0\le x\le\theta\rho_1,\quad
\frac{\rho_2x}{\theta\rho_1}
\le y\le
\frac{\rho_2x}{\theta\rho_1}+\frac{\rho_2}{\theta}
\right\}.
\]

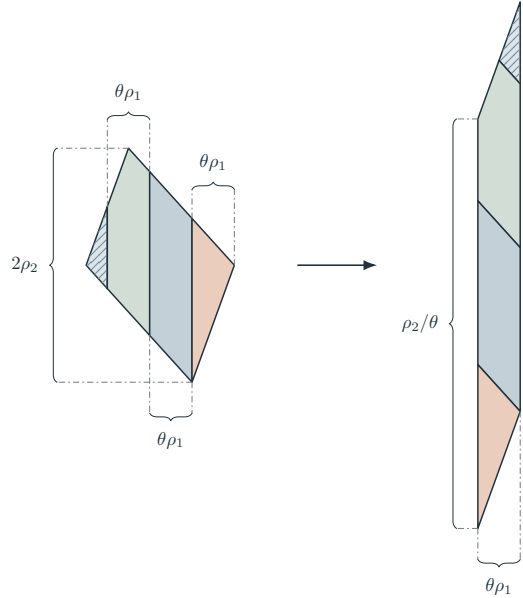
\begin{figure}[H]
\centering
\scalebox{0.7}{%
\begin{tikzpicture}[
  y=1.3cm, font=\small, line cap=round, line join=round,
  piece/.style={draw=cpInk,line width=.8pt},
  pieceA/.style={
    piece,fill=cpBlue!18,
    postaction={
      pattern=north east lines,
      pattern color=cpBlue!75
    }
  },
  pieceB/.style={piece,fill=cpGreen!32},
  pieceC/.style={piece,fill=cpBlue!32},
  pieceD/.style={piece,fill=cpOrange!38},
  guide/.style={
    draw=cpInk!65,line width=.5pt,
    dash pattern=on 3pt off 2pt on .6pt off 2pt
  },
  brace/.style={
    draw=cpInk,line width=.5pt,decorate,
    decoration={brace,amplitude=4pt}
  },
  every node/.style={text=cpInk}
]
\definecolor{cpBlue}{HTML}{466C86}
\definecolor{cpOrange}{HTML}{C77C50}
\definecolor{cpGreen}{HTML}{73977D}
\colorlet{cpInk}{cpBlue!45!black}

\def\L{2.8}
\def\H{1.7}
\pgfmathsetmacro{\w}{2*\L/7}
\pgfmathsetmacro{\res}{\L-3*\w}
\pgfmathsetmacro{\k}{\H*\w/(\L-\w)}
\pgfmathsetmacro{\sideHeight}{\H*\L/\w}
\pgfmathsetmacro{\sourceY}{(\sideHeight-\H)/2}
\pgfmathsetmacro{\arrowY}{(\sideHeight+\H)/2}

\def\pieceApath{
  (0,\H)--(\res,{2*\k})
  --(\res,{\H*(1+\res/\w)})--cycle}
\def\pieceBpath{
  (\res,{2*\k})--({\res+\w},\k)
  --({\res+\w},{\H+2*\k})--(\w,{2*\H})
  --(\res,{\H*(1+\res/\w)})--cycle}
\def\pieceCpath{
  ({\res+\w},\k)--({\L-\w},0)
  --({\L-\w},{\H+\k})
  --({\res+\w},{\H+2*\k})--cycle}
\def\pieceDpath{
  ({\L-\w},0)--(\L,\H)
  --({\L-\w},{\H+\k})--cycle}

\begin{scope}[shift={(0,\sourceY)}]
  \path[pieceA] \pieceApath;
  \path[pieceB] \pieceBpath;
  \path[pieceC] \pieceCpath;
  \path[pieceD] \pieceDpath;

  \draw[guide] (\w,{2*\H})--(-.55,{2*\H});
  \draw[guide] ({\L-\w},0)--(-.55,0);
  \draw[brace] (-.55,0)--(-.55,{2*\H})
    node[midway,left=6pt] {$2\rho_2$};

  \draw[guide]
    (\res,{\H*(1+\res/\w)})
    --(\res,{2*\H+.42});
  \draw[guide]
    ({\res+\w},{\H+2*\k})
    --({\res+\w},{2*\H+.42});
  \draw[brace]
    (\res,{2*\H+.42})
    --({\res+\w},{2*\H+.42})
    node[midway,above=6pt] {$\theta\rho_1$};

  \draw[guide]
    ({\L-\w},{\H+\k})
    --({\L-\w},{\H+\k+.42});
  \draw[guide]
    (\L,\H)--(\L,{\H+\k+.42});
  \draw[brace]
    ({\L-\w},{\H+\k+.42})
    --(\L,{\H+\k+.42})
    node[midway,above=6pt] {$\theta\rho_1$};

  \draw[guide]
    ({\res+\w},\k)--({\res+\w},-.45);
  \draw[guide]
    ({\L-\w},0)--({\L-\w},-.45);
  \draw[brace]
    ({\L-\w},-.45)--({\res+\w},-.45)
    node[midway,below=6pt] {$\theta\rho_1$};
\end{scope}

\draw[-{Latex[length=2.6mm]},draw=cpInk,line width=.85pt]
  (4.0,\arrowY)--(5.5,\arrowY);

\begin{scope}[xshift=7.4cm]
  \begin{scope}[shift={({\w-\res},{3*\H})}]
    \path[pieceA] \pieceApath;
  \end{scope}
  \begin{scope}[shift={({-\res},{2*\H})}]
    \path[pieceB] \pieceBpath;
  \end{scope}
  \begin{scope}[shift={({-(\L-2*\w)},\H)}]
    \path[pieceC] \pieceCpath;
  \end{scope}
  \begin{scope}[shift={({-(\L-\w)},0)}]
    \path[pieceD] \pieceDpath;
  \end{scope}

  \draw[guide] (0,0)--(-.4,0);
  \draw[guide] (0,\sideHeight)--(-.4,\sideHeight);
  \draw[brace] (-.4,0)--(-.4,\sideHeight)
    node[midway,left=6pt] {$\rho_2/\theta$};

  \draw[guide] (0,0)--(0,-.45);
  \draw[guide] (\w,\H)--(\w,-.45);
  \draw[brace] (\w,-.45)--(0,-.45)
    node[midway,below=6pt] {$\theta\rho_1$};
\end{scope}
\end{tikzpicture}
}
\caption[The multiple cuts and translations of
$\tau^{(3)}_{\rho_1,\rho_2;\theta}$.]{%
The multiple cuts and translations of
$\tau^{(3)}_{\rho_1,\rho_2;\theta}$.\protect\footnotemark}
\label{fig:cut-and-paste-3}
\end{figure}
\footnotetext{This figure was converted from a hand-drawn
sketch to TikZ code using ChatGPT.}

On the parallelogram obtained in the third step, define (see Figure~\ref{fig:cut-and-paste-4})
\[
\tau_{\rho_1,\rho_2;\theta}^{(4)}(x,y)
=
\begin{cases}
(x,y),
& y<\rho_2/\theta,\\[2pt]
(x,y-\rho_2/\theta),
& y>\rho_2/\theta.
\end{cases}
\]
The horizontal cut is $y=\rho_2/\theta$.
The triangular cap above this line is translated downward
by $(0,-\rho_2/\theta)$, filling the triangular gap at the bottom.
The resulting rectangle is
\[
[0,\theta\rho_1]\times
\left[0,\frac{\rho_2}{\theta}\right].
\]

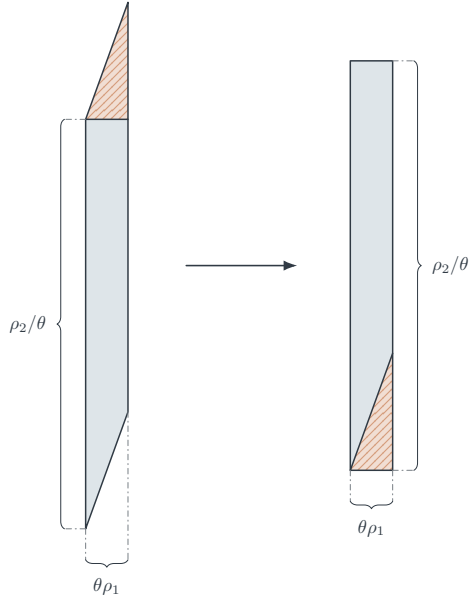
\begin{figure}[H]
\centering
\scalebox{0.7}{%
\begin{tikzpicture}[
  y=1.3cm,
  font=\small, line cap=round, line join=round,
  outline/.style={draw=cpInk,line width=.8pt},
  guide/.style={
    draw=cpInk!65,line width=.5pt,
    dash pattern=on 3pt off 2pt on .6pt off 2pt
  },
  brace/.style={
    draw=cpInk,line width=.5pt,decorate,
    decoration={brace,amplitude=4pt}
  },
  every node/.style={text=cpInk}
]
\definecolor{cpBlue}{HTML}{466C86}
\definecolor{cpOrange}{HTML}{C77C50}
\definecolor{cpInk}{HTML}{303A44}

\def\W{0.8}
\def\V{5.95}
\def\H{1.7}
\pgfmathsetmacro{\midY}{(\V+\H)/2}

\fill[cpBlue!18]
  (0,0)--(\W,\H)--(\W,\V)--(0,\V)--cycle;
\fill[cpOrange!25]
  (0,\V)--(\W,\V)--(\W,{\V+\H})--cycle;
\path[pattern=north east lines,pattern color=cpOrange!75]
  (0,\V)--(\W,\V)--(\W,{\V+\H})--cycle;

\draw[outline]
  (0,0)--(\W,\H)--(\W,{\V+\H})--(0,\V)--cycle;
\draw[outline] (0,\V)--(\W,\V);

\draw[guide] (0,0)--(-.4,0);
\draw[guide] (0,\V)--(-.4,\V);
\draw[brace] (-.4,0)--(-.4,\V)
  node[midway,left=6pt] {$\rho_2/\theta$};

\draw[guide] (0,0)--(0,-.45);
\draw[guide] (\W,\H)--(\W,-.45);
\draw[brace] (\W,-.45)--(0,-.45)
  node[midway,below=6pt] {$\theta\rho_1$};

\draw[-{Latex[length=2.6mm]},draw=cpInk,line width=.85pt]
  (1.9,\midY)--(4.0,\midY);

\begin{scope}[shift={(5,{\H/2})}]
  \fill[cpBlue!18]
    (0,0)--(\W,\H)--(\W,\V)--(0,\V)--cycle;
  \fill[cpOrange!25]
    (0,0)--(\W,0)--(\W,\H)--cycle;
  \path[pattern=north east lines,pattern color=cpOrange!75]
    (0,0)--(\W,0)--(\W,\H)--cycle;

  \draw[outline] (0,0) rectangle (\W,\V);
  \draw[outline] (0,0)--(\W,\H);

  \draw[guide] (\W,0)--({\W+.4},0);
  \draw[guide] (\W,\V)--({\W+.4},\V);
  \draw[brace] ({\W+.4},\V)--({\W+.4},0)
    node[midway,right=6pt] {$\rho_2/\theta$};

  \draw[guide] (0,0)--(0,-.45);
  \draw[guide] (\W,0)--(\W,-.45);
  \draw[brace] (\W,-.45)--(0,-.45)
    node[midway,below=6pt] {$\theta\rho_1$};
\end{scope}
\end{tikzpicture}%
}
\caption[The final cut and translation of
$\tau^{(4)}_{\rho_1,\rho_2;\theta}$ turns the parallelogram
into a rectangle.]{%
The final cut and translation of
$\tau^{(4)}_{\rho_1,\rho_2;\theta}$ turns the parallelogram
into a rectangle.\protect\footnotemark}
\label{fig:cut-and-paste-4}
\end{figure}
\footnotetext{This figure was converted from a hand-drawn
sketch to TikZ code using ChatGPT.}

The resulting planar map
\begin{align}
    \mathcal{T}_{\rho_1,\rho_2;\theta}:
    [0,\rho_1]\times[0,\rho_2]
    \longrightarrow[0,\theta\rho_1]\times[0,\rho_2/\theta]
\end{align}
is defined almost everywhere by
\begin{align}
    \mathcal{T}_{\rho_1,\rho_2;\theta}
    :=\tau^{(4)}_{\rho_1,\rho_2;\theta}
    \circ\tau^{(3)}_{\rho_1,\rho_2;\theta}
    \circ\tau^{(2)}_{\rho_1,\rho_2;\theta}
    \circ\tau^{(1)}_{\rho_1,\rho_2;\theta}.
\end{align}

We now define $T_\nu$. For $0<\nu<1$, set
\begin{align}
    \rho:=\nu^{1/d},
    \qquad
    \phi_\nu:[0,1]^d\longrightarrow[0,\rho]^d,
    \qquad \phi_\nu(x):=\rho x.
\end{align}
The first rearrangement is
\begin{align}
    T_\nu^1(x_1,\ldots,x_d)
    :=\bigl(
    [\mathcal{T}_{\rho,\rho;\rho}(x_1,x_2)]_1,
    [\mathcal{T}_{\rho,\rho;\rho}(x_1,x_2)]_2,
    x_3,\ldots,x_d\bigr),
\end{align}
where $[\cdot]_1$ and $[\cdot]_2$ denote the two coordinate
components. Thus $T_\nu^1\circ\phi_\nu$ sends $[0,1]^d$ onto
$[0,\rho^2]\times[0,1]\times[0,\rho]^{d-2}$ modulo null sets.
When $d\ge3$, the second rearrangement is
$$
T_\nu^2(x_1,\ldots,x_d):=\bigl([\mathcal{T}_{\rho^2,\rho;\rho}(x_1,x_3)]_1,x_2,[\mathcal{T}_{\rho^2,\rho;\rho}(x_1,x_3)]_2,
    x_4,\ldots,x_d\bigr).
$$
The composition $T_\nu^2\circ T_\nu^1\circ\phi_\nu$ sends
$[0,1]^d$ onto
$[0,\rho^3]\times[0,1]^2\times[0,\rho]^{d-3}$ modulo null sets.
In general, for $k=1,\ldots,d-1$, define
$$
\begin{aligned}
    T_\nu^k(x_1,\ldots,x_d)
&:=\bigl([\mathcal{T}_{\rho^k,\rho;\rho}(x_1,x_{k+1})]_1,
    x_2,\ldots,x_k,\\
&\qquad[\mathcal{T}_{\rho^k,\rho;\rho}(x_1,x_{k+1})]_2,
    x_{k+2},\ldots,x_d\bigr),
\end{aligned}
$$
with empty coordinate lists omitted. At step $k$, the side
length in the $e_1$ direction is multiplied by $\rho$, and
the side length in the $e_{k+1}$ direction is multiplied by
$\rho^{-1}$. The successive boxes are therefore
\begin{align}
\begin{aligned}
    {}[0,1]^d
    &\xrightarrow{\phi_\nu}[0,\rho]^d
    \xrightarrow{T_\nu^1}
    [0,\rho^2]\times[0,1]\times[0,\rho]^{d-2}\\
    &\xrightarrow{T_\nu^2}\cdots
    \xrightarrow{T_\nu^k}
    [0,\rho^{k+1}]\times[0,1]^k\times[0,\rho]^{d-k-1}\\
    &\xrightarrow{T_\nu^{k+1}}\cdots
    \xrightarrow{T_\nu^{d-1}}
    [0,\rho^d]\times[0,1]^{d-1}
    =[0,\nu]\times[0,1]^{d-1}.
\end{aligned}
\end{align}
Define the map:
\[
F_\nu:=T_\nu^{d-1}\circ\cdots\circ T_\nu^1\circ\phi_\nu
\]
be the underlying point map. Then for any set of finite perimeter $ S \in\mathcal{C}([0,1]^d)$,
define
\begin{align}
    T_\nu(S):=(F_\nu)_{\#}S,
\end{align}
where the image is taken piecewise on the domains of the
defining formulas, modulo null sets. At the endpoints, set
$T_0(S):=\varnothing$, $T_1(S):=S$, and
$F_1:=\operatorname{id}_{[0,1]^d}$. The equivariance follows immediately since the operation is a set map.

We prove the joint continuity of $\iota_\nu\circ T_\nu$.
Throughout, we identify sets modulo null sets and regard
$T_\nu(S)$ as a subset of $[0,1]^d$ via $\iota_\nu$.
For these top-dimensional fillings, flat convergence agrees
with convergence of characteristic functions in $L^1$.

Since the rearrangements are
translations on finitely many pieces and are bijective almost
everywhere, for $0<\nu<1$ we have
\[
|T_\nu(S)|=\nu|S|,
\qquad
|T_\nu(S)\triangle T_\nu(S')|
=\nu|S\triangle S'|,
\]
where $|\cdot|$ denotes Lebesgue measure.

Fix $\nu_0\in(0,1]$. Locally near $\nu_0$, the number of pieces
is uniformly bounded. Away from the finitely many limiting
cutting hyperplanes, the branch inequalities and the integers
in the floor functions stabilize, while the corresponding
translations vary continuously. At $\nu_0=1$, the limiting
second and fourth translations cancel, and the composite
is the identity. Consequently,
\[
F_\nu(x)\longrightarrow F_{\nu_0}(x)
\quad\text{for almost every }x\in[0,1]^d.
\]
For fixed $S$ and every continuous function $h$ on $[0,1]^d$,
the change-of-variables formula and dominated convergence give
\[
\int_{T_\nu(S)}h(y)\,dy
=\nu\int_S h(F_\nu(x))\,dx
\longrightarrow
\nu_0\int_S h(F_{\nu_0}(x))\,dx
=\int_{T_{\nu_0}(S)}h(y)\,dy.
\]
Since the characteristic functions are bounded by one,
approximating $\mathbf 1_{T_{\nu_0}(S)}$ in $L^1$ by continuous
functions shows that
$|T_\nu(S)\cap T_{\nu_0}(S)|\to |T_{\nu_0}(S)|$.
Together with $|T_\nu(S)|=\nu|S|$, this proves
$|T_\nu(S)\triangle T_{\nu_0}(S)|\to0$.
At $\nu_0=0$, continuity follows uniformly in $S$ from
$|T_\nu(S)|\leq\nu$. Finally,
\[
|T_\nu(S)\triangle T_{\nu_0}(S_0)|
\leq \nu|S\triangle S_0|
+|T_\nu(S_0)\triangle T_{\nu_0}(S_0)|
\]
proves joint continuity.

It remains to prove the estimate in point (3). Write
\[
P_H(E;U):=\M_H(\pa E\llcorner U).
\]
Since $H$ is autonomous, each translation preserves the
anisotropic perimeter inside the translated piece. Moreover in the whole operations, we are never using any rotation, hence every perimeter is transferred identically, and some cut pieces will be added which we will subsequently estimate.

By the BV gluing formula and the comparability assumption,
a new gluing interface contributes at most $C(\lambda)$
times its $(d-1)$-dimensional measure.

In each planar operation, the newly cut edges become exterior
edges of the intermediate target domain. The total length of
the new interior gluing interfaces over the four operations
is bounded by $C\,(\rho_1+\rho_2)$, uniformly in $\theta$.
Indeed, the first two interfaces have lengths $\rho_2$ and
$\rho_1$, respectively. In the third operation, the vertical
cuts become exterior edges, while the gluing interfaces are
paired translated portions of the original slanted sides.
The fourth operation glues the two remaining slanted sides.
Their lengths are bounded by $C\,(\rho_1+\rho_2)$.

Recall that $\rho=\nu^{1/d}$. At step $k=1,\ldots,d-1$, we use
$\rho_1=\rho^k$, $\rho_2=\rho$, and $\theta=\rho$.
The untouched coordinates have total transverse volume
$\rho^{d-k-1}$, so this step increases the relative perimeter
by at most
\[
C(\lambda)(\rho^k+\rho)\rho^{d-k-1}.
\]
The initial homothety multiplies anisotropic perimeter by
$\rho^{d-1}$. Summing over the rearrangements therefore gives
\[
\begin{aligned}
&P_H\bigl(T_\nu(S);(0,\nu)\times(0,1)^{d-1}\bigr)\\
&\quad\le
\rho^{d-1}P_H(S;(0,1)^d)
+C(\lambda)\sum_{k=1}^{d-1}
(\rho^k+\rho)\rho^{d-k-1}\\
&\quad\le
\nu^{\frac{d-1}{d}}P_H(S;(0,1)^d)
+C(d,\lambda)\nu^{1/d}.
\end{aligned}
\]
Finally, the boundary of the target rectangle has anisotropic
mass bounded by $C(d,\lambda)$, uniformly in $\nu$.
Adding its possible contribution yields
\[
\M_H(\pa T_\nu(S))
\le
\nu^{\frac{d-1}{d}}\M_H(\pa S)+C(d,\lambda).
\]
The endpoint cases follow from $T_0(S)=\varnothing$ and
$T_1(S)=S$.

We finally prove (4). The hypothesis and a finite covering
of $[0,1]^d$ by small balls imply
\[
 \sup_{S\in\mathscr E}\|\partial S\|((0,1)^d)<\infty.
\]
Fix $\delta\in(0,1)$. For $\nu\in[\delta,1]$, the number
of pieces in the construction is uniformly bounded.
On each piece, the underlying map is a homothety of
factor $\nu^{1/d}$ followed by a translation.
All additional boundary contributions lie in a uniformly
bounded number of affine hyperplanes. Hence
\begin{align*}
&\sup_{\substack{\nu\in[\delta,1]\\
                  S\in\mathscr E\\y\in\R^d}}
 \|\partial T_\nu(S)\|(B_r(y))\\
&\quad\leq
 C(d,\delta)
 \sup_{\substack{S\in\mathscr E\\y\in\R^d}}
 \|\partial S\|
 \bigl(B_{r\delta^{-1/d}}(y)\cap(0,1)^d\bigr)
 +C(d,\delta)r^{d-1}.
\end{align*}
For fixed $\delta$, this tends to zero as $r\downarrow0$.

For $0<\nu\leq\delta$, we use the relative perimeter
estimate proved above with $H=1$.
The boundary of the target slab has $2d$ faces, each
contributing at most $C(d)r^{d-1}$ inside any ball
of radius $r$. Therefore
\begin{align*}
 \sup_{\substack{\nu\in[0,\delta]\\
                  S\in\mathscr E\\y\in\R^d}}
 \|\partial T_\nu(S)\|(B_r(y))
\leq
 \delta^{(d-1)/d}
 \sup_{S\in\mathscr E}\|\partial S\|((0,1)^d)
 +C(d)\delta^{1/d}+C(d)r^{d-1},
\end{align*}
where $\nu=0$ contributes nothing.
Combining the two estimates, first letting $r\downarrow0$
and then $\delta\downarrow0$, proves (4).
\end{proof}

\begin{rmk}
    There are many ways to construct the map in Proposition~\ref{prop:cube-to-slab-map} and the construction already done above has the advantage that the multiplier behind the perimeter in the estimate, is indeed one which makes subsequent computations in Lemma~\ref{lem:join-construction-area-calculation} more intuitive to follow. One can imagine to cut the square after appropriate rescaling into a series of smaller $\nu$ rectangles with some leftovers and then place a scaling of the whole sweepout in the left over that is not covered in the target. This would indeed introduce a multiplicative error $(1+o(1))$ depending on the rationality of some power of $\nu$ which then makes the computations later much harder to follow. We chose to present the preceeding construction, since the heart of the proof which is the \textit{integer allocation} estimate in Lemma~\ref{lem:join-construction-area-calculation} would be better presented.
\end{rmk}

\subsection{Competitor sweepouts in a long rectangle}\label{subsec:long-rectangle-join}
In this section we prove the upper bound using a variational construction. 
We first recall the definition of the join of two topological spaces $X,Y$:
$$
    X*Y:=\frac{X\times Y\times[0,1]}{\sim},
$$
where the equivalence relation is generated by
\begin{align*}
    \begin{cases}
        (x,y_1,0)\sim(x,y_2,0)
        &\text{for all }x\in X,\ y_1,y_2\in Y,\\
        (x_1,y,1)\sim(x_2,y,1)
        &\text{for all }x_1,x_2\in X,\ y\in Y.
    \end{cases}
\end{align*}
The use of joins of filling double covers to construct higher
sweepouts is also present in \cite[Section~3]{Dey2025}. Now for $N\in\mathbb{N}$, consider the rectangle
\[
    Q=[0,N]\times[0,1]^{d-1}=\bigcup_{i=1}^N C_i,
\]
where $\{C_i\}_{i=1}^N$ are unit cubes, equipped respectively with frozen
integrands $H_i$ satisfying \ref{H:comparability} and \ref{H:ellipticity}. Let
\[
    \Phi_i:X_i\longrightarrow
    \mathcal Z^0_{d-1,\mathrm{rel}}(C_i,\de C_i;\mathbb Z_2),
    \qquad i=1,\dots,N,
\]
be $p$-sweepouts which are almost optimal for $H_i$:
\begin{align}\label{eq:almost-optimality-sweepout}
    \sup_{x\in X_i} \M_{H_i}(\Phi_i(x)) \leq \omega_p^{d-1}(C_i,H_i) + \epsilon.
\end{align}
To work with sets of finite perimeter and set operations, we first pass to the pullback double covers for all $i=1,\dots,N$
$$
\widetilde X_i:=\left\{(x,E)\in X_i\times\mathcal C(C_i):\partial_{\mathrm{rel}}[E]=\Phi_i(x)\right\},
$$
where $\mathcal C(C_i)$ denotes the space of sets of finite perimeter
in $C_i$, identified up to null sets. The covering projection is
$\pi_i(x,E)=x$, and its deck involution is
\[
    \tau_i(x,E)=(x,C_i\setminus E).
\]
The lifted filling maps
\[
    \widetilde\Phi_i:\widetilde X_i\longrightarrow\mathcal C(C_i),
    \qquad
    \widetilde\Phi_i(x,E)=E,
\]
satisfy
\[
    \partial_{\mathrm{rel}}
    [\widetilde\Phi_i(\widetilde x)]
    =\Phi_i(\pi_i(\widetilde x)),
    \qquad
    \widetilde\Phi_i(\tau_i\widetilde x)
    =C_i\setminus\widetilde\Phi_i(\widetilde x).
\]

We equip the product space $\widetilde X_1\times\dots\times\widetilde X_N$ with the diagonal
involution $$(\widetilde x_1,\dots,\widetilde x_N)\mapsto
(\tau_1\widetilde x_1,\dots,\tau_N\widetilde x_N)$$ and form the join
$$
    \widetilde Z_K
    :=
    (\widetilde X_1\times\dots\times\widetilde X_N)^{*K}
    =
    \underbrace{
        (\widetilde X_1\times\dots\times\widetilde X_N)*\cdots*
        (\widetilde X_1\times\dots\times\widetilde X_N)
    }_{K\text{ times}}.
$$
A point of $\widetilde Z_K$ can be formally represented by
$$
    z=\begin{bmatrix}
        (t_1;\widetilde x^{1}_1,\dots,\widetilde x_N^1)\\
        (t_2;\widetilde x^{2}_1,\dots,\widetilde x_N^2)\\
        \cdots\\
        (t_K;\widetilde x^{K}_1,\dots,\widetilde x_N^K).
    \end{bmatrix},
    \qquad
    t_\ell\ge0,\quad \sum_{\ell=1}^K t_\ell=N,
$$
with $\widetilde x_i^\ell\in\widetilde X_i$. Whenever $t_\ell=0$,
the tuple $(\widetilde x_1^\ell,\ldots,\widetilde x_N^\ell)$ is forgotten.
The induced free involution on the join is
$$
    \tau^{(K)} z
    =
    \begin{bmatrix}
        (t_1;\tau_1\widetilde x^{1}_1,\dots,\tau_N\widetilde x_N^1)\\
        (t_2;\tau_1\widetilde x^{2}_1,\dots,\tau_N\widetilde x_N^2)\\
        \cdots\\
        (t_K;\tau_1\widetilde x^{K}_1,\dots,\tau_N\widetilde x_N^K).
    \end{bmatrix}.
$$
The parameter space for the resulting family of relative cycles
will be $Z_K:=\widetilde Z_K/\langle\tau^{(K)}\rangle$.

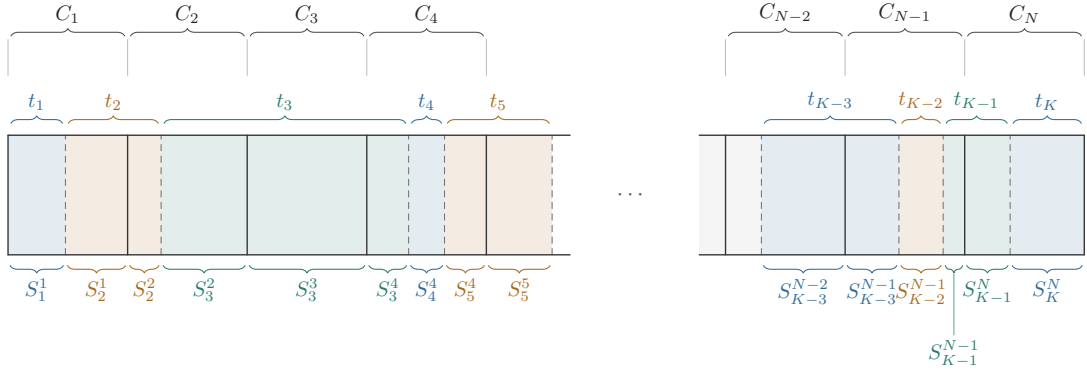
\begin{figure}[H]
\centering
\definecolor{slabBlue}{HTML}{316A9B}
\definecolor{slabOchre}{HTML}{AB681D}
\definecolor{slabTeal}{HTML}{287B70}

\resizebox{0.95\linewidth}{!}{%
\begin{tikzpicture}[
  x=2cm,y=2cm,
  font=\small,
  line cap=round,line join=round,
  cell edge/.style={draw=black!80,line width=.65pt},
  cut/.style={
    draw=black!55,line width=.5pt,
    dash pattern=on 3pt off 2.5pt
  },
  upper brace/.style={
    decorate,decoration={brace,amplitude=3.5pt},
    line width=.5pt
  },
  lower brace/.style={
    decorate,decoration={brace,mirror,amplitude=3.5pt},
    line width=.5pt
  },
  guide/.style={draw=black!27,line width=.35pt},
  every node/.style={inner sep=1pt}
]

\foreach \a/\b/\j/\col in {
  0/.48/1/slabBlue,
  .48/1.28/2/slabOchre,
  1.28/3.35/3/slabTeal,
  3.35/3.65/4/slabBlue,
  3.65/4.55/5/slabOchre}{
  \fill[\col!13] (\a,0) rectangle (\b,1);
  \draw[upper brace,draw=\col]
    ({\a+.025},1.10)--({\b-.025},1.10)
    node[midway,above=5pt,text=\col] {$t_{\j}$};
}

\foreach \a in {.48,1.28,3.35,3.65,4.55}
  \draw[cut] (\a,0)--(\a,1);

\draw[cell edge]
  (0,0)--(4.70,0)
  (0,1)--(4.70,1);

\foreach \a in {0,1,2,3,4}
  \draw[cell edge] (\a,0)--(\a,1);

\foreach \a/\b/\i/\j/\col in {
  0/.48/1/1/slabBlue,
  .48/1/1/2/slabOchre,
  1/1.28/2/2/slabOchre,
  1.28/2/2/3/slabTeal,
  2/3/3/3/slabTeal,
  3/3.35/4/3/slabTeal,
  3.35/3.65/4/4/slabBlue,
  3.65/4/4/5/slabOchre,
  4/4.55/5/5/slabOchre}{
  \draw[lower brace,draw=\col]
    ({\a+.018},-.09)--({\b-.018},-.09)
    node[midway,below=5pt,text=\col] {$S_{\j}^{\i}$};
}

\foreach \a in {0,1,2,3,4}
  \draw[guide] (\a,1.50)--(\a,1.80);

\foreach \i in {1,2,3,4}{
  \draw[upper brace,draw=black!70]
    ({\i-1+.025},1.83)--({\i-.025},1.83)
    node[midway,above=5pt,text=black!85] {$C_{\i}$};
}

\node[text=black!65,font=\normalsize]
  at (5.22,.50) {$\cdots$};

\begin{scope}[shift={(6,0)}]

  \fill[black!4] (-.22,0) rectangle (.30,1);

  \foreach \a/\b/\j/\col in {
    .30/1.45/{K-3}/slabBlue,
    1.45/1.82/{K-2}/slabOchre,
    1.82/2.38/{K-1}/slabTeal,
    2.38/3/K/slabBlue}{
    \fill[\col!13] (\a,0) rectangle (\b,1);
    \draw[upper brace,draw=\col]
      ({\a+.025},1.10)--({\b-.025},1.10)
      node[midway,above=5pt,text=\col] {$t_{\j}$};
  }

  \foreach \a in {.30,1.45,1.82,2.38}
    \draw[cut] (\a,0)--(\a,1);

  \draw[cell edge]
    (-.22,0)--(3,0)
    (-.22,1)--(3,1);

  \foreach \a in {0,1,2,3}
    \draw[cell edge] (\a,0)--(\a,1);

  \foreach \a/\b/\i/\j/\col in {
    .30/1/{N-2}/{K-3}/slabBlue,
    1/1.45/{N-1}/{K-3}/slabBlue,
    1.45/1.82/{N-1}/{K-2}/slabOchre,
    2/2.38/N/{K-1}/slabTeal,
    2.38/3/N/K/slabBlue}{
    \draw[lower brace,draw=\col]
      ({\a+.018},-.09)--({\b-.018},-.09)
      node[midway,below=5pt,text=\col] {$S_{\j}^{\i}$};
  }

  \draw[lower brace,draw=slabTeal]
    (1.838,-.09)--(1.982,-.09);
  \draw[draw=slabTeal!75,line width=.4pt]
    (1.91,-.16)--(1.91,-.67);
  \node[below,text=slabTeal]
    at (1.91,-.69) {$S_{K-1}^{N-1}$};

  \foreach \a in {0,1,2,3}
    \draw[guide] (\a,1.50)--(\a,1.80);

  \foreach \a/\b/\i in {
    0/1/{N-2},
    1/2/{N-1},
    2/3/N}{
    \draw[upper brace,draw=black!70]
      ({\a+.025},1.83)--({\b-.025},1.83)
      node[midway,above=5pt,text=black!85] {$C_{\i}$};
  }

\end{scope}
\end{tikzpicture}%
}
\caption[Subdivision of the cells into the pieces $S_\ell^i$.
The dependence on $z$ is suppressed for simplicity.]{%
Subdivision of the cells into the pieces $S_\ell^i$.
The dependence on $z$ is suppressed for simplicity.%
\protect\footnotemark}
\label{fig:join-subdivision}
\end{figure}

For $z\in\widetilde Z_K$, define the incremental slabs
$$
    S_\ell(z):=[s_{\ell-1},s_\ell)\times[0,1]^{d-1},
    \qquad s_0=0,
    \qquad s_\ell=\sum_{j=1}^{\ell}t_j,
    \qquad\ell=1,\dots,K.
$$
All partitions are understood up to null sets. If $t_\ell=0$,
then $S_\ell$ is empty, consistently with the identifications
in the join. We further set (see figure~\ref{fig:join-subdivision})
$$
    S_\ell^i(z):=S_\ell(z)\cap C_i,
    \qquad i=1,\dots,N.
$$

\footnotetext{This figure was converted from a hand-drawn
sketch to TikZ code using ChatGPT.}

Now on each nonempty slab $S_\ell^i$, we use the maps of
Proposition~\ref{prop:cube-to-slab-map} and assign sets. First consider $\phi_{\eta}$ to be the set map of shifting in the $e_1$ direction by $\eta$. Further consider the map:
$$
    T_{\ell}^i := \phi_{\max\{i-1,s_{\ell-1}\}}
\circ T_{|S_\ell^i(z)|}
\circ\phi_{1-i} : \mathcal C(C_i) \to \mathcal C(S_{\ell}^i).
$$
Here $T_{\nu}$ is the map defined in Proposition~\ref{prop:cube-to-slab-map} and $|S_\ell^i(z)|$ denotes its volume, equivalently its
width in the first coordinate (here empty slabs contribute
the empty set). Then assign the set:
\[
E_\ell^i(z)
:=
T_\ell^i
\bigl(\widetilde\Phi_i(\widetilde x_i^\ell)\bigr),
\]
Pasting these sets of finite perimeter together gives a set
$\widetilde\Psi(z)\in\mathcal C(Q)$, precisely defined as:
$$
    \widetilde\Psi(z) := \sum_{i=1}^{N}\sum_{\ell=1}^K E_\ell^i(z).
$$
By the equivariance in Proposition~\ref{prop:cube-to-slab-map}(2) The construction respects complementation:
\begin{align*}
\begin{aligned}
\widetilde\Psi(\tau^{(K)}z)
&=\sum_{i=1}^N\sum_{\ell=1}^K
T_\ell^i\bigl(\widetilde\Phi_i(\tau_i\widetilde x_i^\ell)\bigr)=\sum_{i=1}^N\sum_{\ell=1}^K
T_\ell^i\bigl(C_i\setminus\widetilde\Phi_i(\widetilde x_i^\ell)\bigr)\\
&=\sum_{i=1}^N\sum_{\ell=1}^K
\bigl(S_\ell^i\setminus E_\ell^i(z)\bigr)
=Q\setminus\widetilde\Psi(z),
\end{aligned}
\end{align*}
where in the second-to-last equality we used Proposition~\ref{prop:cube-to-slab-map}(2).
Consequently, its relative boundary defines a map on the quotient,
$$
\Psi:Z_K\longrightarrow \mathcal Z^0_{d-1,\mathrm{rel}}(Q,\de Q;\mathbb Z_2),\qquad
\Psi([z]):=\partial_{\mathrm{rel}}[\widetilde\Psi(z)].
$$
\begin{lemma}\label{lem:join-construction-is-higher-sweepout}
    The map $\Psi$ constructed above is indeed a $K(p+1)-1$ sweepout, continuous in flat topology with no concentration of mass.
\end{lemma}
\begin{proof}
By Proposition~\ref{prop:cube-to-slab-map}(1), continuity
of translations in $L^1$, and continuity of the lifted
filling maps, each $E_\ell^i(z)$ depends continuously on
the parameters whenever $S_\ell^i(z)$ has positive volume.
Continuity also holds when this slab collapses, since
\[
 |E_\ell^i(z)|\leq |S_\ell^i(z)|\longrightarrow0.
\]
Moreover, as $t_\ell\to0$,
\[
 \sum_{i=1}^N|E_\ell^i(z)|
 \leq\sum_{i=1}^N|S_\ell^i(z)|
 =t_\ell\longrightarrow0,
\]
uniformly in the coordinates forgotten by the join.
Consequently, $\widetilde\Psi$ is well defined and
continuous in $L^1$ on $\widetilde Z_K$.
The estimate
\[
 \flat\bigl(
 \partial_{\mathrm{rel}}[E],
 \partial_{\mathrm{rel}}[E']\bigr)
 \leq |E\triangle E'|
\]
and the equivariance of $\widetilde\Psi$ then show that
$\Psi$ is well defined and flat continuous on $Z_K$.

The family $\Psi$ also has no concentration of mass.
Indeed, each $\Phi_i$ is mass continuous and therefore
has no concentration of mass by
Proposition~\ref{prop:mass-continuity-no-concentration}.
Applying Proposition~\ref{prop:cube-to-slab-map}(4) to
the translated lifted fillings
$\phi_{1-i}(\widetilde\Phi_i(\widetilde x_i))$,
and using translation invariance of the supremum over
ball centers, gives
\[
 \lim_{r\downarrow0}
 \sup_{\substack{z\in\widetilde Z_K\\y\in\R^d}}
 \|\partial E_\ell^i(z)\|(B_r(y))=0
 \qquad\text{for every }i,\ell.
\]
Since the pieces have disjoint interiors, subadditivity
of boundary mass yields
\[
 \sup_{\substack{[z]\in Z_K\\y\in Q}}
 \|\Psi([z])\|(B_r(y))
 \leq
 \sum_{i=1}^N\sum_{\ell=1}^K
 \sup_{\substack{z\in\widetilde Z_K\\y\in\R^d}}
 \|\partial E_\ell^i(z)\|(B_r(y))
 \longrightarrow0.
\]
Here the full boundaries of the pieces already include
all contributions along the gluing interfaces.

We next verify the sweepout condition. All cohomology
groups below have coefficients in $\Z_2$.
The characteristic class of the filling double cover
is $\bar\lambda$: a loop has nontrivial monodromy precisely
when its lifted fillings differ at the endpoints by the
whole domain, equivalently when its Almgren class is the
relative fundamental class. This identification follows from Almgren's isomorphism
\cite{Almgren-homotopy}; see \cite[Section~2.5]{Weyl-LMN}
for the relative formulation and
\cite[Section~5]{MarquesNeves2021}
for the filling double-cover description in codimension one.
Thus, writing $w_1$ for the first Stiefel-Whitney class
of a double cover, we have
\[
 w_1(\pi_i)^p
 =(\Phi_i^*\bar\lambda)^p\ne0
 \qquad\text{for every }i.
\]

For a finite polyhedron with a free $\Z_2$-action, its
mod-$2$ cohomological height is the largest exponent for
which the corresponding power of $w_1$ is nonzero.
The product formula in \cite[Corollary~2.5]{DKV2023} states that the height of a diagonal product is the minimum
of the heights of its factors. Hence
$\widetilde X_1\times\cdots\times\widetilde X_N$
has height at least $p$.

The join formula in \cite[Theorem~2.4]{DKV2023} states
that the height of a join is the sum of the heights
of its factors plus one. This additivity already appears,
in the projective-space formulation, in
\cite{Givental1990}; see also
\cite[Section~3 and Appendix~B]{GKPS2021}.
Iterating the join formula gives
\[
 w_1(\widetilde Z_K\longrightarrow Z_K)^{K(p+1)-1}
 \ne0.
\]

Finally, the equivariant filling map $\widetilde\Psi$
identifies $\widetilde Z_K\to Z_K$ with the pullback
under $\Psi$ of the filling double cover of the relative
cycle space of $Q$. Naturality therefore gives
\[
 \Psi^*\bar\lambda
 =w_1(\widetilde Z_K\longrightarrow Z_K),
\]
and consequently
\[
 \Psi^*(\bar\lambda^{K(p+1)-1})
 =
 w_1(\widetilde Z_K\longrightarrow Z_K)^{K(p+1)-1}
 \ne0.
\]
Since $Z_K$ is a finite polyhedron, it admits a finite
cubical structure, so this is a $K(p+1)-1$ sweepout
in the sense of our definition.
\end{proof}
Define the integrand $H$ as follows:
$$
    H\big|_{\operatorname{int}C_i} = H_i.
$$
On cube faces we set $H=1$; their contributions
will be bounded separately.
With the definition above we state the following lemma:
\begin{lemma}\label{lem:join-construction-area-calculation}
    The following estimate is true:
    \begin{align*}
    \begin{aligned}
        \sup_{x\in\mathrm{dmn}(\Psi)} \M_H(\Psi(x)) \leq  p^{1/d}&(K+N-1)^{\frac1d}\left(\sum_{i=1}^N a(H_i)^{\frac{d}{d-1}}\right)^{\frac{d-1}{d}} \\ &+ o(p^{1/d}) + C(\lambda,N,d)K + C(N,K)\epsilon.
    \end{aligned}
    \end{align*}
\end{lemma}
\begin{proof}

    First note that, bounding the boundaries of the cubes with some constant $C(N,d,\lambda)$ we get the estimate:
$$
        \sup_{x\in\mathrm{dmn}(\Psi)} \M_H(\Psi(x)) \leq C(N,d,\lambda) +  \sup_{x\in\mathrm{dmn}(\Psi)} \sum_{i=1}^{N}\sum_{\ell=1}^{K}\M_{H_i}(\pa E_\ell^i(x)).
$$
    Here the full boundaries of the pieces include all interfaces
    between adjacent slabs.

    Consider the width of each $S_{\ell}^i(x)$ to be:
$$
        w_{\ell}^i(x) = |S_{\ell}^i(x)|,\qquad \sum_{i=1}^N\sum_{\ell=1}^K w_{\ell}^i = N.
$$

    Then we have, by the relative-perimeter estimate proved in
    Proposition~\ref{prop:cube-to-slab-map}, adding the boundary
    of each target slab, and almost optimality of the sweepouts
    $\Phi_i$ in \eqref{eq:almost-optimality-sweepout}:
$$
\sum_{i=1}^{N}\sum_{\ell=1}^{K}\M_{H_i}(\pa E_\ell^i(x))\leq \sum_{i=1}^{N}\sum_{\ell=1}^{K}(w_\ell^i)^{\frac{d-1}{d}}\omega_p^{d-1}(C_i,H_i)+C(\lambda,d)NK+NK\epsilon.
$$
Empty slabs contribute zero. Now for $p$ large enough, we have:
$$
\omega_p^{d-1}(C_i,H_i) = p^{1/d}a(H_i) + o(p^{1/d}).
$$
    Hence we see that:

    \begin{align*}
    \begin{aligned}
        \sum_{i=1}^{N}\sum_{\ell=1}^{K}\M_{H_i}(\pa E_\ell^i(x))\leq\, &p^{1/d}\sum_{i=1}^{N}\sum_{\ell=1}^{K}(w_\ell^i)^{\frac{d-1}{d}}a(H_i) \\ +&  C(\lambda,d)NK+C(N,K) \left(\epsilon + o(p^{1/d})\right).
    \end{aligned}
    \end{align*}
    Inside each cube $i$ we define $m_i$ to be the number of slabs
    $S_\ell$ whose intersection with $C_i$ has positive volume:
$$
m_i := |\{1\leq \ell \leq K\,:\, |S_\ell\cap C_i|>0\}|=|\{1\leq \ell \leq K\,:\,w_\ell^i>0\}|.
$$
    Since these pieces partition $C_i$ modulo null sets,
    $\sum_{\ell=1}^K w_\ell^i=1$. Then using H\"older inequality in each cube and then with all the slabs, we have:
$$
    \sum_{i=1}^{N}\sum_{\ell=1}^{K}(w_\ell^i)^{\frac{d-1}{d}}a(H_i)\leq \sum_{i=1}^{N} m_i^{\frac1d} a(H_i) \leq \left(\sum_{i=1}^{N} m_i\right)^{1/d}\left(\sum_{i=1}^N a(H_i)^{\frac{d}{d-1}}\right)^{\frac{d-1}{d}}.
$$
    This is indeed an optimal integer allocation problem. Note that the total number of positive-volume slabs $S_\ell^i$ is at most $K+N-1$, since the total number of slabs is at most the number of cuts plus one (counting the cuts the cubes already have), which is indeed $K+N-1$. Then we have:
$$
\sum_{i=1}^{N}\sum_{\ell=1}^{K}(w_\ell^i)^{\frac{d-1}{d}}a(H_i) \leq (K+N-1)^{\frac1d}\left(\sum_{i=1}^N a(H_i)^{\frac{d}{d-1}}\right)^{\frac{d-1}{d}}.
$$
    Putting everything together, we finish the proof.
\end{proof}

\subsection{The upper bound}
In this section we prove the general upper bound for manifolds with variable anisotropies. To do this we subdivide the manifold into small cubes and map them into a cube in a long rectangle. Then we use the constructions of the last section and build admissible sweepouts with correct energy bounds.
\begin{thm}
\label{thm:variable-weyl-upper}
Let $(M^d,g)$ be a smooth, compact, connected manifold
with smooth boundary (possibly empty),
$d\geq2$, and let $F$ satisfy
(\ref{H:comparability}-\ref{H:ellipticity}). Then
\[
\limsup_{p\to\infty}
p^{-1/d}\omega_p^{d-1}(M,F)
\leq
\left(
\int_M a(F(x,\cdot))^{\frac{d}{d-1}}\,dV_g(x)
\right)^{\frac{d-1}{d}}.
\]
\end{thm}

\begin{proof}
We first assume that $\partial M=\varnothing$.
Fix $\epsilon>0$. Partition $M$, up to boundaries, into
finitely many piecewise smooth regions whose closures
lie in small normal-coordinate charts. Cover their
coordinate images by grid cubes, retaining the whole
cubes meeting each region. For sufficiently fine grids,
all retained cubes lie in their coordinate domains,
and the total excess volume tends to zero as the grids
are refined.

Enumerate the retained coordinate cubes as
\[
B_i=b_i+\ell_i[0,1]^d,
\qquad i=1,\ldots,N,
\]
where $b_i\in\mathbb R^d$ and $\ell_i>0$. For each $i$,
let
\[
\varphi_i:U_i\subset\mathbb R^d\longrightarrow M
\]
be the normal-coordinate parametrization used for
$B_i$, so that $B_i\subset U_i$. The same chart may
occur for several indices.

We define
\[
C_i=[i-1,i]\times[0,1]^{d-1},
\qquad
Q:=\bigcup_{i=1}^N C_i
=[0,N]\times[0,1]^{d-1}.
\]
Thus $C_i$ is the normalized copy of $B_i$, placed in
the $i$-th position of the long rectangle $Q$.
The corresponding map back to $M$ is
\[
f_i:C_i\longrightarrow M,
\qquad
f_i(y_1,\ldots,y_d)
:=
\varphi_i\!\left(
b_i+\ell_i(y_1-i+1,y_2,\ldots,y_d)
\right).
\]
Indeed, this formula first translates $C_i$ onto
$[0,1]^d$, then rescales it onto $B_i$, and finally
applies the coordinate parametrization $\varphi_i$.
In particular,
\[
f_i(C_i)=\varphi_i(B_i).
\]

We set
\[
x_i=f_i\!\left(i-\tfrac12,\tfrac12,\ldots,\tfrac12\right),
\qquad
H_i=F(x_i,\cdot),
\qquad
v_i=\vol_g(f_i(C_i)),
\]
using the orthonormal frame obtained from the coordinate
frame at $x_i$ to represent $H_i$. Here $v_i$ records
the volume of the whole coordinate cube in $M$,
whereas each normalized cube $C_i$ has unit Euclidean
volume; see Figure~\ref{fig:manifold-cubes-to-rectangle}.

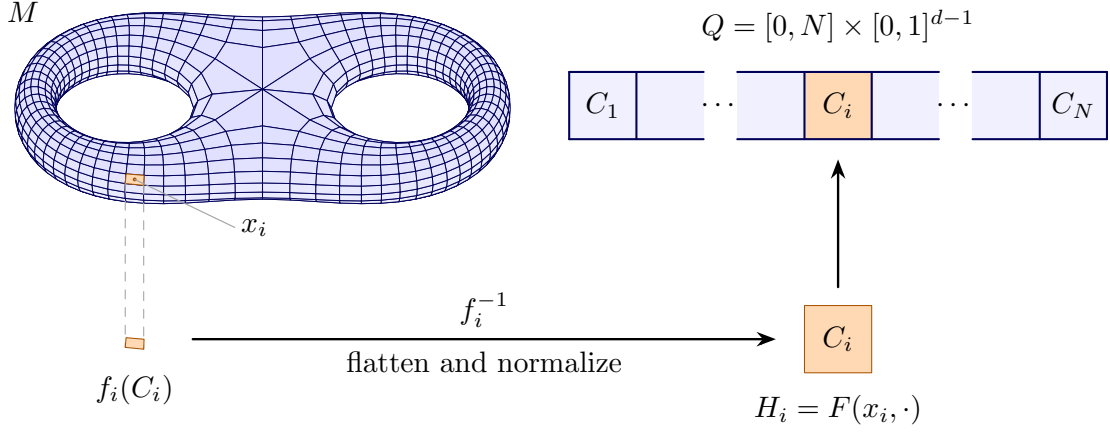
\begin{figure}[H]
\centering
\resizebox{0.98\linewidth}{!}{%
\begin{tikzpicture}[
    line cap=round,line join=round,font=\small,
    outline/.style={draw=blue!35!black,line width=0.7pt},
    selected/.style={draw=orange!65!black,fill=orange!30,
                     line width=0.35pt},
    arrow/.style={-{Stealth[length=2.3mm]},thick},
    guide/.style={gray!65,dashed,thin}
]

\pgfmathdeclarefunction{manifoldX}{2}{%
  \pgfmathparse{sqrt(max(0,
    (sqrt(1+(1.2+.4*cos(#2))^2
      +2*(1.2+.4*cos(#2))*cos(#1))
      +1+(1.2+.4*cos(#2))*cos(#1))/2))}%
}
\pgfmathdeclarefunction{manifoldY}{2}{%
  \pgfmathparse{sign(#1)*sqrt(max(0,
    (sqrt(1+(1.2+.4*cos(#2))^2
      +2*(1.2+.4*cos(#2))*cos(#1))
      -1-(1.2+.4*cos(#2))*cos(#1))/2))}%
}

\begin{axis}[
    hide axis,view={0}{45},
    at={(3.05cm,3.40cm)},anchor=origin,
    x={(1.83cm,0cm)},
    y={(0cm,1.29cm)},z={(0cm,1.29cm)},
    colormap={surfaceblue}{
      color(0cm)=(blue!5);color(1cm)=(blue!12)},
    shader=faceted,faceted color=blue!35!black,
    z buffer=sort,
    every axis plot/.append style={line width=0.16pt}
]
    \addplot3[surf,domain=-180:180,y domain=0:360,
              samples=37,samples y=19]
      ({manifoldX(x,y)},
       {manifoldY(x,y)},
       {0.20*sin(y)});

    \addplot3[surf,domain=-180:180,y domain=0:360,
              samples=37,samples y=19]
      ({-manifoldX(x,y)},
       {-manifoldY(x,y)},
       {0.20*sin(y)});

    \coordinate (A) at
      (axis cs:{-manifoldX(120,20)},
               {-manifoldY(120,20)},{0.20*sin(20)});
    \coordinate (B) at
      (axis cs:{-manifoldX(130,20)},
               {-manifoldY(130,20)},{0.20*sin(20)});
    \coordinate (C) at
      (axis cs:{-manifoldX(130,40)},
               {-manifoldY(130,40)},{0.20*sin(40)});
    \coordinate (D) at
      (axis cs:{-manifoldX(120,40)},
               {-manifoldY(120,40)},{0.20*sin(40)});
\end{axis}

\node at (0.20,4.55) {$M$};

\path[selected] (A)--(B)--(C)--(D)--cycle;
\coordinate (xi) at ($(A)!0.5!(C)$);
\fill[orange!65!black] (xi) circle (0.65pt);
\draw[gray!70,thin] (xi) -- (2.72,1.98);
\node[right,inner sep=2pt] at (2.72,1.98) {$x_i$};

\foreach \P in {A,B,C,D}{
    \coordinate (\P-copy) at ($(\P)+(0,-1.95)$);
}
\draw[guide] (A) -- (A-copy);
\draw[guide] (B) -- (B-copy);
\path[selected]
  (A-copy)--(B-copy)--(C-copy)--(D-copy)--cycle;
\node[below=5pt] at ($(A-copy)!0.5!(B-copy)$)
  {$f_i(C_i)$};

\draw[arrow] (2.22,0.65) -- (9.16,0.65)
  node[midway,above] {$f_i^{-1}$}
  node[midway,below] {flatten and normalize};

\path[selected] (9.50,0.25) rectangle ++(0.80,0.80);
\node at (9.90,0.65) {$C_i$};
\node[below] at (9.90,0.12) {$H_i=F(x_i,\cdot)$};

\begin{scope}[shift={(6.70,3.03)}]
    \fill[blue!6] (0,0) rectangle (1.60,0.80);
    \fill[blue!6] (2.00,0) rectangle (4.40,0.80);
    \fill[blue!6] (4.80,0) rectangle (6.40,0.80);

    \path[selected] (2.80,0) rectangle (3.60,0.80);

    \foreach \a/\b in {0/1.60,2.00/4.40,4.80/6.40}{
        \draw[outline] (\a,0) -- (\b,0);
        \draw[outline] (\a,0.80) -- (\b,0.80);
    }
    \foreach \x in {0,0.80,2.80,3.60,5.60,6.40}{
        \draw[outline] (\x,0) -- (\x,0.80);
    }

    \node at (1.80,0.40) {$\cdots$};
    \node at (4.60,0.40) {$\cdots$};
    \node at (0.40,0.40) {$C_1$};
    \node at (3.20,0.40) {$C_i$};
    \node at (6.00,0.40) {$C_N$};

    \node[above] at (3.20,1.02)
      {$Q=[0,N]\times[0,1]^{d-1}$};
\end{scope}

\draw[arrow] (9.90,1.27) -- (9.90,2.79);

\end{tikzpicture}%
}
\caption[The coordinate cubes are flattened, normalized,
and arranged as $Q$. The highlighted cube carries the
frozen anisotropy $H_i=F(x_i,\cdot)$.]{%
The coordinate cubes are flattened, normalized,
and arranged as $Q$. The highlighted cube carries the
frozen anisotropy $H_i=F(x_i,\cdot)$.\protect\footnotemark}
\label{fig:manifold-cubes-to-rectangle}
\end{figure}
\footnotetext{This figure was converted from a hand-drawn
sketch to TikZ code using ChatGPT.}

By choosing the charts sufficiently small and then
refining the grids, normal-coordinate estimates and
continuity of $F$ give
\[
\F((f_i)_\#T)
\leq
(1+\epsilon)v_i^{\frac{d-1}{d}}\M_{H_i}(T)
\]
for every rectifiable $(d-1)$-chain $T$ in $C_i$.
Moreover, continuity of the frozen coefficient and
the vanishing excess volume give
\[
\sum_{i=1}^N v_i a(H_i)^{\frac{d}{d-1}}
\leq
\int_M a(F(x,\cdot))^{\frac{d}{d-1}}\,dV_g(x)
+\epsilon.
\]

Fix this cover and $K\geq1$. Apply the preceding join
construction on $Q$ to the frozen integrands
$v_i^{(d-1)/d}H_i$. By homogeneity,
\[
a\bigl(v_i^{\frac{d-1}{d}}H_i\bigr)^{\frac{d}{d-1}}
=
v_i a(H_i)^{\frac{d}{d-1}}.
\]
Transfer the resulting fillings to $M$ by setting
\[
D_i:=f_i(C_i)\setminus\bigcup_{j<i}f_j(C_j),
\qquad
\widetilde\Psi_M(z)
:=
\bigcup_{i=1}^N
\left(
f_i(\widetilde\Psi(z)\cap C_i)\cap D_i
\right).
\]
Since the $D_i$ partition $M$ up to null sets, this map
is $L^1$ continuous and satisfies
\[
\widetilde\Psi_M(\tau^{(K)}z)
=
M\setminus\widetilde\Psi_M(z).
\]
Consequently,
$\Psi_M([z]):=\partial[\widetilde\Psi_M(z)]$
is flat continuous and pulls back the same filling
double cover as $\Psi$. Hence it is a
$K(p+1)-1$ sweepout by
Lemma~\ref{lem:join-construction-is-higher-sweepout}.

The additional boundaries are supported on the fixed
interfaces $\bigcup_i\partial D_i$, whose total area
bounds their mass independently of $p$ and $K$.
Their mass measures have no atoms, so these interfaces
introduce no concentration of mass. The same property
for the remaining pieces follows from that of $\Psi$
and the bi-Lipschitz coordinate maps.

We may therefore apply
Proposition~\ref{prop:mass-approximation} on $M$.
The freezing estimate and
Lemma~\ref{lem:join-construction-area-calculation}
then yield
\[
\begin{aligned}
\omega_{K(p+1)-1}^{d-1}(M,F)
&\leq
(1+\epsilon)p^{1/d}(K+N-1)^{1/d}
\left(
\sum_{i=1}^N v_i a(H_i)^{\frac{d}{d-1}}
\right)^{\frac{d-1}{d}}
\\
&\qquad +o(p^{1/d})+C_{\epsilon,K},
\end{aligned}
\]
where the cover and $K$ are fixed as $p\to\infty$.

Divide by $(K(p+1)-1)^{1/d}$ and let $p\to\infty$.
Monotonicity of the widths extends the resulting
estimate to all indices, since
\[
\frac{K(p+2)-1}{K(p+1)-1}
=
1+\frac{K}{K(p+1)-1}
\longrightarrow1.
\]
Thus
\[
\limsup_{p\to\infty}
p^{-1/d}\omega_p^{d-1}(M,F)
\leq
(1+\epsilon)
\left(1+\frac{N-1}{K}\right)^{1/d}
\left(
\sum_{i=1}^N v_i a(H_i)^{\frac{d}{d-1}}
\right)^{\frac{d-1}{d}}.
\]
Letting $K\to\infty$ with the cover fixed, and using
the weighted Riemann-sum estimate above, gives
\[
\limsup_{p\to\infty}
p^{-1/d}\omega_p^{d-1}(M,F)
\leq
(1+\epsilon)
\left(
\int_M a(F(x,\cdot))^{\frac{d}{d-1}}\,dV_g(x)
+\epsilon
\right)^{\frac{d-1}{d}}.
\]
Finally, let $\epsilon\downarrow0$.

For a nonempty smooth boundary, embed $M$ as a compact domain
in a smooth closed connected manifold $\widehat M$, extending
$g$ smoothly. Given $\eta>0$, extend $F$ continuously to a thin
exterior collar, preserving positivity and convexity of the
normal extension. Such an extension is obtained by transporting
the normal integrands along the collar through linear
isomorphisms of the tangent spaces. Blend this extension with
a small positive constant multiple of the area integrand
outside the collar, using a partition of unity depending only
on the base point. Convexity and positivity are preserved.
The collar can be chosen sufficiently thin and the constant
sufficiently small that the resulting integrand $\widehat F$
agrees with $F$ on $M$ and satisfies
\[
\int_{\widehat M\setminus M}
 a(\widehat F(x,\cdot))^{\frac d{d-1}}\,dV_{\widehat g}(x)
<\eta.
\]
Here boundedness on the collar follows from comparability,
and homogeneity of $a$ controls the remaining region.

To compare widths, lift any mass-continuous $p$-sweepout on
$\widehat M$ to its filling double cover and restrict each
filling $E$ to $E\cap M$. This restriction is $L^1$ continuous
and commutes with complementation in the respective domains.
Its relative boundary therefore pulls back the same covering
class, and its anisotropic mass is at most that of the original
cycle. Indeed, the new trace term lies on $\partial M$ and is
ignored in the relative mass. The relative boundary measures
inside $M$ are restrictions of the original boundary measures,
so no concentration of mass is preserved. Applying
Proposition~\ref{prop:mass-approximation} on $M$ and then taking
the infimum over the original sweepouts yields
\[
\omega_p^{d-1}(M,F)
\leq\omega_p^{d-1}(\widehat M,\widehat F).
\]
The closed case just proved gives
\[
\limsup_{p\to\infty}p^{-1/d}\omega_p^{d-1}(M,F)
\leq
\left(
\int_M a(F(x,\cdot))^{\frac d{d-1}}\,dV_g(x)+\eta
\right)^{\frac{d-1}{d}}.
\]
Letting $\eta\downarrow0$ proves the result with boundary.
\end{proof}

At this point we are ready to prove the full variable Weyl law combining the lower and upper bound:
\begin{proof}[Proof of Theorem~\ref{thm:main-weyl}]
    The proof follows immediately from combining the lower bound of Corollary~\ref{cor:variable-weyl-lower} and the matching upper bound in Theorem~\ref{thm:variable-weyl-upper}.
\end{proof}

\begin{rmk}
    Using the existence of the anisotropic Weyl Law in Theorem~\ref{thm:main-weyl} and the homogeneity $a(cH)=c\,a(H)$ for all $c>0$ and frozen anisotropies $H$, one may also adapt the arguments of \cite{IMN2018} to anisotropies with generic background metrics, to show the density of anisotropic stationary hypersurfaces in the base manifold.
\end{rmk}

\part{Grassmannian equidistribution}\label{part:equidistribution}

\section{The smooth realization assumption and the conclusions}
\label{sec:equidistribution-scope}

Throughout this part, $(M^d,g)$ is closed, connected, and smooth,
with $3\leq d\leq7$. We will keep the metric $g$ fixed and we will instead perturb the
anisotropic integrand. We write
$G:=G_{d-1}(M)$, for the bundle of unoriented tangent hyperplanes.  For a varifold $V$, we write
\[
 \F(V):=\F(V,M)=\int_G F(x,T)\,dV(x,T).
\]

To prove equidistribution in the Grassmannian bundle, we
combine the Weyl law with smooth realization of every width
for integrands sufficiently close to area.
We assume throughout
Part~\ref{part:equidistribution} this realization property in Assumption \ref{hyp:near-area-realization}.
Theorem~\ref{thm:main-realization} proves Assumption \ref{hyp:near-area-realization} for
$3\leq d \leq 6$. 

\begin{ipotesi}\label{hyp:near-area-realization}
There exist $\delta_0>0$ and an integer $k_0\geq6$, with the following property. For every smooth integrand
$F:G\to(0,\infty)$ satisfying $\|F-1\|_{C^{k_0}(G)}<\delta_0$
and every $p\geq1$, there exists an integral varifold
\[
 V=\sum_{i=1}^s m_i|\Sigma_i|,
 \qquad m_i\in\N,
\]
such that $\Sigma_i$ are smooth, closed, embedded, connected,
pairwise disjoint $F$-stationary hypersurfaces and
\[
 \F(V)=\omega_p^{d-1}(M,F).
\]
\end{ipotesi}
We can decrease $\delta_0$ so that the normal
representations of the integrands in this neighborhood of the area are uniformly strictly elliptic. 

The following theorem is the main result of Part~\ref{part:equidistribution}. 
\begin{thm}
\label{thm:fixed-metric-anisotropic-equidistribution}
Assume that $(M^d,g)$ satisfies
Assumption~\ref{hyp:near-area-realization}.
There exist a sequence of smooth positive strictly elliptic integrands
$F_j\to1$ in $C^\infty(G)$ and a sequence of smooth, closed, embedded,
connected $F_j$-stationary hypersurfaces $\{\Sigma_j\}_j$ such that
\begin{equation}\label{eq:fixed-g-anisotropic-sequence}
 \lim_{q\to\infty}\frac{\displaystyle\sum_{j=1}^q\int_{\Sigma_j}h(x,T_x\Sigma_j)\,dA_g}{\displaystyle\sum_{j=1}^q\operatorname{area}_g(\Sigma_j)}
 =
 \int_G h\,d\Lambda_g
 \qquad\text{for every }h\in C^0(G).
\end{equation}
Each $\Sigma_j$ can be chosen nondegenerate for $F_j$ and as
a component of a smooth $F_j$-stationary varifold realizing
$\omega_{p_j}^{d-1}(M,F_j)$ for some $p_j\geq1$.
\end{thm}

Here nondegeneracy means that the Jacobi operator has trivial
kernel on smooth sections of the normal bundle. Combining Theorem~\ref{thm:fixed-metric-anisotropic-equidistribution} with Theorem~\ref{thm:main-realization}, we deduce Theorem \ref{thm:main-equidistribution}.

The proof follows the finite-dimensional deformation and
averaging argument of Marques-Neves-Song~\cite{MNS},
with the anisotropy as the parameter. We first use rotation
invariance and homogeneity of the frozen Weyl constant to
identify the limit of averaged derivatives of the Weyl
functional. This step does not require differentiability
of $a$ at $1$ or an explicit formula for $a$. We then adapt
White's transversality argument to variations of $F$ with
$g$ fixed. Together with
Assumption~\ref{hyp:near-area-realization}, this gives a
formula for derivatives of widths in terms of smooth
stationary hypersurfaces. We obtain finite weighted
collections approximating any prescribed finite family
of moments of $\Lambda_g$, and then assemble these
collections into the sequence in
Theorem~\ref{thm:fixed-metric-anisotropic-equidistribution}.

\section{Averaged derivatives of the Weyl functional}
\label{sec:grassmann-symmetry}

We work in the neighborhood of area from
Assumption~\ref{hyp:near-area-realization}. We use the Weyl functional
$\mathcal W$ and the probability measure $\Lambda_g$ introduced in
\eqref{eq:intro-anisotropic-weyl} and
\eqref{eq:intro-grassmann-measure}, and set
\begin{equation}\label{eq:grassmann-target-functional}
 L_g(h):=\mathcal W(1)\int_G h\,d\Lambda_g.
\end{equation}
Our first goal is to approximate $L_g$ by averaged derivatives of
$\mathcal W$ at integrands approaching area. We will not assume differentiability of $a$ at $1$, and we will rely instead on its rotation invariance and homogeneity.

\subsection{Finite-dimensional families of test functions}
\label{subsec:grassmann-test-families}

For $m\geq0$, let $\mathscr E_m\to M$ be the vector bundle of even
spherical harmonics of degrees at most $2m$ on the unit spheres
of $TM$, equipped with the fiberwise $L^2(\sigma_x)$ inner product.
We identify its sections with functions on $G$.
The union, over $m\geq0$, of the spaces of smooth sections of $\mathscr E_m$ is uniformly dense in $C(G)$.
Indeed, one can first approximate a continuous function by a smooth one, then apply uniform polynomial approximation in local orthonormal frames and combine the approximations with a partition of unity.
Taking the even parts gives a finite harmonic expansion in each fiber.

The following lemma provides finite families whose evaluation maps
respect the fiberwise inner products. It allows us to average over
the fibers without choosing a global orthonormal frame of $TM$.

\begin{lemma}\label{lem:grassmann-parseval}
Let $\mathscr E\to M$ be a smooth Euclidean vector bundle
of rank $r$, and let $h_0,\ldots,h_s$ be smooth sections.
There exist smooth sections $s_1,\ldots,s_L$ such that
\[
\begin{gathered}
 \sum_{a=1}^L s_a(x)\otimes s_a(x)
 =\operatorname{Id}_{\mathscr E_x}
 \;\text{for every }x\in M,
\\
\mbox{ and }
\quad
 h_0,\ldots,h_s
 \in\operatorname{span}_{\R}\{s_1,\ldots,s_L\}.
\end{gathered}
\]
In particular, the maps
\[
 A_x:\R^L\longrightarrow\mathscr E_x,
 \qquad
 A_xt:=\sum_{a=1}^L t_as_a(x),
\]
satisfy $A_xA_x^*=\operatorname{Id}_{\mathscr E_x}$.
\end{lemma}

\begin{proof}
Choose a finite open cover $U_1,\ldots,U_N$ of $M$
such that $\mathscr E$ admits a smooth orthonormal frame
$e_{i,1},\ldots,e_{i,r}$ on each $U_i$.
Choose smooth functions $\psi_1,\ldots,\psi_N$ on $M$
such that
\[
 \operatorname{spt}\psi_i\subset U_i
 \quad\text{for }i=1,\ldots,N,
 \qquad
 \sum_{i=1}^N\psi_i^2=1.
\]

Since $M$ is compact, we can choose $\eta>0$ sufficiently
small that
\[
 B:=\operatorname{Id}_{\mathscr E}
     -\eta^2\sum_{j=0}^s h_j\otimes h_j
\]
is positive definite. Its positive square root $B^{1/2}$
is a smooth bundle endomorphism.
We take $s_1,\ldots,s_L$, with $L=Nr+s+1$, to be the sections
\[
 \psi_i B^{1/2}e_{i,\alpha}
 \quad (1\leq i\leq N,\ 1\leq\alpha\leq r),
 \qquad
 \eta h_0,\ldots,\eta h_s,
\]
where the sections supported in $U_i$ are extended by zero
outside $U_i$.

By orthonormality of the local frames,
\[
 \sum_{a=1}^L s_a\otimes s_a
 =\sum_{i=1}^N\psi_i^2 B
   +\eta^2\sum_{j=0}^s h_j\otimes h_j
 =\operatorname{Id}_{\mathscr E}.
\]
Moreover, each $h_j$ is a constant multiple of one of the
chosen sections. Finally,
\[
 A_xA_x^*
 =\sum_{a=1}^L s_a(x)\otimes s_a(x)
 =\operatorname{Id}_{\mathscr E_x},
\]
which completes the proof.
\end{proof}

\subsection{Averaging the Weyl derivatives}
\label{subsec:averaged-weyl-derivatives}
From Lemma~\ref{lem:frozen-weyl-lipschitz} we know that $a$ is locally Lipschitz.
The same argument gives local Lipschitz continuity
of $\mathcal W$. We now identify the limit of its
averaged derivatives near $1$ (the area functional).

\begin{proposition}\label{prop:global-averaged-weyl-derivative}
Let $h_0=1,h_1,\ldots,h_s$ be smooth sections of $\mathscr E_m$.
There exist a finite-dimensional subspace $U\subset C^\infty(G)$
containing these sections and a surjective linear map
$A:\R^L\to U$ such that
\[
 w(v):=\mathcal W(1+Av)
\]
is well-defined and Lipschitz in a neighborhood of $0$.
Moreover, for every nonnegative smooth radial function
$\rho$ compactly supported in the unit ball of $\R^L$ with
$\int_{\R^L}\rho(t)\,dt=1$, we have
\begin{equation}\label{eq:averaged-weyl-derivative}
 \lim_{\eps\downarrow0}
 \int_{\R^L}\rho(t)\,Dw(\eps t)[b]\,dt
 =L_g(Ab)
 \qquad\text{for every }b\in\R^L.
\end{equation}
Consequently, $L_g|_U$ belongs to the closed convex hull
of the limits of derivatives of $\mathcal W$ along $1+U$
at differentiability points approaching $1$.
The corresponding integrands converge to $1$ in $C^\infty(G)$.
\end{proposition}

\begin{proof}
Let $s_1,\ldots,s_L$ be the sections given by
Lemma~\ref{lem:grassmann-parseval}, and set
\[
 U:=\operatorname{span}_{\R}\{s_1,\ldots,s_L\},
 \qquad
 At:=\sum_{a=1}^L t_as_a.
\]
Then $U$ contains $h_0,\ldots,h_s$, and $A$ is surjective.
Since $A$ has smooth image and finite-dimensional domain,
$1+Av$ belongs to the prescribed neighborhood of area
for $v$ sufficiently small. The local Lipschitz bound
for $\mathcal W$ gives the corresponding bound for $w$.

The evaluation maps $A_xt=(At)(x,\cdot)$ satisfy
$A_xA_x^*=\operatorname{Id}_{\mathscr E_{m,x}}$.
Since $\rho$ is radial, the pushforward of $\rho(t)\,dt$
under $A_x$ has a smooth density which is the same in every
fiber after orthogonal identification.
Indeed, in orthonormal coordinates $A_x$ is an orthogonal
projection, and this density is obtained by integrating
$\rho$ over the remaining coordinates.
In particular, the projected measure is invariant under
the action induced by $O(d)$.

We evaluate $a$ in orthonormal frames, as in the definition
of $\mathcal W$, and take its derivatives within the
finite-dimensional fibers $\mathscr E_{m,x}$.
By the local Lipschitz bound, $a$ is differentiable almost
everywhere in each fiber near $1$.
Since $A_x$ is surjective, Fubini's theorem implies that,
for almost every $t$ in the support of $\rho$,
$a$ is differentiable at $1+\eps A_xt$ for almost every $x$.
The uniform Lipschitz bound allows differentiation under
the integral. Applying the chain rule, we obtain
\begin{equation}\label{eq:integrated-fiber-derivative}
 Dw(\eps t)[b]
 =
 \int_M
 \left(
 \frac{a(1+\eps A_xt)}{w(\eps t)}
 \right)^{1/(d-1)}
 Da(1+\eps A_xt)[A_xb]\,dV_g(x)
\end{equation}
for almost every $t$.

By rotation invariance of $a$ and of the projected measure,
the linear functional
\[
 h\longmapsto
 \int_{\R^L}\rho(t)\,Da(1+\eps A_xt)[h]\,dt
 \qquad\text{on }\mathscr E_{m,x}
\]
is invariant under $O(d)$.
Averaging $h$ over $O(d)$ gives its mean times the constant
function. Hence
\[
\begin{aligned}
 &\int_{\R^L}\rho(t)\,Da(1+\eps A_xt)[h]\,dt
 \\
 &\quad =
 \left(\int_{G_{d-1}(T_xM)}h\,d\sigma_x\right)
 \int_{\R^L}\rho(t)\,Da(1+\eps A_xt)[1]\,dt.
\end{aligned}
\]
On the other hand, homogeneity gives
\[
 Da(H)[1]=a(H)-Da(H)[H-1]
\]
at every differentiability point.
Since $Da$ is locally bounded and $a(1)=a_d$, it follows that
\begin{equation}\label{eq:fiber-average-derivative}
 \lim_{\eps\downarrow0}
 \int_{\R^L}\rho(t)\,Da(1+\eps A_xt)[h]\,dt
 =
 a_d\int_{G_{d-1}(T_xM)}h\,d\sigma_x.
\end{equation}
The convergence is uniform in $x$ and for $h$ in the
unit balls of the fibers.

Moreover, $w(0)=a_d\vol_g(M)^{(d-1)/d}$, so
\[
 \left(
 \frac{a(1+\eps A_xt)}{w(\eps t)}
 \right)^{1/(d-1)}
 \longrightarrow \vol_g(M)^{-1/d}
\]
uniformly for $x\in M$ and $t\in\operatorname{supp}\rho$.
Using the uniform bound for $Da$, we may therefore replace
this factor in the average of
\eqref{eq:integrated-fiber-derivative}
by $\vol_g(M)^{-1/d}$ with an error tending to zero.
Equation~\eqref{eq:fiber-average-derivative} then gives
\begin{align}
\begin{aligned}
 &\lim_{\eps\downarrow0}
 \int_{\R^L}\rho(t)\,Dw(\eps t)[b]\,dt
 \\ &\quad =
 a_d\vol_g(M)^{-1/d}
 \int_M\int_{G_{d-1}(T_xM)}
 (Ab)(x,T)\,d\sigma_x(T)\,dV_g(x)
 =
 L_g(Ab).
\end{aligned}
\end{align}
This proves \eqref{eq:averaged-weyl-derivative}.

To prove the convex hull assertion, recall that $Dw$ is
bounded near $0$.
If $L_g\circ A$ were outside the closed convex hull of
the limits of $Dw(v)$ as $v\to0$ through differentiability
points, separation and boundedness would give
$b\in\R^L$ and $\delta>0$ such that
\[
 Dw(v)[b]\leq L_g(Ab)-\delta
\]
at every differentiability point sufficiently close to $0$.
Integrating this inequality with $v=\eps t$ against
$\rho(t)\,dt$ contradicts
\eqref{eq:averaged-weyl-derivative}.

Since $A$ is surjective, differentiability of $w$ at $v$
is equivalent to differentiability of
$\mathcal W|_{1+U}$ at $1+Av$, and the derivatives are
related by composition with $A$.
Thus the preceding conclusion gives the asserted convex hull property for the derivatives of
$\mathcal W$ restricted to $1+U$.
Finally, for every $k$,
\[
 \|Av\|_{C^k}
 \leq \|A\|_{\R^L\to C^k}|v|
 \longrightarrow0
 \qquad\text{as }v\to0,
\]
which proves the last assertion.
\end{proof}

\subsection{Passing to the volume spectrum}
\label{sec:grassmann-spectrum-derivatives}

For the map $A$ given by
Proposition~\ref{prop:global-averaged-weyl-derivative}, set
\[
 u_p(v):=\widetilde\omega_p^{d-1}(M,1+Av).
\]
The Weyl law gives $u_p\to w$ pointwise.
Comparison of the integrands and the uniform bound for
$\widetilde\omega_p^{d-1}(M,1)$ give a Lipschitz bound for
$u_p$, independent of $p$, on a sufficiently small parameter
ball. Hence the convergence is uniform on compact subsets.
The following lemma allows us to pass to the limit in
averages of the derivatives.

\begin{lemma}\label{lem:averaged-spectral-derivatives}
Let $A$ and $\rho$ be as in
Proposition~\ref{prop:global-averaged-weyl-derivative}.
For every sufficiently small fixed $\eps>0$, we have
\begin{equation}\label{eq:spectral-derivative-convergence}
 \lim_{p\to\infty}
 \int_{\R^L}\rho(t)\,Du_p(\eps t)[b]\,dt
 =
 \int_{\R^L}\rho(t)\,Dw(\eps t)[b]\,dt
 \qquad\text{for every }b\in\R^L.
\end{equation}
For each $u_p$ and for $w$, the corresponding averaged
derivative can be approximated, uniformly for $|b|\leq1$,
by a convex combination of at most $L+1$ derivatives
evaluated at points of $\eps\operatorname{supp}\rho$.
These points may be chosen in any prescribed subset
of full Lebesgue measure in the parameter ball on which
the corresponding function is differentiable.
\end{lemma}

\begin{proof}
Since $u_p$ is Lipschitz and $\rho$ is smooth and compactly
supported, integration by parts gives
\[
 \int_{\R^L}\rho(t)\,Du_p(\eps t)[b]\,dt
 =
 -\frac1\eps
 \int_{\R^L}u_p(\eps t)\,D\rho(t)[b]\,dt.
\]
The same identity holds for $w$. Subtracting the two
identities, we bound the difference of the averaged
derivatives by
\[
 \frac{|b|}{\eps}
 \sup_{t\in\operatorname{supp}\rho}
 |u_p(\eps t)-w(\eps t)|
 \int_{\R^L}|D\rho(t)|\,dt.
\]
Since $\eps$ is fixed, this tends to zero uniformly for
$|b|\leq1$, proving
\eqref{eq:spectral-derivative-convergence}.

Finally, each averaged derivative belongs to the closed
convex hull of the derivatives evaluated at points of
$\eps\operatorname{supp}\rho$ in the prescribed full-measure
set. It can therefore be approximated uniformly for
$|b|\leq1$ by a finite convex combination of such derivatives.
Since the space of linear functionals on $\R^L$ has
dimension $L$, Carath\'eodory's theorem reduces each
combination to at most $L+1$ terms.
\end{proof}

\section{Geometric realization and equidistribution}
\label{sec:grassmann-geometric-realization}

We now use the smooth realizations in
Assumption~\ref{hyp:near-area-realization}
and deformations of the integrand to represent the
averaged spectral derivatives by stationary hypersurfaces.

\subsection{Deformations of stationary hypersurfaces}
\label{sec:anisotropic-white}

We adapt White's deformation argument~\cite{White1991}
by varying the integrand while keeping the metric fixed.
All perturbations remain in the neighborhood specified in
Assumption~\ref{hyp:near-area-realization}.
We first record the variation formulas that we will use.

\begin{lemma}\label{lem:anisotropic-first-variation-white}
Let $F$ be a smooth strictly elliptic integrand, and let
$\Sigma\subset M$ be a smooth, closed, embedded hypersurface.
For a smooth normal variation, write its velocity locally
as $u\nu$, where $\nu$ is a unit normal. Then
\begin{equation}\label{eq:anisotropic-euler-white}
 \begin{aligned}
 \delta\F(|\Sigma|)[u\nu]=\int_\Sigma uH_F\,dA_g,
 \qquad H_F=\operatorname{div}_\Sigma D_\xi\varphi_x(\nu)+\langle\nabla_x\varphi_x(\nu),\nu\rangle.
 \end{aligned}
\end{equation}
Here $\nabla_x$ denotes differentiation in the base point
while parallel transporting $\xi$ with respect to the
Levi-Civita connection, and $\operatorname{div}_\Sigma$
is the tangential trace of the ambient covariant derivative.

If $\Sigma$ is $F$-stationary, the linearization $J_F$ of
$H_F$ under normal variations has the local expression
\begin{equation}\label{eq:anisotropic-jacobi-white}
 J_Fu
 =-\operatorname{div}_\Sigma(A_F\nabla_\Sigma u)+c_Fu,
 \qquad
 A_F(x):=
 D^2_{\xi\xi}\varphi_x(\nu)|_{T_x\Sigma},
\end{equation}
where $c_F$ is a smooth function on $\Sigma$.
The operator is elliptic and formally self-adjoint with
respect to $dA_g$. For every $0<\alpha<1$, it is Fredholm
of index zero from $C^{2,\alpha}$ sections of the normal
bundle to $C^{0,\alpha}$ sections.
For $F=1$, it reduces to
\[
 J_Fu
 =-\Delta_\Sigma u
  -\bigl(|A_\Sigma|^2+\operatorname{Ric}_g(\nu,\nu)\bigr)u.
\]

Finally, let $f\in C^\infty(M)$ and set
$F_s(x,T)=e^{sf(x)}F(x,T)$.
On an $F$-stationary hypersurface,
\begin{equation}\label{eq:anisotropic-forcing-white}
 \left.\frac{d}{ds}\right|_{s=0}H_{F_s}
 =\langle D_\xi\varphi_x(\nu),\nabla f\rangle.
\end{equation}
If $f|_\Sigma=0$, the right-hand side is
$F(x,T_x\Sigma)\partial_\nu f$.
All these formulas also apply to one-sided hypersurfaces,
with the scalar expressions interpreted using local
unit normals.
\end{lemma}

\begin{proof}
Set $H=\operatorname{div}_\Sigma\nu$.
Under a normal variation with velocity $u\nu$, the
derivatives of the normal and the area element are
$-\nabla_\Sigma u$ and $uH\,dA_g$, respectively.
Thus
\[
 \delta\F(|\Sigma|)[u\nu]
 =
 \int_\Sigma
 \left(
 u\langle\nabla_x\varphi_x(\nu),\nu\rangle
 -\langle D_\xi\varphi_x(\nu),\nabla_\Sigma u\rangle
 +uH\varphi_x(\nu)
 \right)\,dA_g.
\]
By homogeneity, $\langle D_\xi\varphi_x(\nu),\nu\rangle
 =\varphi_x(\nu)$.
Integrating the tangential part of
$D_\xi\varphi_x(\nu)$ by parts therefore proves
\eqref{eq:anisotropic-euler-white}.

The second-order terms in the normal linearization of
$H_F$ agree with those of the operator
$-\operatorname{div}_\Sigma(A_F\nabla_\Sigma u)$.
At an $F$-stationary hypersurface, symmetry of the second
variation makes $J_F$ formally self-adjoint with respect
to $dA_g$. The divergence-form operator is also formally
self-adjoint, so their difference has no first-order term.
This proves \eqref{eq:anisotropic-jacobi-white}.
Strict ellipticity of $F$ gives ellipticity of $J_F$.
Since $\Sigma$ is closed, elliptic theory gives the
Fredholm property, and self-adjointness gives index zero.
When $F=1$, we have $H_F=H$, and the stated expression
for $J_F$ is the usual normal variation formula for
mean curvature.

The normal representation of $F_s$ is
$e^{sf(x)}\varphi_x(\xi)$.
Substituting this into
\eqref{eq:anisotropic-euler-white} gives
\[
 H_{F_s}
 =
 e^{sf}
 \left(
 H_F+s\langle D_\xi\varphi_x(\nu),\nabla f\rangle
 \right).
\]
Since $H_F=0$, differentiation at $s=0$ proves
\eqref{eq:anisotropic-forcing-white}.
If $f|_\Sigma=0$, then
$\nabla f=(\partial_\nu f)\nu$ along $\Sigma$,
and homogeneity gives
\[
 \langle D_\xi\varphi_x(\nu),\nabla f\rangle
 =\varphi_x(\nu)\partial_\nu f
 =F(x,T_x\Sigma)\partial_\nu f.
\]

Finally, $\varphi_x$ is even. Replacing $\nu$ by $-\nu$
changes the signs of $u$, $H_F$, and $J_Fu$.
Hence $uH_F$ and the normal fields $H_F\nu$ and
$(J_Fu)\nu$ are independent of the local choice of
unit normal. This proves that the formulas also apply
when $\Sigma$ is one-sided.
\end{proof}

\begin{proposition}\label{prop:anisotropic-white-structure}
Fix an integer $q\geq k_0$ and $0<\alpha<1$.
Consider pairs $(F,\Sigma)$, where $F$ is a $C^q$ integrand
in the strictly elliptic neighborhood of area fixed above and
$\Sigma$ is a closed embedded $C^{2,\alpha}$
$F$-stationary hypersurface.
The space of these pairs is a separable $C^{q-3}$ Banach
manifold admitting a countable atlas, and the projection
\[
 (F,\Sigma)\longmapsto F
\]
is Fredholm of index zero.
Its critical points are precisely the pairs for which $J_F$ has a nontrivial kernel on $C^{2,\alpha}$ sections
of the normal bundle.

Within the same neighborhood, the smooth integrands for which
every closed embedded stationary hypersurface is nondegenerate
form a residual set in the $C^\infty$ topology.
Moreover, every smooth family $t\mapsto F(t)$, parametrized
by an open subset of $\R^b$, can be perturbed arbitrarily
slightly in the $C^\infty$ topology to a family transverse
to the projection, taking
$q\geq\max\{k_0,b+5\}$.
For such a family, every closed embedded
$F(t)$-stationary hypersurface is nondegenerate
for almost every $t$.

All assertions hold for both two-sided and one-sided
hypersurfaces, with $J_F$ acting on sections of the normal bundle.
\end{proposition}

\begin{proof}
We use the local construction in the proof of
\cite[Theorem~2.1]{White1991}, with the integrand as parameter.
The variation formulas and Fredholm properties in
Lemma~\ref{lem:anisotropic-first-variation-white}
remain valid for $C^q$ integrands and $C^{2,\alpha}$
stationary hypersurfaces, by the same local computations
and elliptic regularity.

Choose a fixed smooth reference hypersurface and represent
nearby hypersurfaces as normal graphs of $C^{2,\alpha}$ sections.
After identifying the normal bundles and pulling back to the reference hypersurface, the stationarity equation defines a $C^{q-3}$ map into the fixed space of $C^{0,\alpha}$ normal sections. At a stationary graph, its derivative in the graph variable is identified with $J_F$: the graph velocity is first projected onto the normal bundle of the current hypersurface, and the resulting mean-curvature variation is pulled back to the reference hypersurface. The terms obtained by differentiating these identifications vanish because $H_F=0$. We first show that variations of the integrand complement
the range of this operator.

By self-adjointness, it suffices to show that no nonzero
$w\in\ker J_F$ is orthogonal to all variations of $H_F$
induced by changing $F$.
Choose a patch where a unit normal is defined and $w$
is nonzero, and represent $w$ by its scalar coefficient
with respect to this normal.
Let $\chi$ be a nonnegative smooth cutoff supported in
the patch, with $\chi w\not\equiv0$.
In a tubular neighborhood of the patch, choose an ambient
function $f$ satisfying
$f|_\Sigma=0$, and $\partial_\nu f=\chi w$, and extend it by zero away from that neighborhood.
By \eqref{eq:anisotropic-forcing-white},
\[
 \int_\Sigma
 w\left.\frac{d}{ds}\right|_{s=0}H_{e^{sf}F}\,dA_g
 =
 \int_\Sigma \chi F(x,T_x\Sigma)w^2\,dA_g
 >0.
\]
The scalar expressions in this identity are supported in
the chosen patch.
If $\Sigma$ has only finite regularity, we approximate $f$
in $C^1$ by smooth ambient functions supported in the same
neighborhood. The strict inequality persists, since the
variation of $H_F$ depends continuously on $f$ in $C^1$.

Thus the integrand variations span the cokernel of $J_F$.
Since this cokernel is finite-dimensional, finitely many
smooth variations suffice.
The full linearization is therefore surjective with
complemented kernel, and the implicit function theorem
gives the local solution manifold.
The kernel of the derivative of the projection is
$\ker J_F$, while its range has codimension
$\dim\ker J_F$.
This proves the index and critical point assertions.
Since the preceding perturbations are supported in a
single patch, the argument also applies when the normal
bundle is nontrivial.

We next justify the existence of a countable atlas for the space of pairs $(F,\Sigma)$.
For each diffeomorphism type, choose a countable family
of smooth reference embeddings into $M$ that is dense
in the $C^3$ topology.
Elliptic regularity gives $C^{3,\alpha}$ local graph
representations for every stationary hypersurface in the space. Each such hypersurface can therefore be
represented as a normal graph over a sufficiently close
reference hypersurface from one of these countable families.

In each reference chart, the stationary graph sections
belong to the closure of the smooth sections in
$C^{2,\alpha}$.
This closure is a separable metric space, as is $C^q(G)$.
The solution set in the chart is therefore second countable
in the induced topology.
Since there are only countably many diffeomorphism types
and reference embeddings, these charts give a countable
open cover of the solution space.
Consequently, the solution manifold is second countable
and admits a countable Banach manifold atlas.

We next prove the assertions about smooth families.
Fix the parameter dimension $b$ and take
$q\geq\max\{k_0,b+5\}$.
At a fixed solution, split the normal sections into
$\ker J_F$ and its orthogonal complement.
Solving the equation in the complement by the implicit
function theorem reduces stationarity locally to
\[
 K(F,z)=0,
 \qquad z\in\R^r,
 \qquad K(F,z)\in\R^r,
\]
where $r$ is the dimension of the Jacobi kernel at the
fixed solution.
The map $K$ is $C^{q-3}$, and its full derivative is
surjective at every zero in a sufficiently small neighborhood.
A family $t\mapsto F(t)$ is transverse to the projection
precisely when
$D_{(t,z)}\bigl(K(F(t),z)\bigr):\R^b\times\R^r\rightarrow\R^r$
is surjective at every zero.

We verify this condition on compact subsets of the
$t$ and $z$ domains, restricting $F(t)$ to a closed
$C^q$ ball contained strictly inside the integrand domain
of the chart.
For any fixed family, the corresponding set of zeros
in $(t,z)$ is compact.
Surjectivity at these zeros is therefore open under
$C^1$ perturbations of the family with values in $C^q$.

To prove density, first choose slightly larger restrictions
of the same kind.
The preceding cokernel argument and compactness give
finitely many smooth functions $f_1,\ldots,f_N$ such that,
for
\[
 F_{\mathbf a}(t,x,T)
 =
 \exp\left(\sum_{i=1}^N a_i f_i(x)\right)F(t,x,T),
 \qquad
 \mathbf a=(a_1,\ldots,a_N)\in\R^N,
\]
the derivative of $K(F_{\mathbf a}(t),z)$ with respect to
$(\mathbf a,t,z)$ is surjective near the relevant zeros.
The larger restrictions ensure that every zero in the
smaller region remains covered after a sufficiently small
perturbation.
A cutoff in $t$, equal to one on the larger compact
parameter set, allows these perturbations to be supported
in a compact subset of the parameter domain.

The resulting zero set has dimension $b+N$.
Since $q-3>b$, Sard's theorem applied to its projection
onto $\R^N$ shows that arbitrarily small choices of
$\mathbf a$ give a transverse family on the restricted
region.
This is the local parametric transversality argument
of \cite{Smale}.
Because all $f_i$ are smooth, the perturbation can be made
arbitrarily small in any prescribed finite collection
of $C^\infty$ seminorms.
It also preserves strict ellipticity and the fixed
$C^{k_0}$ neighborhood.

We equip the space of smooth families taking values in
the prescribed neighborhood with the $C^\infty$ topology
on compact subsets of the parameter domain.
This is a Baire space.
Indeed, a countable compact exhaustion of the parameter
domain expresses the condition of taking values in the
prescribed neighborhood as a countable intersection
of open conditions in the Fr\'echet space of smooth
functions on the product of the parameter domain with $G$.
The resulting space is therefore completely metrizable.

Using the countable atlas, we choose countably many
nested restrictions of the preceding kind whose interiors
cover all solutions.
Each localized transversality condition is open and dense
in this space of smooth families, by the preceding
perturbation argument.
Baire's theorem therefore gives a residual, and hence
dense, set of families transverse in every chart.
Taking $b=0$ proves the residual assertion for smooth
integrands.

Finally, for a transverse $b$-parameter family, the space
of pairs $(t,\Sigma)$ with $\Sigma$ $F(t)$-stationary is
a $b$-dimensional $C^{q-3}$ manifold with a countable atlas.
Sard's theorem implies that almost every $t$ is a regular
value of the projection
$(t,\Sigma)\longmapsto t$.
The kernel of its derivative at $(t,\Sigma)$ is precisely
$\ker J_{F(t)}$.
Since the domain and target have the same dimension,
regularity forces this kernel to vanish.
Thus every closed embedded stationary hypersurface is
nondegenerate at almost every parameter.
\end{proof}

\subsection{Derivatives of the widths}
\label{subsec:countable-spectral-branches}

We now combine Assumption~\ref{hyp:near-area-realization}
with Proposition~\ref{prop:anisotropic-white-structure}
to express the derivatives of the widths in terms of
their realizing hypersurfaces.
We first record an elementary lemma.

\begin{lemma}\label{lem:countable-smooth-selection}
Let $\Omega\subset\R^N$ be open, and let
$w:\Omega\to\R$ be locally Lipschitz.
Suppose that there are countably many open sets
$U_\ell\subset\Omega$ and functions $e_\ell\in C^1(U_\ell)$
such that, for almost every $t\in\Omega$,
\[
 w(t)=e_\ell(t)
 \qquad\text{for some }\ell\text{ with }t\in U_\ell.
\]
Then, for almost every $t\in\Omega$, such an $\ell$
can be chosen so that
$Dw(t)=De_\ell(t)$.
\end{lemma}

\begin{proof}
For each $\ell$, set $E_\ell:=\{t\in U_\ell:w(t)=e_\ell(t)\}$.
Almost every point of $E_\ell$ is a density point
of $E_\ell$ at which $w$ is differentiable.
At any such point $t$, we have
$D(w-e_\ell)(t)=0$.
Indeed, a nonzero differential would imply that
$w-e_\ell$ is nonzero in a cone with vertex $t$,
contradicting the density-one property of $E_\ell$.
Since the collection is countable and
$\bigcup_\ell E_\ell$ has full measure in $\Omega$,
the conclusion follows.
\end{proof}

\begin{proposition}\label{prop:anisotropic-countable-derivative}
Suppose that Assumption~\ref{hyp:near-area-realization}
holds, and let $t\mapsto F(t)$ be a smooth family of smooth
integrands in the realization neighborhood of area, defined near
a compact cube in $\R^N$.
After an arbitrarily small perturbation in the $C^\infty$
topology, there is a full-measure set of parameters in
the cube with the following property.

For every $p\geq1$, the function
$t\mapsto\omega_p^{d-1}(M,F(t))$ is differentiable at $t$,
and there is a realizing integral varifold $V_{p,t}$
with smooth nondegenerate connected components such that
\begin{equation}\label{eq:anisotropic-spectral-derivative}
 D_t\omega_p^{d-1}(M,F(t))[v]
 =
 \int_G D_tF(t)[v](x,T)\,dV_{p,t}(x,T)
 \qquad\text{for every }v\in\R^N.
\end{equation}
The full-measure set is independent of $p$.
\end{proposition}

\begin{proof}
Choose $q\geq\max\{k_0,N+5\}$ and apply
Proposition~\ref{prop:anisotropic-white-structure}
to perturb the family to one transverse to the projection.
We keep the perturbation inside the neighborhood in
Assumption~\ref{hyp:near-area-realization}.
For the resulting family, all closed embedded
$F(t)$-stationary hypersurfaces are nondegenerate
for almost every parameter $t$.

By Proposition~\ref{prop:anisotropic-white-structure},
the derivative of the projection onto the integrand space
is an isomorphism at every nondegenerate stationary
hypersurface. The inverse function theorem therefore
gives a neighborhood of each such pair on which the
projection is a diffeomorphism onto an open set of
integrands.
Since the solution manifold admits a countable atlas,
these neighborhoods have a countable subcover.
Composing the corresponding local inverses with
$t\mapsto F(t)$ gives countably many $C^1$ branches,
defined on open parameter sets, which cover all
nondegenerate connected stationary hypersurfaces
in the family.

A varifold realizing a width may have several connected
components, each with a positive integer multiplicity.
We therefore consider all finite collections of these
branches, assigning a positive integer multiplicity
to each member. Each collection is defined on the
intersection of the parameter domains of its branches,
restricted to the open set where the corresponding
hypersurfaces are pairwise disjoint.
Adding their energies with the chosen multiplicities
gives a $C^1$ function.
There are countably many such choices, so the resulting
energy functions form a countable collection,
which we denote by $e_\ell$.
At every parameter where all stationary hypersurfaces
are nondegenerate, each realization supplied by
Assumption~\ref{hyp:near-area-realization}
is represented by one of these collections.

Fix $p\geq1$.
Comparison of the integrands shows that
$t\mapsto\omega_p^{d-1}(M,F(t))$ is locally Lipschitz.
By the preceding paragraph, it agrees almost everywhere
with at least one of the functions $e_\ell$.
Lemma~\ref{lem:countable-smooth-selection} therefore gives,
at almost every parameter, a realizing collection whose
energy has the same differential as the width.
Let $V_{p,t}$ be the varifold determined by this collection
at $t$.

When differentiating its energy, the contribution from
the motion of each hypersurface vanishes by stationarity.
Only the derivative of the integrand remains, giving
\eqref{eq:anisotropic-spectral-derivative}.
Since the lemma identifies the full differentials,
the same varifold $V_{p,t}$ gives this identity for every
$v\in\R^N$.

Finally, intersect these full-measure sets over
$p\geq1$ and with the full-measure set where all stationary
hypersurfaces are nondegenerate.
This gives all the conclusions on a single full-measure
set of parameters.
\end{proof}

\subsection{Finite stationary collections}
\label{sec:geometric-averaged-derivatives}

We apply Proposition~\ref{prop:anisotropic-countable-derivative}
to a small perturbation of the affine family.
The following lemma controls the resulting error in the
averaged derivatives.

\begin{lemma}\label{lem:generic-perturbation-averaging}
Fix $A$ and $\rho$ as above, $k\geq k_0$, and a compact
parameter cube containing $\operatorname{supp}\rho$.
Let $\eps>0$ be sufficiently small that
\[
 H(t):=1+\eps At
\]
belongs to the neighborhood in
Assumption~\ref{hyp:near-area-realization}
on a neighborhood of this cube.
For every sufficiently small $\eta>0$, there exists a
smooth family $\widetilde H$ satisfying the conclusions
of Proposition~\ref{prop:anisotropic-countable-derivative},
with
\[
 \|\widetilde H-H\|_{C^1_tC^k_G}<\eta.
\]
The realizing varifolds $V_{p,t}$ given by that proposition
can be chosen so that
$
 t\mapsto\int_G h\,dV_{p,t}
$
is measurable for every $h\in C(G)$ and every $p\geq1$.
Moreover,
\begin{equation}\label{eq:generic-averaging-error}
 \left|
 \int\rho(t)
 \left(
 p^{-1/d}\int_G Ab\,dV_{p,t}
 -Du_p(\eps t)[b]
 \right)\,dt
 \right|
 \leq C\frac{\eta}{\eps}|b|
 \qquad\text{for every }b\in\R^L,
\end{equation}
where $C$ is independent of $p$ and $\eta$.
\end{lemma}

\begin{proof}
Apply Proposition~\ref{prop:anisotropic-countable-derivative}
to choose $\widetilde H$ with the required smallness.
We first explain how to choose the realizing varifolds
measurably.

For each $p$, enumerate the finite collections of stationary
branches constructed in the proof of that proposition.
At each parameter in the full-measure set given there,
choose the first collection whose energy and differential
agree with those of the width.
The sets on which these identities hold are measurable:
the branch energies are $C^1$, the width is locally
Lipschitz, and equality of differentials can be checked
in the finitely many coordinate directions.
For each $h\in C(G)$, its integral over any fixed branch
collection varies continuously with the parameter.
Thus the selected varifolds have the asserted measurability.
On the exceptional null set, set $V_{p,t}=0$.

Write
\[
 v_p(t):=
 \widetilde\omega_p^{d-1}(M,\widetilde H(t)).
\]
Since the integrands lie between fixed positive constants,
comparison of the widths and the isotropic upper bound give
\[
 p^{-1/d}\|V_{p,t}\|(M)\leq C,
 \qquad
 \sup_t|v_p(t)-u_p(\eps t)|\leq C\eta.
\]
The constants are independent of $p$ and $\eta$.

By \eqref{eq:anisotropic-spectral-derivative},
\[
 Dv_p(t)[b]
 =
 p^{-1/d}
 \int_G D_t\widetilde H(t)[b]\,dV_{p,t}
\]
for almost every $t$ and every $b\in\R^L$.
Since $D_tH(t)[b]=\eps Ab$, the assumed smallness and
the mass bound imply
\[
 \left|
 \frac1\eps Dv_p(t)[b]
 -p^{-1/d}\int_G Ab\,dV_{p,t}
 \right|
 \leq C\frac{\eta}{\eps}|b|.
\]

On the other hand, integration by parts gives
\begin{align*}
&\frac1\eps\int\rho(t)Dv_p(t)[b]\,dt
 -\int\rho(t)Du_p(\eps t)[b]\,dt
 \\ &\qquad =-\frac1\eps
 \int\bigl(v_p(t)-u_p(\eps t)\bigr)D\rho(t)[b]\,dt.
\end{align*}
Consequently,
\[
 \left|
 \frac1\eps\int\rho(t)Dv_p(t)[b]\,dt
 -\int\rho(t)Du_p(\eps t)[b]\,dt
 \right|
 \leq
 C\frac{\eta}{\eps}|b|
 \int|D\rho(t)|\,dt.
\]
Combining the last two estimates and using
$\rho\geq0$ and $\int\rho=1$ proves
\eqref{eq:generic-averaging-error}.
\end{proof}

\begin{proposition}\label{prop:geometric-grassmann-collections}
Suppose that Assumption~\ref{hyp:near-area-realization} holds.
There exist integers $p_m\uparrow\infty$ and, for every $m\geq1$,
finite collections of smooth anisotropies $F_{m,j}$ in the
realization neighborhood, integral varifolds $V_{m,j}$,
and weights $\lambda_{m,j}\geq0$, indexed by $1\leq j\leq N_m$,
with
\[
 \sum_{j=1}^{N_m}\lambda_{m,j}=1,
\]
such that each $V_{m,j}$ realizes
$\omega_{p_m}^{d-1}(M,F_{m,j})$ and its connected components
are smooth, closed, embedded, $F_{m,j}$-stationary and nondegenerate.
Moreover,
\[
 \max_{1\leq j\leq N_m}\|F_{m,j}-1\|_{C^k(G)}
 \longrightarrow0
 \qquad\text{for every }k,
\]
and the measures
\[
\mu_m:=p_m^{-1/d}\sum_{j=1}^{N_m}\lambda_{m,j}V_{m,j}
\]
satisfy
$\mu_m\rightharpoonup\mathcal W(1)\Lambda_g$ on $G$.
In particular,
\[
\mu_m(G)\longrightarrow\mathcal W(1)>0,
 \qquad
 \frac{\mu_m}{\mu_m(G)}\rightharpoonup\Lambda_g.
\]
\end{proposition}

\begin{proof}
Choose a uniformly dense sequence
$h_0=1,h_1,h_2,\ldots$ in $C(G)$ consisting of smooth
functions with finite fiberwise harmonic expansions.

Fix $m\geq1$.
Apply Proposition~\ref{prop:global-averaged-weyl-derivative}
to $h_0,\ldots,h_m$, taking a harmonic bundle containing
all these functions.
Write $A:\R^L\to C^\infty(G)$ for the resulting map,
and choose $b_0,\ldots,b_m\in\R^L$ such that
$Ab_\ell=h_\ell$ for $0\leq\ell\leq m$.
Choose a smooth radial probability density $\rho$
supported in the unit ball of $\R^L$, and fix a compact
parameter cube containing its support in its interior.
All these choices, as well as the corresponding functions
$w$ and $u_p$, depend on $m$; we suppress this dependence
in the notation.
All suprema below are taken over the chosen cube.

First choose $\eps>0$ sufficiently small that
\[
 \left|
 \int_{\R^L}\rho(t)\,Dw(\eps t)[b_\ell]\,dt
 -L_g(h_\ell)
 \right|
 <\frac1m
 \qquad\text{for }0\leq\ell\leq m,
\]
and
\[
 \sup_t\|\eps At\|_{C^{\max\{k_0,m\}}(G)}
 <\min\left\{\frac{\delta_0}{2},\frac1{2m}\right\}.
\]
By Lemma~\ref{lem:generic-perturbation-averaging}, we can
choose a smooth perturbation $\widetilde H$ of
$t\mapsto1+\eps At$ such that
\[
 \sup_t\|\widetilde H(t)-1\|_{C^{\max\{k_0,m\}}(G)}
 <\min\left\{\delta_0,\frac1m\right\},
\]
and the error in \eqref{eq:generic-averaging-error}
is less than $1/m$ in each of the directions
$b_0,\ldots,b_m$.
Since that estimate is uniform in $p$, the same
perturbation gives this bound for every $p\geq1$.

Keeping $\eps$ and $\widetilde H$ fixed, choose $p_m$
sufficiently large, with $p_m>p_{m-1}$ when $m>1$,
that the error in
\eqref{eq:spectral-derivative-convergence}
is less than $1/m$ in the same directions.
Combining the three estimates gives
\[
 \left|
 \int_{\R^L}\rho(t)\,
 p_m^{-1/d}\int_G h_\ell\,dV_{p_m,t}\,dt
 -L_g(h_\ell)
 \right|
 <\frac3m
 \qquad\text{for }0\leq\ell\leq m.
\]

By Lemma~\ref{lem:generic-perturbation-averaging}, the map
\[
 t\longmapsto
 \left(
 p_m^{-1/d}\int_Gh_\ell\,dV_{p_m,t}
 \right)_{\ell=0}^m
\]
is bounded and measurable.
Its average against $\rho(t)\,dt$ belongs to the closed
convex hull of its values at parameters in
$\operatorname{supp}\rho$ where the varifolds realize
the widths and have nondegenerate components.
By Carath\'eodory's theorem, this average can be approximated
within $1/m$ in each coordinate by a convex combination
of at most $m+2$ such values.

Let $t_1,\ldots,t_{N_m}$ be the selected parameters and
$\lambda_{m,1},\ldots,\lambda_{m,N_m}$ the corresponding
weights. Set
\[
 F_{m,j}:=\widetilde H(t_j),
 \qquad
 V_{m,j}:=V_{p_m,t_j}
 \qquad\text{for }1\leq j\leq N_m.
\]
These varifolds have the required realization and
nondegeneracy properties, and the corresponding measure
$\mu_m$ satisfies
\[
 \left|
 \int_Gh_\ell\,d\mu_m-L_g(h_\ell)
 \right|
 \leq\frac4m
 \qquad\text{for }0\leq\ell\leq m.
\]

Recall that
\[
 L_g(h)=\mathcal W(1)\int_Gh\,d\Lambda_g.
\]
Taking $h_0=1$ therefore gives
$\mu_m(G)\to\mathcal W(1)$.
The measures $\mu_m$ have uniformly bounded mass, so
density of the functions $h_\ell$ extends the convergence
to every $h\in C(G)$.
This proves the asserted weak convergence.
Since $\mathcal W(1)>0$, the normalized measures converge
to $\Lambda_g$ as well.

Finally, for every fixed $k$ and all sufficiently large $m$,
the estimate for $\widetilde H$ gives
\[
 \max_{1\leq j\leq N_m}\|F_{m,j}-1\|_{C^k(G)}<\frac1m.
\]
This proves the convergence of the anisotropies.
\end{proof}

\subsection{Proof of equidistribution}
\label{subsec:proof-grassmann-equidistribution}

We now pass from the finite weighted collections of
Proposition~\ref{prop:geometric-grassmann-collections}
to a single sequence of hypersurfaces.

\begin{lemma}\label{lem:grassmann-block-sequence}
There exists a sequence of connected components of the varifolds
$V_{m,j}$, with repetitions allowed, whose cumulative area
measures, normalized by their total mass, converge weakly
to $\Lambda_g$ on $G$.
Moreover, the total area diverges, and the indices $m$
of the selected components tend to infinity.
\end{lemma}

\begin{proof}
Recall that $\mu_m/\mu_m(G)\rightharpoonup\Lambda_g$.
For every $m$, discard the components with zero coefficients
and approximate the remaining coefficients in
$\sum_j\lambda_{m,j}V_{m,j}$ by positive rational numbers,
so that the corresponding normalized measure changes by
at most $1/m$ in total variation.
Clearing denominators gives a finite list $B_m$ of
hypersurfaces, with repetitions, whose normalized area
measures converge weakly to $\Lambda_g$.

Fix these lists and concatenate sufficiently many copies
of each $B_m$, choosing the number of repetitions so that
the total area accumulated at the end of stage $m$
is at least $m$ times both the area of all preceding stages
and the area of $B_{m+1}$, and is at least $m$.

At the end of stage $m$, the contribution of all preceding
stages is at most $1/m$ of the total area. Hence the
normalized measures converge to $\Lambda_g$ at these
endpoints.

We claim that convergence holds for the entire sequence,
and not only at the end of each stage.
Indeed, during stage $m+1$, the normalized area measure
of the preceding stages together with the completed copies
of $B_{m+1}$ is a convex combination of two measures
converging to $\Lambda_g$.
The possible unfinished copy has area at most $1/m$
of the area already accumulated, and therefore changes
the normalized integral of any $h\in C(G)$ by at most
$2\|h\|_\infty/m$.

Since each stage is finite, the indices $m$ tend to
infinity. The construction also ensures that the total
area diverges.
\end{proof}

\begin{proof}[Proof of Theorem~\ref{thm:fixed-metric-anisotropic-equidistribution}]
By Proposition~\ref{prop:geometric-grassmann-collections}
and Lemma~\ref{lem:grassmann-block-sequence}, we obtain a
sequence of smooth, closed, embedded, connected hypersurfaces
$\{\Sigma_j\}_j$ satisfying
\eqref{eq:fixed-g-anisotropic-sequence}, with divergent total
area. Each $\Sigma_j$ is a nondegenerate component of some
$V_{m,i}$ realizing the width
$\omega_{p_m}^{d-1}(M,F_{m,i})$.
We denote by $F_j$ the corresponding anisotropy.
In particular, $\Sigma_j$ is $F_j$-stationary.
Since the indices $m$ tend to infinity, Proposition~
\ref{prop:geometric-grassmann-collections} also gives
$F_j\to 1$ in $C^\infty(G)$.

Finally, suppose that only finitely many distinct hypersurfaces
occur in the sequence. The union of their lifts to $G$ is a
closed set of $\Lambda_g$-measure zero, since its projection
to $M$ has zero $d$-dimensional volume.
Since all the measures appearing on the left-hand side of
\eqref{eq:fixed-g-anisotropic-sequence} are supported on this
set, they cannot converge weakly to $\Lambda_g$.
This contradiction completes the proof.
\end{proof}

\part{Smooth realization of the anisotropic widths}
\label{part:realization}

\section{Introduction and preliminaries}
In this part of the paper, we prove Theorem~\ref{thm:main-realization}. Indeed we show that, using min-max theory for anisotropies sufficiently close to the area, one can realize $p$-widths with a collection of smooth hypersurfaces, thus establishing Assumption~\ref{hyp:near-area-realization} for $3\leq d\leq 6$. We follow the Almgren-Pitts scheme developed for anisotropies in \cite{DPDRLi2024} and describe the modifications needed therein to adapt the situation to multi-parameter min-max and establish regularity.

The first main adaptation is the sweepout; in particular we need to work with multiparameter settings, we cannot perform the surgeries in \cite{DPDRLi2024}, since they are one-parameter in nature, as they rely on the construction of nested sweepouts introduced in~\cite{ChambersRotman}, and subsequently refined in \cite{Chambers-Liokumovich2020,Chodosh-Liokumovich-Spolaor2022}. This means that we cannot expect an a priori bound on the density ratio and we rely on Schoen-Simon type curvature estimates for anisotropies near the area functional. For this we use the Bernstein theorems of Lin~\cite{Lin1990},
Chodosh-Li~\cite[Corollary~1.5]{ChodoshLi2023},
and Li-Xia~\cite[Corollary~3]{LiXia2026} for anisotropies near the area functional and apriori smoothness for a $C^3$-neighbourhood of the area in \cite{SSA1977}, with a contradiction argument.

Second, gaining the interior smoothness of replacements from the above theorem, we prove the smooth gluing of consecutive replacements using a boundary-barrier argument inspired by \cite{Hardt1977} and the replacement ideas of \cite{DPDRLi2024}.

The remaining possible singularities are isolated, and
the hypersurface is stable in a punctured neighborhood
of each such point.
In the isotropic case, monotonicity supplies the density
bounds, which in turn imply the isolated singularities have zero capacity. This means they can be ignored in the stability inequality.
Since here the monotonicity formula is out of reach, we use instead the removable singularity strategy employing the separation function developed in \cite{DRHW2026} and the near-area curvature
inequalities of \cite{Winklmann2005}.

\subsection{The setup}

Throughout this part, $(M^d,g)$ is closed, connected, and smooth,
with $3\leq d\leq6$. We consider smooth integrands satisfying
\[
 \|F-1\|_{C^{k_0}(G_{d-1}(M))}<\delta_0,
 \qquad k_0\geq6,
\]
where $\delta_0>0$ will be chosen sufficiently small,
independently of the number of parameters.
In particular, these integrands are uniformly strictly elliptic.
We retain the notation and conventions of
Section~\ref{subsec:sweepouts-widths}.

For a varifold $V$ and a smooth vector field $Y$, with flow
$\psi_t$, we write
\[
 \delta\F(V)[Y]
 :=
 \left.\frac{d}{dt}\right|_{t=0}
 \F\bigl((\psi_t)_\#V\bigr).
\]
We call $V$ \emph{$F$-stationary} in an open set $U\subset M$
if $\delta\F(V)[Y]=0$ for every smooth vector field compactly
supported in $U$.
A smooth $F$-stationary hypersurface is \emph{stable} in $U$
if its second variation is nonnegative for every smooth
normal variation compactly supported in $U$.
For a one-sided hypersurface, variations are understood
as sections of its normal bundle.

We work with the space $\mathcal Z_{d-1}^0(M;\Z_2)$ of
cycles in the flat path component of zero.
We denote the usual varifold metric by $\bF$ and equip
this cycle space with the metric
\[
 \bF(T,T')
 :=
 \flat(T,T')+\bF(|T|,|T'|).
\]
The energy $T\mapsto\F(T)$ is continuous in this topology.
Moreover, a $\bF$-continuous family on a compact parameter
space has no concentration of mass: its associated mass
measures form a weakly compact family of nonatomic measures.

For the min-max construction, we allow $\bF$-continuous
families. This is consistent with the strong mass
formulation used earlier: strong mass continuity implies
$\bF$-continuity, while
Proposition~\ref{prop:mass-approximation} approximates every
such family by a mass-continuous family in the same flat
homotopy class, with arbitrarily small increase in the
maximal anisotropic energy.

We use the continuous formulation of the Almgren-Pitts
homotopy classes in
\cite[Section~2]{MarquesNeves2021}.

\begin{definition}\label{def:cycle-homotopy-class}
Let $X$ be a finite cubical subcomplex of some $I(m,j)$,
and let
\[
 \Phi_0:X\longrightarrow\mathcal Z_{d-1}^0(M;\Z_2)
\]
be $\bF$-continuous.
We define its homotopy class $\Pi$ to consist of all
$\bF$-continuous maps
$\Phi:X\to\mathcal Z_{d-1}^0(M;\Z_2)$ which are homotopic
to $\Phi_0$ in the flat topology.
Thus there is a flat-continuous map
\[
 \Psi:[0,1]\times X\longrightarrow
 \mathcal Z_{d-1}^0(M;\Z_2)
\]
such that $\Psi(0,\cdot)=\Phi_0$ and $\Psi(1,\cdot)=\Phi$.
\end{definition}

The associated min-max value is
\begin{equation}\label{eq:cycle-minmax-value}
 \mathbf L(\Pi)
 :=
 \inf_{\Phi\in\Pi}\sup_{x\in X}\F(\Phi(x)).
\end{equation}
For a sequence $S=\{\Phi_i\}_{i\in\N}\subset\Pi$, we set
\[
 \mathbf L(S)
 :=
 \limsup_{i\to\infty}\sup_{x\in X}\F(\Phi_i(x)).
\]
We call $S$ a \emph{min-max sequence} if
$\mathbf L(S)=\mathbf L(\Pi)$.
Such sequences exist by the definition of the infimum.
We will apply the construction when
\begin{equation}\label{eq:cycle-nontriviality}
 \mathbf L(\Pi)>0.
\end{equation}

Given a min-max sequence $S$, a sequence of slices
$\{\Phi_{i_j}(x_j)\}$, with $i_j\to\infty$, is
\emph{critical} if
\[
 \F(\Phi_{i_j}(x_j))\longrightarrow\mathbf L(\Pi).
\]
We define the critical set $\mathbf C(S)$ to be the
collection of varifold limits of critical sequences:
\[
 \mathbf C(S)
 :=
 \left\{
 V:
 |\Phi_{i_j}(x_j)|\rightharpoonup V,\quad
 \F(\Phi_{i_j}(x_j))\to\mathbf L(\Pi),\quad i_j\to\infty
 \right\}.
\]
By comparability,
\begin{equation}\label{eq:cycle-mass-bound}
 \limsup_{i\to\infty}\sup_{x\in X}\M(\Phi_i(x))
 \leq\lambda^{-1}\mathbf L(\Pi).
\end{equation}
Consequently, varifold compactness and continuity of
$\F$ under varifold convergence imply that
$\mathbf C(S)$ is nonempty and compact, and
\begin{equation}\label{eq:cycle-critical-energy}
 \F(V)=\mathbf L(\Pi)
 \qquad\text{for every }V\in\mathbf C(S).
\end{equation}

Finally, if $\Phi_0$ is a $p$-sweepout, every map in
$\Pi$ is also a $p$-sweepout, by flat homotopy invariance
of $\bar\lambda^p$.
Proposition~\ref{prop:mass-approximation} and
\eqref{eq:cycle-minmax-value} therefore give
\begin{equation}\label{eq:cycle-width-comparison}
 \omega_p^{d-1}(M,F)
 \leq\mathbf L(\Pi)
 \leq\sup_{x\in X}\F(\Phi_0(x)).
\end{equation}

\subsection{Main results}

We state the min-max theorem for sequences of homotopy
classes whose parameter complexes lie in a fixed cube.
This includes the case of a single homotopy class and
allows us to deduce smooth realization of every width
directly.

\begin{thm}\label{thm:near-area-minmax}
Let $(M^d,g)$ be closed, connected, and smooth, with
$3\leq d\leq6$, and fix $k_0\geq6$.
There exists $\delta_0>0$ with the following property.
Let $F$ be a smooth positive integrand satisfying
\[
 \|F-1\|_{C^{k_0}(G_{d-1}(M))}<\delta_0,
\]
and let $\Pi_j$ be homotopy classes as in
Definition~\ref{def:cycle-homotopy-class}, with parameter
complexes $X_j$ contained in subdivisions of a fixed
cube $I^m$. Suppose that
\[
 \mathbf L(\Pi_j)\longrightarrow L\in(0,\infty).
\]
Then there exist maps $\Phi_j\in\Pi_j$ satisfying
\[
 \sup_{x\in X_j}\F(\Phi_j(x))
 \leq\mathbf L(\Pi_j)+j^{-1}
\]
and, after passing to a subsequence, points $x_j\in X_j$
such that
\[
 |\Phi_j(x_j)|\rightharpoonup
 V=\sum_{i=1}^s m_i|\Sigma_i|,
 \qquad
 \F(V)=L.
\]
Here $m_i\in\N$, and the $\Sigma_i$ are smooth, closed,
embedded, connected, pairwise disjoint $F$-stationary
hypersurfaces.

The constant $\delta_0$ is independent of $m$, the
homotopy classes, and their min-max values.
\end{thm}

In particular, taking $\Pi_j=\Pi$ with
$\mathbf L(\Pi)>0$ gives a min-max sequence
$S\subset\Pi$ and a smooth $F$-stationary integral
varifold $V\in\mathbf C(S)$, satisfying $\F(V)=\mathbf L(\Pi)$.
Theorem~\ref{thm:near-area-minmax} will be proved in
the remaining sections of this part. We first derive
the smooth realization theorem.

\begin{proof}[Proof of Theorem~\ref{thm:main-realization}]
Choose $\delta_0$ as in
Theorem~\ref{thm:near-area-minmax}, and fix an integrand
$F$ in this neighborhood and $p\geq1$.
By the definition of the width, there are
mass-continuous $p$-sweepouts
\[
 \Phi_{0,j}:X_j\longrightarrow
 \mathcal Z_{d-1}^0(M;\Z_2)
\]
such that
\[
 \sup_{x\in X_j}\F(\Phi_{0,j}(x))
 \leq\omega_p^{d-1}(M,F)+j^{-1}.
\]

We may choose the $X_j$ as cubical subcomplexes of
subdivisions of the fixed cube $I^{2p+1}$, by the
parameter-space reduction in
\cite[proof of Proposition~3.2]{Li2023}.
Indeed, restricting a sweepout to the $p$-skeleton
preserves the sweepout condition, since restriction
is injective on $H^p(\,\cdot\,;\Z_2)$.
After a piecewise-linear embedding of this skeleton
into the interior of $I^{2p+1}$, take a cubical
neighborhood contained in a regular neighborhood
and precompose the restricted sweepout with the
retraction onto the skeleton.
The retraction has a section, so the pullback of
$\bar\lambda^p$ remains nonzero.
These operations preserve mass continuity and do
not increase the maximal anisotropic energy.

Let $\Pi_j$ be the homotopy class of $\Phi_{0,j}$.
By \eqref{eq:cycle-width-comparison},
\[
 0<\omega_p^{d-1}(M,F)
 \leq\mathbf L(\Pi_j)
 \leq\sup_{x\in X_j}\F(\Phi_{0,j}(x))
 \leq\omega_p^{d-1}(M,F)+j^{-1}.
\]
Consequently,
\[
 \mathbf L(\Pi_j)\longrightarrow
 \omega_p^{d-1}(M,F).
\]
The parameter complexes lie in the fixed cube
$I^{2p+1}$, so Theorem~\ref{thm:near-area-minmax}
provides an integral varifold
\[
 V=\sum_{i=1}^s m_i|\Sigma_i|,
 \qquad
 \F(V)=\omega_p^{d-1}(M,F),
\]
where the $\Sigma_i$ have all the asserted
smoothness, embeddedness, and stationarity
properties.
Since the constant $\delta_0$ in that theorem is
independent of the parameter dimension, the same
$\delta_0$ works for every $p$.
\end{proof}
    \begin{rmk}
The stable Bernstein theorem in $\mathbb R^7$
of \cite{HLW-bernstein-7} motivates studying whether
the present result extends to ambient dimension seven.
Such an extension would require suitable near-area
analogues of the stable curvature and
removable-singularity estimates used below.
\end{rmk}
\section{Almost minimizing varifolds}\label{sec:almost-minimizing}

\subsection{Pull tight}\label{subsec:pull-tight}

Fix the homotopy classes $\{\Pi_j\}_j$, parameter complexes
$X_j\subset I^m$, and limiting value
\[
 \mathbf L(\Pi_j)\longrightarrow L\in(0,\infty)
\]
from the hypotheses of
Theorem~\ref{thm:near-area-minmax}.
For sequences $S=\{\Phi_j\}$ with $\Phi_j\in\Pi_j$
and $\sup_{X_j}\F(\Phi_j)\to L$, we extend the previous
notation by setting
\[
 \mathbf C(S)
 :=
 \left\{
 V:
 |\Phi_{j_k}(x_k)|\rightharpoonup V,\quad
 \F(\Phi_{j_k}(x_k))\to L, \quad 
 j_k\to\infty,\quad x_k\in X_{j_k}
 \right\}.
\]
Thus every $V\in\mathbf C(S)$ satisfies $\F(V)=L$.

\begingroup
\emergencystretch=1em
We use the cycle formulation of almost minimizing.
A cycle $T\in\mathcal Z_{d-1}^0(M;\Z_2)$ is
\emph{$(\varepsilon,\delta)$-almost minimizing for $F$
in $U$} if every finite sequence
$T=T_0,\ldots,T_q\in\mathcal Z_{d-1}^0(M;\Z_2)$ satisfying
\[
 \operatorname{spt}(T_\ell-T)\subset U,\qquad
 \flat(T_\ell,T_{\ell-1})\leq\delta,\qquad
 \F(T_\ell)\leq\F(T)+\delta
 \quad(1\leq\ell\leq q)
\]
also satisfies
\[
 \F(T_q)\geq\F(T)-\varepsilon.
\]
A varifold $V$ is \emph{almost minimizing for $F$ in $U$}
if there exist $\varepsilon_i,\delta_i\downarrow0$ and
$(\varepsilon_i,\delta_i)$-almost minimizing cycles
$T_i$ in $U$ such that $|T_i|\rightharpoonup V$.
\par\endgroup

\begin{proposition}
\label{prop:pull-tight-almost-minimizing}
There exist representatives $\Phi_j\in\Pi_j$ satisfying
\[
 \sup_{x\in X_j}\F(\Phi_j(x))
 \leq\mathbf L(\Pi_j)+j^{-1}
\]
such that, writing $S=\{\Phi_j\}$, the following hold:
\begin{itemize}[leftmargin=0.6cm]
\item[(i)]
Every $W\in\mathbf C(S)$ is $F$-stationary:
\begin{equation}\label{eq:bounded-anisotropic-first-variation}
 \delta\F(W)[Y]=0
 \qquad\text{for every smooth vector field }Y.
\end{equation}

\item[(ii)]
There exists $V\in\mathbf C(S)$ which is almost
minimizing for $F$ in small annuli: for every $x\in M$,
there is $r_x>0$ such that $V$ is almost minimizing in
\[
 B_t(x)\setminus\overline{B_s(x)}
 \qquad\text{whenever }0<s<t<r_x.
\]
In particular, $V\neq0$ and $\F(V)=L$.

Moreover, there exist finitely many points
$p_1,\ldots,p_N\in M$ such that, setting
\[
 P:=\{p_1,\ldots,p_N\}\subset M,
\]
the varifold $V$ is almost minimizing in a sufficiently
small ball about every point of $M\setminus P$.
The ball radii can be chosen uniformly on each compact
subset of $M\setminus P$.
\end{itemize}
\end{proposition}

\begin{proof}
Choose representatives $\Phi_j^0\in\Pi_j$ with maximal
energy at most $\mathbf L(\Pi_j)+(2j)^{-1}$.
Comparability gives a uniform mass bound. Apply the
pull-tight construction of
\cite[Section~2.8]{MarquesNeves2021}, replacing the area
functional and its first variation by $\F$ and $\delta\F$,
as in \cite[Section~4.1]{DPDR2024}.
The resulting deformation preserves the homotopy classes,
does not increase $\F$, and makes every critical varifold
$F$-stationary.

By \cite[Proposition~3.7]{MarquesNeves2021}, we may further
approximate these families uniformly in the $\bF$-metric
by mass-continuous maps in the same classes. Choose the
approximation errors tending to zero and the energy
increase at most $(2j)^{-1}$. Since $\F$ is continuous
under varifold convergence, this preserves the critical
set. Denoting the resulting sequence by $S=\{\Phi_j\}$,
we obtain the asserted energy bound and \textup{(i)}.

For \textup{(ii)}, use the combinatorial deformation
argument of \cite[Theorem~4.10]{Pitts1981};
see \cite[proof of Lemma~3.2]{Li2023} for varying
parameter complexes and
\cite[proof of Theorem~4.6]{DPDRLi2024} for the
anisotropic adaptation of the existence argument.
The required modification is to measure the retained
pieces by their $\F$-energy, use local lower
semicontinuity for pieces substituted from a flat limit,
and estimate the added slicing chains by
$\lambda^{-1}$ times their mass, as in
Proposition~\ref{prop:mass-approximation}.
The ordinary mass estimates continue to control
fineness, while
\begin{align}
 |\F(T)-\F(T')|
 \leq \lambda^{-1}\M(T-T')
\end{align}
controls the energy errors in the extensions.

The same modification applies to the localized
interpolation underlying
\cite[Theorem~3.9]{Pitts1981}: almost minimization
with respect to mass-small steps in $U$ implies
almost minimization in the stated flat-step sense
in every $W\Subset U$, with possibly different
witnessing cycles and thresholds.
Thus mass-step almost minimization in sufficiently
small annuli gives the required flat-step property,
by applying this implication in slightly larger annuli.

All combinatorial constants are uniform because the
$X_j$ lie in subdivisions of the fixed cube $I^m$.
Consequently, if no element of $\mathbf C(S)$ were
almost minimizing in small annuli, the cited argument
would produce maps $\Psi_j\in\Pi_j$ satisfying
\[
 \limsup_{j\to\infty}\sup_{x\in X_j}\F(\Psi_j(x))<L,
\]
contrary to $\mathbf L(\Pi_j)\to L$.
Hence such a $V\in\mathbf C(S)$ exists; moreover,
$\F(V)=L>0$, so $V\neq0$.

Finally, choose a finite subcover
\[
 M=\bigcup_{a=1}^N B_{r_{p_a}/2}(p_a),
 \qquad P=\{p_1,\ldots,p_N\}.
\]
For every $y\in M\setminus P$, some sufficiently small
ball about $y$ is contained in an annulus
$B_t(p_a)\setminus\overline{B_s(p_a)}$ with
$0<s<t<r_{p_a}$.
Almost minimization passes to smaller open sets,
giving the ball conclusion.
Covering each compact subset of $M\setminus P$
by finitely many smaller such balls gives a uniform
choice of radius there.
\end{proof}

\subsection{Curvature estimates and regularity of replacements}\label{subsec:curvature-estimate-regularity-replacements}
For the area functional, Bellettini~\cite{Bellettini2025}
obtained extensions of the classical curvature estimates
of Schoen--Simon--Yau~\cite{SSY1975} and the regularity
theory of Schoen--Simon~\cite{SS1981}
using De Giorgi iteration. In the anisotropic setting, Allard~\cite{Allard1983}
established a priori estimates for the oscillation of
the normal under control of the first and second
variations of an elliptic integrand.

We first establish interior curvature estimates and compactness for
smooth stable hypersurfaces. We then construct replacements
and prove their interior regularity. Boundary regularity
requires quantitative control of the prescribed boundary,
which we state separately.

\begin{proposition}
\label{prop:stable-curvature-compactness}
After decreasing $\delta_0>0$, there exist $r_0>0$ and
$C>0$, depending only on $(M,g)$ and $k_0$, with the
following properties.
\begin{itemize}[leftmargin=0.6cm]
\item[(i)]
Suppose that $\|F-1\|_{C^{k_0}}\leq\delta_0$ and that
$\Sigma$ is a smooth, properly immersed, two-sided,
stable $F$-stationary hypersurface without boundary
in $U$. Then
\begin{equation}\label{eq:stable-curvature-estimate}
 \sup_{\Sigma\cap B_r(x)}|A_\Sigma|
 \leq \frac{C}{r}
\end{equation}
whenever $B_{2r}(x)\Subset U$ and $0<r\leq r_0$.
The estimate is independent of the area of $\Sigma$.

\item[(ii)]
Suppose that $F_j\to F_\infty$ smoothly, with
$\|F_j-1\|_{C^{k_0}}\leq\delta_0$, and that
$\Sigma_j$ are smooth, properly embedded, two-sided,
stable $F_j$-stationary hypersurfaces without boundary
in $U$. Assume that their areas are locally uniformly
bounded in $U$.

After passing to a subsequence, $|\Sigma_j|$
converges to an integral varifold whose support is a
smooth, embedded, stable $F_\infty$-stationary
hypersurface. The convergence is locally smooth as
graphical sheets, with positive integer multiplicities
constant on connected components. The limit may be empty.
\end{itemize}
\end{proposition}

\begin{proof}
The Bernstein theorems of Lin~\cite{Lin1990},
Chodosh-Li~\cite[Corollary~1.5]{ChodoshLi2023},
and Li-Xia~\cite[Corollary~3]{LiXia2026} imply that
every complete, two-sided, stable immersion in
$\R^d$, $3\leq d\leq6$, stationary for an autonomous
integrand sufficiently $C^4$-close to area is flat. (For the isotropic stable Bernstein theorem, see
\cite{doCarmoPeng1979,FischerColbrieSchoen1980,Pogorelov1981}
in $\R^3$,
\cite{ChodoshLi2024,CatinoMastroliaRoncoroni2024}
in $\R^4$,
\cite{ChodoshLiMinterStryker2026} in $\R^5$,
and \cite{Mazet2024} in $\R^6$.) These are the main tools used in this proof.

Choose
$\delta_0$ so that the frozen integrands below and
their limits belong to this neighborhood.

If \eqref{eq:stable-curvature-estimate} fails,
point-picking gives points $y_j$, numbers
$\lambda_j=|A_{\Sigma_j}|(y_j)\to\infty$, and radii
$\rho_j$ such that
\[
 \lambda_j\rho_j\to\infty,
 \qquad
 \sup_{\Sigma_j\cap B_{\rho_j}(y_j)}|A_{\Sigma_j}|
 \leq 2\lambda_j.
\]
Rescaling normal coordinates by $\lambda_j$, the
ambient metrics converge to the Euclidean metric
and the integrands converge locally in $C^4$ to
an autonomous integrand in the Bernstein neighborhood.
Uniform intrinsic graphical charts and elliptic
estimates give a pointed limiting immersion satisfying
$|A|(0)=1$, and $|A|\leq2$.
It is complete, two-sided, stationary, and stable:
the rescaled domains exhaust intrinsic neighborhoods
of arbitrarily large radius, and the equation and
stability inequality pass to the limit on compact
intrinsic subsets. This contradicts the Bernstein
theorem. Only a selected pointed sheet is used,
so no area bound enters the argument.

For \textup{(ii)}, \textup{(i)} gives uniform curvature
bounds on compact subsets. Elliptic estimates give
bounds for all derivatives, and the local area bounds
control the number of sheets in smaller balls.
A diagonal argument therefore gives smooth graphical
convergence with integer multiplicities.
Embeddedness prevents transverse intersections, and
the strong maximum principle implies that touching
limiting sheets coincide. Thus the support is embedded
and the multiplicity is constant on each component.
Stationarity and stability pass to the limit.
For a one-sided limiting component, stability is
understood for compactly supported sections of its
normal bundle.
\end{proof}

We next record a boundary version for a quantitatively
controlled prescribed boundary. Smooth embedded
convergence below includes a common local graphical
radius; in particular, distinct prescribed boundary
sheets do not collide.

\begin{proposition}
\label{prop:boundary-multiplicity}
Assume the hypotheses of
Proposition~\ref{prop:stable-curvature-compactness}.
Work in local coordinates $(x,z)\in\mathbb R^{d-1}\times\mathbb R$,
and let $\Sigma_0=\{z=u_0(x)\}$ be a smooth $F$-minimal hypersurface.
For $D=B_r^{d-1}$ and $a>0$, consider the lens
\[
 \begin{gathered}
 C=\left\{(x,z):x\in D,\quad
 |z-u_0(x)|<a\,\frac{r^2-|x|^2}{r}\right\},
 \\
 \Gamma=\{(x,u_0(x)):x\in\partial D\}.
 \end{gathered}
\]
Choose the coordinates and then $r$ sufficiently small that
$\|Du_0\|+r\|D^2u_0\|\le a$ and both open faces of $C$
are strictly $F$-convex.

Let $p\in\Gamma$, and let $W$ be an $F$-stationary varifold
in a neighborhood $U$ of $p$. Suppose that, for some
$m\in\mathbb N$,
\[
 W=m|\Sigma_0|
 \quad\text{in }U\setminus\overline C,
\]
and that $W\llcorner C$ is integral and has smooth embedded support which is $F$-stable with respect to every
compactly supported normal variation in $U\cap C$.

If $a$ is sufficiently small, then near $p$ the interior
restriction of $W$ is a finite sum, with positive integer
multiplicities, of smooth hypersurfaces extending smoothly
to $\Gamma$. Moreover, if $a$ is chosen sufficiently small depending also
on $m$, then
\[
 W=m|\Sigma_0|
 \qquad\text{in a neighborhood of }p.
\]
\end{proposition}

\begin{proof}
The boundary-barrier strategy is inspired by Hardt's
work~\cite{Hardt1977}; see also the proof of
\cite[Lemma~5.2]{DPDRLi2024}.
Here we use the interior curvature estimate and the
narrow-lens geometry to obtain graphicality and smooth
extension to the boundary. We then use the anisotropic
conormal balance and integrality to identify the
multiplicity and prove smooth matching with the
exterior hypersurface.

We first exclude contact with the open faces of the lens.
At every point of $\partial C\setminus\Gamma$, the exterior
sheet $\Sigma_0$ is absent from a sufficiently small
neighborhood. Consequently, stationarity of $W$ across the
face and its strict $F$-convexity exclude support contact,
by the anisotropic barrier maximum principle; see
\cite[Section~3]{DDH2019}.

Shrink $U$ to a coordinate neighborhood for which
$U\cap C$ is simply connected. Every properly embedded
component of the interior support is then two-sided in
$U\cap C$. Its choice of normal also orients its extension
by the empty set across the open faces.
Thus the interior support, extended by the empty set across
the faces, is a smooth hypersurface without boundary in
$U\setminus\Gamma$. It is stable there: the restriction of
any compactly supported variation to this support can be
cut off inside $C$, since its intersection with the support
stays a positive distance from the faces.
The interior curvature estimate therefore gives, after
shrinking $U$,
\begin{equation}\label{eq:lens-curvature}
 |A|(q)\le \frac{C}{\operatorname{dist}(q,\Gamma)}.
\end{equation}
Crucially, $C$ is independent of the opening $a$.

Put $f(x,z)=z-u_0(x)$. The definition of the lens gives
\[
 |f(x,z)|
 \le 2a\,\operatorname{dist}(x,\partial D)
 \le Ca\,\operatorname{dist}((x,z),\Gamma).
\]
At a point whose distance from $\Gamma$ is $s$,
\eqref{eq:lens-curvature} bounds the second derivative of
$f$ along unit-speed intrinsic geodesics by $C/s$.
Applying Taylor's formula along geodesic segments of length
$s\sqrt a$ yields $|\nabla_\Sigma f|\le C\sqrt a$.
Hence, for small $a$, projection onto the $x$-coordinates
is a local diffeomorphism.

Over a smaller simply connected boundary patch of $D$,
this projection is proper: the lens bounds prevent vertical
escape, and the exclusion of face contact prevents any
additional boundary. It is therefore a covering.
Local finiteness of the mass makes its degree finite.
Accordingly, the interior support consists of finitely many
ordered graphs $z=u_\ell(x)$, with constant positive integer
multiplicities $q_\ell$, satisfying
\[
 u_\ell=u_0\quad\text{on }\partial D,
 \qquad
 |Du_\ell-Du_0|\le C\sqrt a.
\]
Each graph solves the uniformly elliptic scalar equation
for an $F$-minimal graph. Boundary $C^{1,\alpha}$ estimates
for the Dirichlet problem, followed by Schauder estimates,
give smooth extension to $\Gamma$.

There is no additional varifold mass concentrated on
$\Gamma$. Indeed, test stationarity with cutoffs of the
normal displacement field from $\Gamma$, supported in
shrinking tubular neighborhoods. On $\Gamma$ the field
vanishes and its derivative is projection onto the
two-dimensional normal space. For the area integrand,
its tangential divergence on every hyperplane is at least
one; this remains positive when $F$ is sufficiently close
to area. Contributions off $\Gamma$ tend to zero, proving
$\|W\|(\Gamma)=0$.

We finally prove the last assertion. Along $\Gamma$,
let $\beta_0$ and $\beta_\ell$ denote the anisotropic
conormal vectors of $\Sigma_0\cap C$ and the interior
graphs, with the same outward convention. Stationarity gives
\[
 m\beta_0=\sum_\ell q_\ell\beta_\ell.
\]
More explicitly, at $p\in\Gamma$ let
$N_p\Gamma=(T_p\Gamma)^{\perp_g}$, let $\eta_\ell$ be the
outward unit conormal, and choose unit normals $\nu_\ell$
so that one quarter-turn $R_p$ of this two-plane satisfies
$R_p\nu_\ell=\eta_\ell$, including $\ell=0$. In the normal
representation introduced above,
\[
 \beta_\ell
 =\varphi_p(\nu_\ell)\eta_\ell
  -\langle D\varphi_p(\nu_\ell),\eta_\ell\rangle_g\nu_\ell
 =R_pD(\varphi_p|_{N_p\Gamma})(\nu_\ell).
\]
Here the last derivative is the gradient in $N_p\Gamma$.
These vectors lie on the rotated Wulff curve of
$\varphi_p|_{N_p\Gamma}$, which is a smooth strictly convex
curve enclosing the origin by ellipticity.
Let $L$ be its supporting linear functional at $\beta_0$,
normalized by $L(\beta_0)=1$. Then
$L(\beta_\ell)\le1$, with equality only when
$\beta_\ell=\beta_0$.

Since the graphs are smooth up to the boundary, the lens
height bound also gives
$|Du_\ell-Du_0|\le Ca$ on $\partial D$.
We may therefore choose $a$, depending on $m$, so that
\[
 L(\beta_\ell)>\frac{m}{m+1}
 \qquad\text{for every }\ell.
\]
Writing $N=\sum_\ell q_\ell$ and applying $L$ to the
conormal identity gives
\[
 \frac{m}{m+1}N<m\le N.
\]
Thus $N=m$, and equality in the supporting inequality
forces every $\beta_\ell=\beta_0$.

Each interior graph consequently joins $\Sigma_0$ with
matching first derivatives. The joined graph is a weak
solution of the graphical Euler-Lagrange equation and
is smooth by elliptic regularity. Unique continuation
identifies it with $\Sigma_0$. Summing its multiplicities
proves $W=m|\Sigma_0|$ near $p$.
\end{proof}

\textit{Smooth replacements.}
Let $V$ be an $F$-stationary varifold in $M$ and let
$K\subset M$ be compact. A \emph{smooth $F$-replacement}
of $V$ in $K$ is an $F$-stationary varifold $V^*$ such that
\[
 V^*\res(M\setminus K)=V\res(M\setminus K),
 \qquad
 \F(V^*)=\F(V),
\]
and whose restriction to $\operatorname{int}K$ is integral,
with smooth embedded $F$-stable support and positive integer
multiplicities constant on connected components.
The interior restriction may be zero. Stability is understood
with respect to compactly supported normal variations;
the limiting support need not be two-sided. We prove below existence and regularity of replacements.

\begin{proposition}
\label{prop:smooth-replacements}
Assume the dimensional and near-area hypotheses of
Proposition~\ref{prop:stable-curvature-compactness},
decreasing $\delta_0$ if necessary.
Let $V$ be $F$-stationary in $M$ and almost minimizing
for $F$ in an open set $U$.

\begin{itemize}[leftmargin=0.6cm]
\item[(i)]
For every compact $K\Subset U$, there exists a smooth
$F$-replacement $V^*$ of $V$ in $K$ which remains
almost minimizing in $U$.
Its interior support satisfies the curvature estimate
\eqref{eq:stable-curvature-estimate} whenever
$B_{2r}(x)\Subset\operatorname{int}K$ and $r\le r_0$,
with the same constant $C$.

\item[(ii)]
The construction can be iterated any finite number of times
in compact subsets of the same open set $U$.
At every stage, stationarity, almost minimization in $U$,
and the total anisotropic energy are preserved.
The replacement regions may have corners; in particular,
they may be closed round annuli or closed adapted lens-annuli.

\item[(iii)]
Suppose that, near an outer-edge point $p\in\Gamma$,
the region $K$ agrees with the closure of an adapted lens
$C$ from Proposition~\ref{prop:boundary-multiplicity},
and that
\[
 V=m|\Sigma_0|
 \quad\text{outside }\overline C
 \quad\text{near }p,
 \qquad m\in\mathbb N.
\]
Then the boundary conclusions of that proposition apply
to $V^*$. In particular, for sufficiently small opening,
the interior sheets extend smoothly to $\Gamma$ with
integer multiplicities. If the opening is sufficiently
small depending also on $m$, then
\[
 V^*=m|\Sigma_0|
 \qquad\text{near }p.
\]
This applies at the outer edge of an adapted lens-annulus,
since its inner boundary is absent from a sufficiently
small neighborhood of $p$.
\end{itemize}

The threshold $\delta_0$ is independent of multiplicity
and of the number of replacements. In item \textup{(iii)}, only
the opening and size of the lens may depend on $m$.
\end{proposition}

\begin{proof}
Choose $(\varepsilon_i,\delta_i)$-almost minimizing
cycles $T_i$ in $U$ with $|T_i|\rightharpoonup V$.
We use constrained minimization among admissible endpoints;
compare the replacement construction in
\cite[Lemma~4.8 and Proposition~4.11]{DPDRLi2024}.

For each $i$, minimize $\F$ among endpoints of admissible
chains starting at $T_i$ whose changes are supported in $K$.
Compactness for mod-$2$ cycles and lower semicontinuity
give a flat limit of a minimizing sequence.
This limit is itself admissible: append it to a sufficiently
late endpoint by one flat step of size less than $\delta_i$.
Thus the minimum is attained by a cycle $T_i^*$, and
\begin{equation}\label{eq:replacement-endpoint-energy}
 \operatorname{spt}(T_i^*-T_i)\subset K,
 \qquad
 \F(T_i)-\varepsilon_i
 \le \F(T_i^*)\le \F(T_i).
\end{equation}

Concatenating a chain from $T_i$ to $T_i^*$ with any
admissible chain starting at $T_i^*$ shows that $T_i^*$
is still $(\varepsilon_i,\delta_i)$-almost minimizing
in $U$. Indeed, the latter chain has energy at most
$\F(T_i^*)+\delta_i\le\F(T_i)+\delta_i$, and its endpoint
has energy at least
$\F(T_i)-\varepsilon_i\ge\F(T_i^*)-\varepsilon_i$.

Moreover, $T_i^*$ locally minimizes $\F$ in
$\operatorname{int}K$. To see this, take a sufficiently
small coordinate ball $B\Subset\operatorname{int}K$
with $\operatorname{Vol}(B)<\delta_i$.
A cycle competitor agreeing with $T_i^*$ outside $B$
differs from it by the boundary of a mod-$2$ filling
supported in $B$, of mass at most $\operatorname{Vol}(B)$.
An energy-decreasing competitor would therefore be an
admissible additional endpoint, contradicting minimality.

Although these minimizing neighborhoods may shrink with
$i$, stationarity and stability hold throughout
$\operatorname{int}K$. Indeed, any sufficiently short
smooth variation compactly supported there can be
subdivided into flat steps smaller than $\delta_i$,
while its energy stays below $\F(T_i)+\delta_i$.
Endpoint minimality then makes $t=0$ a local minimum
along every such variation.

After decreasing $\delta_0$, the near-area minimizing
regularity theorem
\cite[Theorem~II.8 and Remark~II.9]{SSA1977}
implies that the support of $T_i^*$ is a $C^2$ embedded
hypersurface in $\operatorname{int}K$; the dimensional
restriction eliminates its singular set.
Indeed, in sufficiently small normal-coordinate
charts the $C^{k_0}$ bound, with $k_0\ge6$, places the frozen
normal integrands in a fixed compact family in the $C^3$
topology. After decreasing the coordinate radius and
$\delta_0$, the frozen integrands lie in the neighborhood of area
within this family required by Theorem~II.8. With hypersurface dimension
$d-1$, choose $t=8-d>0$ in that theorem; its singular-set
estimate becomes $\mathcal H^0(\operatorname{sing}T_i^*)=0$.
The compact family and the resulting smallness threshold
are independent of $i$, the minimizing radius, and the
number of parameters.
Remark~II.9 includes both the Riemannian and mod-$2$
versions needed here. Elliptic bootstrapping gives
smoothness. Since $T_i^*$ bounds a set of finite perimeter,
its smooth interior support is two-sided and carries
multiplicity one.

The energy bounds give a uniform mass bound. Passing
to a subsequence, write $|T_i^*|\rightharpoonup V^*$.
Proposition~\ref{prop:stable-curvature-compactness}(ii)
gives locally smooth graphical convergence in
$\operatorname{int}K$, with positive integer
multiplicities and smooth embedded stable support.
Its curvature estimate passes to this limit.
In particular, multiplicity in $V^*$ is retained even
when sheets cancel in the flat limit of the cycles.

The support condition and
\eqref{eq:replacement-endpoint-energy} imply that $V^*=V$ on $M\setminus K$, and $\F(V^*)=\F(V)$.
The cycles $T_i^*$ witness almost minimization of $V^*$
in $U$. This implies stationarity in $U$: a strictly
energy-decreasing flow could be discretized into
admissible flat steps and would contradict almost
minimization for large $i$.
Combining this with stationarity of $V$ outside $K$,
using a cutoff equal to one near $K$, proves stationarity
of $V^*$ on all of $M$. This establishes \textup{(i)}.

Since $V^*$ remains almost minimizing in the same $U$,
the construction can be repeated there, proving
\textup{(ii)}.

Finally, under the hypotheses of \textup{(iii)},
$V^*$ is stationary across the lens boundary, agrees
with $m|\Sigma_0|$ on its exterior, and has integral,
smooth embedded stable support on its interior.
These are exactly the hypotheses of
Proposition~\ref{prop:boundary-multiplicity}.
Its conclusions prove \textup{(iii)}.
\end{proof}

\subsection{Smooth gluing of consecutive replacements} At this point we prove the smooth gluing of consecutive replacements in overlapping annuli, and using this show the regularity of the limit varifold in almost minimizing regions.
\begin{proposition}
\label{prop:smooth-gluing}
Assume the dimensional and near-area hypotheses above and
Proposition~\ref{prop:boundary-multiplicity}.
Let $0<s<\rho<t<2\rho$, let $V'$ be an $F$-stationary replacement of $V$ in
$\overline{B_{2\rho}(x)}\setminus B_\rho(x)$, and let $V''$
be an $F$-stationary replacement of $V'$ in
$\overline{B_t(x)}\setminus B_s(x)$.
Assume that these replacements have smooth embedded stable
interiors with integer multiplicities, and that the replacement
construction can be iterated in compact regions contained in
sufficiently small admissible annuli.

Write $\Sigma'$ for the support of $V'$ in
$B_{2\rho}(x)\setminus\overline{B_\rho(x)}$.
If $\partial B_t(x)$ is strictly $F$-convex and transverse
to $\Sigma'$, then for every
$y\in\Sigma'\cap\partial B_t(x)$ there exists $r_y>0$ such that
\[
 V''\llcorner B_{r_y}(y)=V'\llcorner B_{r_y}(y).
\]
Consequently, $V''$ has smooth embedded support in
$B_{2\rho}(x)\setminus\overline{B_s(x)}$, with locally
constant positive integer multiplicities.
\end{proposition}

\begin{proof}
Fix $y\in\Sigma'\cap\partial B_t(x)$.
Choose a neighborhood $U$ of $y$ in which $V'=m|\Sigma_0|$,
 where  $m\in\mathbb N$ and $\Sigma_0$ is a smooth embedded $F$-minimal hypersurface.
In particular,
\[
 V''=m|\Sigma_0|
 \quad\text{in }U\setminus\overline{B_t(x)}.
\]
After passing to local coordinates and absorbing the metric
into the integrand, write
\[
 y=0,\qquad
 \Sigma_0=\{(z,u(z)):z\in B_{r_0}^{d-1}\},
 \qquad u(0)=0,\quad Du(0)=0.
\]
Shrink the neighborhood so that
$\Sigma_0\cap B_t(x)\cap U$ is connected.

Fix a small opening $a>0$ and define
\[
 C_r=
 \left\{(z,h):|z|<r,\quad
 |h-u(z)|<a\,\frac{r^2-|z|^2}{r}\right\},
 \qquad
 A_r=C_r\setminus\overline{C_{r/2}},
\]
with outer edge
\[
 \Gamma_r=\{(z,u(z)):|z|=r\}.
\]
For each fixed $a$, the two outer faces are strictly
$F$-convex when $r$ is sufficiently small.
Moreover, $\overline{A_r}$ is contained in a round annulus
centered at $y$, with inner radius comparable to $ar$ and
outer radius comparable to $r$.
Thus the assumed replacement construction provides a
further replacement $W_r$ of $V''$ in $\overline{A_r}$.

Choose $a$ sufficiently small for
Proposition~\ref{prop:boundary-multiplicity} with exterior
multiplicity $m$.
By transversality, every sufficiently small $\Gamma_r$
contains a point
$z_r\in\Gamma_r\setminus\overline{B_t(x)}$.
Near $z_r$, the exterior of the lens carries precisely
$m|\Sigma_0|$. The lens proposition therefore gives
\[
 W_r=m|\Sigma_0|
 \quad\text{near }z_r.
\]
The inner lens is disjoint from a sufficiently small
neighborhood of $z_r$, so it does not affect this application.

The set
\[
 \Sigma_0\cap A_r
 =
 \{(z,u(z)):r/2<|z|<r\}
\]
is connected, since $d-1\ge2$.
Interior smoothness, unique continuation, and constancy of
multiplicity propagate the preceding identity along this
entire sheet: $\Sigma_0\cap A_r$ is contained in
$\operatorname{spt}\|W_r\|$ and has multiplicity exactly $m$.
Indeed, the set of points where this local identification
holds is open and closed in $\Sigma_0\cap A_r$.

We claim that
\begin{equation}\label{eq:lens-edge-support}
 \Gamma_r\subset\operatorname{spt}\|V''\|.
\end{equation}
Every point $q\in\Gamma_r$ belongs to
$\operatorname{spt}\|W_r\|$ by the preceding continuation.
If $q\notin\operatorname{spt}\|V''\|$, then $V''$ vanishes
in a neighborhood of $q$. Since $W_r=V''$ outside
$\overline{A_r}$, the support of $W_r$ in that neighborhood
lies in $\overline{C_r}$.
The upper face of $C_r$ extends smoothly through $q$ as a
strictly $F$-convex hypersurface. It would therefore touch
the stationary varifold $W_r$ from its containing side,
contrary to the anisotropic barrier maximum principle.
This proves \eqref{eq:lens-edge-support}.

As $r$ varies, the edges $\Gamma_r$ cover the punctured
graph $\Sigma_0$ near $y$. Consequently,
\[
 \Sigma_0\cap U\subset\operatorname{spt}\|V''\|
\]
after shrinking $U$.
Inside $B_t(x)$, the support of $V''$ is smooth and embedded.
Thus it agrees locally with $\Sigma_0$ at every point of
$\Sigma_0\cap B_t(x)\cap U$.
Its multiplicity along this connected sheet is therefore
a fixed integer $k\ge1$. Other components of $V''$ are
not yet excluded.

We now identify $k$. Choose a new opening sufficiently
small for the lens proposition with both exterior
multiplicities $m$ and $k$, and take one new replacement
$W$ in the corresponding lens-annulus.
At an outer-edge point outside $\overline{B_t(x)}$, the
lens proposition and the same continuation argument give
multiplicity $m$ along $\Sigma_0\cap A_r$.
On the other hand, transversality also supplies
$q_r\in\Gamma_r\cap B_t(x)$. Since $q_r$ lies strictly inside the original annulus,
interior regularity gives
\[
 V''=k|\Sigma_0|
 \quad\text{in a whole neighborhood of }q_r.
\]
Applying the lens proposition at $q_r$ therefore gives
$W=k|\Sigma_0|$ near $q_r$.
Comparing the two identifications on the same connected
sheet $\Sigma_0\cap A_r$ proves $k=m$.

Finally, in $U$ the difference
\[
 R:=V''-m|\Sigma_0|
\]
is a nonnegative $F$-stationary varifold.
Nonnegativity follows from the multiplicity just established
inside $B_t(x)$ and the exterior identity; the transverse
intersection $\Sigma_0\cap\partial B_t(x)$ carries no
$(d-1)$-dimensional mass of $|\Sigma_0|$.
Moreover,
\[
 \operatorname{spt}\|R\|\cap U
 \subset\overline{B_t(x)}.
\]
Strict $F$-convexity of $\partial B_t(x)$ and the barrier
maximum principle exclude $y$ from
$\operatorname{spt}\|R\|$.
Hence $R$ vanishes in a neighborhood of $y$, proving
\[
 V''=m|\Sigma_0|=V'
 \quad\text{near }y.
\]

At points of $\partial B_t(x)$ outside
$\operatorname{spt}\|V'\|$, the same barrier principle
excludes support of $V''$.
Together with interior regularity, this proves the final
assertion.
\end{proof}

We now conclude the regularity away from finitely many points:

\begin{corollary}
\label{cor:regularity-away-finite-points}
Under the dimensional and near-area assumptions above,
let $V$ be the $F$-stationary varifold given by
Proposition~\ref{prop:pull-tight-almost-minimizing}.
There is a finite set $P\subset M$ such that
\[
 V\res(M\setminus P)=\theta|\Sigma|,
\]
where $\Sigma\subset M\setminus P$ is a smooth, properly
embedded $F$-stationary hypersurface and $\theta$ is a
locally constant positive integer-valued function on $\Sigma$.
For every $p\in P$, there is $\rho_p>0$ such that
$B_{\rho_p}(p)\cap P=\{p\}$ and $\Sigma$ is stable in
$B_{\rho_p}(p)\setminus\{p\}$.
\end{corollary}

\begin{proof}
Let $P$ be the exceptional set supplied by
Proposition~\ref{prop:pull-tight-almost-minimizing}.
The integrality argument of \cite[Lemma~4.15]{DPDR2024}
applies using Proposition~\ref{prop:smooth-replacements},
the barrier principle, and
Proposition~\ref{prop:stable-curvature-compactness}.
Hence $V$ is integral.

Fix $x\in M\setminus P$ and a neighborhood $U$ of $x$
in which $V$ is almost minimizing. Successive annular
replacements and Proposition~\ref{prop:smooth-gluing},
followed by the identification argument of
\cite[proof of Proposition~4.14, Step~2]{DPDR2024}
with $c=0$, show that $V$ has smooth properly embedded
support in a punctured ball about $x$.
Since almost minimization holds throughout $U$, the
construction admits a uniform positive outer radius
for centers near $x$. Repeating it at a nearby center
therefore gives regularity at $x$ as well.
Thus the support is smooth and properly embedded in
$M\setminus P$, and integrality and the constancy theorem
give the asserted locally constant positive integer
multiplicities.

Finally, choose $\rho_p>0$ so that
$B_{\rho_p}(p)\cap P=\{p\}$ and $V$ is almost minimizing
in every annulus centered at $p$ contained in this ball.
On its smooth support, almost minimization excludes
negative second variation by the deformation argument
of \cite[proof of Proposition~4.14, Step~3, Case~(a)]{DPDR2024},
with $c=0$.
Every compactly supported variation in
$B_{\rho_p}(p)\setminus\{p\}$ is supported in one such
annulus, proving stability there.
\end{proof}

\section{Removing singularities}
We adapt the separation argument of
\cite[Section~6]{DRHW2026} to near-area integrands in
dimensions $3\leq d\leq6$.
We compare the separation between neighboring sheets
with a small multiple of $|x|$.
A strict radial Jacobi barrier and stability, applied
with cutoffs around the puncture, give a linear lower
bound for this separation.
Together with the curvature estimate, this yields
the density bound needed for the removable-singularity
argument.

\begin{thm}
\label{thm:local-removable-singularity}
Let $3\le d\le6$ and $k_0\ge6$. There exists $\delta_*>0$
such that the following holds.
Let $F$ be a smooth positive elliptic integrand on
$G_{d-1}(B_1)$, where $B_1\subset\mathbb R^d$, satisfying
\[
 \|F-1\|_{C^{k_0}(G_{d-1}(B_1))}<\delta_*.
\]
Suppose that $\Sigma\subset B_1\setminus\{0\}$ is a smooth,
properly embedded, two-sided hypersurface without boundary which is $F$-stationary and $F$-stable away from the origin, and satisfies
$\mathcal H^{d-1}(\Sigma)<\infty$.
Then $\overline\Sigma\cap B_1$ is a smooth embedded,
$F$-stationary and $F$-stable hypersurface.
\end{thm}

\begin{proof}
Write $n=d-1$, $r(x)=|x|$, and
$d\mu_\Sigma=d\mathcal H^n\res\Sigma$. By Proposition~\ref{prop:stable-curvature-compactness},
\[
 |A_\Sigma(x)|\leq \frac K{r(x)}
\]
near the origin, with $K$ uniform in a fixed sufficiently
small neighborhood of the area integrand.
With our sign convention,
\[
 J_F=-\operatorname{div}_\Sigma(A_F\nabla_\Sigma)+c_F,
\]
and stability means
\[
Q_F(\varphi,\varphi)
 :=\int_\Sigma \varphi J_F\varphi\,d\mu_\Sigma=\int_\Sigma
 \left(
   \langle A_F\nabla_\Sigma\varphi,\nabla_\Sigma\varphi\rangle
   +c_F\varphi^2
 \right)d\mu_\Sigma
 \geq0
\]
for every $\varphi\in C_c^\infty(\Sigma)$.

Expanding the stationarity equation and the Jacobi
coefficients around the area integrand gives
\[
 \begin{gathered}
 |H|\leq C\delta_*(1+|A_\Sigma|),
 \qquad
 |A_F-\operatorname{Id}|\leq C\delta_*,
 \\
 |\nabla_\Sigma A_F|
 \leq C\delta_*(1+|A_\Sigma|),
 \qquad
 \bigl|c_F+|A_\Sigma|^2\bigr|
 \leq C\delta_*(1+|A_\Sigma|^2).
 \end{gathered}
\]
The formula for $\Delta_\Sigma r$ therefore yields
\[
 J_Fr\leq
 -\frac{n-|\nabla_\Sigma r|^2}{r}
 -r|A_\Sigma|^2
 +C\delta_*
 \left(
   \frac1r+1+|A_\Sigma|+r|A_\Sigma|^2
 \right).
\]
Since $|\nabla_\Sigma r|\leq1$ and $|A_\Sigma|\leq K/r$,
decreasing $\delta_*$ and then choosing $\rho>0$ small
gives
\begin{equation}\label{eq:strict-radial-barrier}
 J_Fr\leq-\frac br
 \qquad\text{on }\Sigma\cap B_\rho
\end{equation}
for some $b>0$.

The curvature estimate and interior elliptic estimates
give graphical neighborhoods of radius $\theta r(x)$,
with small slope and uniform rescaled derivative bounds.
Following the separation construction of
\cite[Propositions~6.7-6.8]{DRHW2026}, choose
$\kappa\ll\theta$ and let $s$ be the least positive
normal height of a neighboring sheet, allowing neighbors
on either side, wherever this height is less than
$\kappa r$. Denote this open set by
$U\subset\Sigma\cap B_\rho$.

For completeness, the estimates needed from this
construction are independent of the dimension.
A neighboring normal graph $u>0$ satisfies a uniformly
elliptic linear equation obtained by subtracting the
two graphical stationarity equations. Harnack and
Schauder estimates on smaller neighborhoods give
\[
 |\nabla_\Sigma u|\leq C\frac ur,
 \qquad
 |\nabla_\Sigma^2u|\leq C\frac u{r^2}.
\]
Taylor expansion of the normal-graph equation, whose
linearization at zero is $J_F$, consequently gives
\[
 |J_Fu|\leq C\left(\frac{u^2}{r^3}+\frac ur\right).
\]
Proper embeddedness makes the relevant family of
neighboring sheets locally finite. Taking its minimum,
and using the weak minimum principle for $J_F$, gives
\[
 \begin{gathered}
 s=\kappa r
 \;\text{on }\partial_\Sigma U\cap B_\rho,
 \qquad
 |\nabla_\Sigma s|\leq C\frac sr,
 \\
 J_Fs\geq-C\left(\frac{s^2}{r^3}+\frac sr\right)
 \;\text{weakly on }U.
 \end{gathered}
\]
Here $\kappa$ is chosen sufficiently small that
neighboring graphs continue across the smaller
graphical neighborhoods while their heights remain
below $\kappa r$.

Proper embeddedness gives positive separation on a
compact outer annulus. Choose $0<a<\kappa$ so that
$s>ar$ near $U\cap\partial B_\rho$.
After decreasing $a$ and $\rho$, on $\{s<ar\}$ we have
\[
 J_F(ar-s)
 \leq-\frac{ab}{r}
      +C\left(\frac{a^2}{r}+a\right)
 \leq-\frac{ab}{2r}.
\]
Set $w=(ar-s)_+$ on $U$ and extend it by zero.
Then $0\leq w\leq ar$, and $w$ vanishes near the
outer boundary.

Let $\chi_\varepsilon$ vanish on $B_\varepsilon$,
equal one outside $B_{2\varepsilon}$, and satisfy
$|\nabla_\Sigma\chi_\varepsilon|\leq C/\varepsilon$.
The weak positive-part inequality and stability give
\[
 \begin{aligned}
 0
 &\leq Q_F(\chi_\varepsilon w,\chi_\varepsilon w)\leq
 -\frac{ab}{2}
 \int_\Sigma
      \chi_\varepsilon^2\frac wr\,d\mu_\Sigma
 +\int_\Sigma
      w^2
      \langle
        A_F\nabla_\Sigma\chi_\varepsilon,
        \nabla_\Sigma\chi_\varepsilon
      \rangle\,d\mu_\Sigma\\
 &\leq
 -\frac{ab}{2}
 \int_\Sigma
      \chi_\varepsilon^2\frac wr\,d\mu_\Sigma
 +Ca^2\mu_\Sigma
      (B_{2\varepsilon}\setminus B_\varepsilon).
 \end{aligned}
\]
Finite area makes the last term tend to zero.
Thus $w=0$, proving
\begin{equation}\label{eq:linear-sheet-separation}
 s(x)\geq a|x|.
\end{equation}

Together with the graphical neighborhoods, this gives
$\eta_0>0$ such that
$\Sigma\cap B_{\eta_0|x|}(x)$ contains only one
graphical sheet, with uniformly bounded slope.
Covering each dyadic annulus by a bounded number of
such balls and summing gives
\[
 \mu_\Sigma(B_t(0))\leq Ct^n.
\]
For balls centered on $\Sigma$, use a graphical
neighborhood when $t\ll|y|$, and inclusion in
$B_{Ct}(0)$ otherwise. Enlarging intersecting balls
also covers arbitrary centers. Hence, in a smaller
neighborhood of the origin,
\begin{equation}\label{eq:removable-density-bound}
 \mu_\Sigma(B_t(y))\leq Ct^n.
\end{equation}
The constant may depend on $\Sigma$ through $a$.

This bound gives cutoffs around the puncture with
Dirichlet energy tending to zero, using logarithmic
cutoffs when $n=2$.
The near-area stability inequality also gives
\[
 \int_\Sigma |A_\Sigma|^2\varphi^2\,d\mu_\Sigma
 \leq C\int_\Sigma
       \left(|\nabla_\Sigma\varphi|^2+\varphi^2\right)
       d\mu_\Sigma.
\]
Applying this to an outer cutoff multiplied by the
cutoffs around the puncture shows that
$|A_\Sigma|^2$ is locally integrable across the origin.
The usual cutoff arguments therefore extend
stationarity and stability across the puncture.

If $n=2$, we have in particular
\[
 \int_{\Sigma\cap B_{\rho/2}}
      |A_\Sigma|^2\,d\mu_\Sigma<\infty.
\]
White's removable-singularity theorem applies as in
\cite[proof of Theorem~5.1, Step~3]{DPDRLi2024},
and gives the conclusion.

Suppose now that $3\leq n\leq5$.
Fix an arbitrary sequence $r_j\downarrow0$ and put
$\Sigma_j=r_j^{-1}\Sigma$.
The rescaled integrands $F(r_jx,T)$ converge smoothly
on compact sets to the autonomous integrand $F_0(T):=F(0,T)$.
By the curvature and density estimates and
Proposition~\ref{prop:stable-curvature-compactness},
a subsequence converges smoothly with integer
multiplicities on compact subsets of
$\mathbb R^{n+1}\setminus\{0\}$ to an $F_0$-stable
hypersurface $\Sigma_\infty$, possibly empty.
Its support is properly embedded and two-sided in
the punctured ambient space. Writing
$d\mu_\infty=d\mathcal H^n\res\Sigma_\infty$, we have
\[
 \mu_\infty(B_R(0))\leq CR^n,
 \qquad
 |A_{\Sigma_\infty}(x)|\leq K|x|^{-1},
 \qquad R>0.
\]

We spell out the metric and measure in the curvature
calculation on $\Sigma_\infty$.
Let $g_\infty$ be its induced Euclidean metric,
$\nu$ a unit normal, and $S=-d\nu$ its Weingarten
endomorphism. Define
\[
 \phi_0(\xi)=|\xi|F_0(\xi^\perp),
 \qquad
 B_0=D_\xi^2\phi_0(\nu)\big|_{T\Sigma_\infty}.
\]
Ellipticity makes $B_0$ positive definite.
The auxiliary metric and its volume measure are
\[
 g_{F_0}(X,Y)=g_\infty(B_0^{-1}X,Y),
 \qquad
 d\mu_{F_0}
   =(\det_{g_\infty}B_0)^{-1/2}\,d\mu_\infty.
\]
In particular, $d\mu_{F_0}$ is distinct from the
anisotropic energy measure
\[
 F_0(\nu^\perp)\,d\mu_\infty.
\]

Write $\nabla_{F_0}$ and $\Delta_{F_0}$ for the
gradient and Laplace-Beltrami operator of $g_{F_0}$.
Explicitly,
\[
 \begin{gathered}
\nabla_{F_0}f=B_0\nabla_{\Sigma_\infty}f,
 \quad 
 \langle\nabla_{F_0}f,\nabla_{F_0}h\rangle_{F_0}
 =g_\infty
   (B_0\nabla_{\Sigma_\infty}f,\nabla_{\Sigma_\infty}h),
 \\
 |\nabla_{F_0}f|_{F_0}^2
 =g_\infty
   (B_0\nabla_{\Sigma_\infty}f,\nabla_{\Sigma_\infty}f).
 \end{gathered}
\]
Define the anisotropic Weingarten endomorphism and
its norm by
\[
 S_{F_0}=B_0S,
 \qquad
 v=|S_{F_0}|_{F_0}.
\]
Since $S_{F_0}$ is self-adjoint for $g_{F_0}$,
\[
 v^2
 =\sum_{\alpha=1}^n
   g_{F_0}(S_{F_0}e_\alpha,S_{F_0}e_\alpha)
 =\operatorname{tr}(B_0SB_0S),
\]
where $\{e_\alpha\}$ is $g_{F_0}$-orthonormal.
All these metrics, measures, and curvature norms
are uniformly equivalent to their Euclidean
counterparts when $\delta_*$ is small.

By \cite[equations~(10),(11)]{Winklmann2005},
the anisotropic curvature $v=|S_{F_0}|_{F_0}$ satisfies
\[
 \begin{gathered}
 \int_{\Sigma_\infty}v^2f^2\,d\mu_{F_0}
 \leq \sigma\int_{\Sigma_\infty}
       |\nabla_{F_0}f|_{F_0}^2\,d\mu_{F_0},
 \\
 \frac12\Delta_{F_0}(v^2)
 \geq(1+\gamma)|\nabla_{F_0}v|_{F_0}^2-\beta v^4,
 \end{gathered}
\]
where $\sigma,\beta,\gamma$ can be chosen arbitrarily
close to $1,1,2/n$, respectively, by taking $F_0$
sufficiently close to area.

We use these inequalities to obtain the critical
exponent $n$, including the cases $n=3,4$ not covered
by the exponent range in
\cite[equation~(3)]{Winklmann2005}.
Set
\[
 q=\frac{n-4}{2},
 \qquad
 v_\varepsilon=(v^2+\varepsilon^2)^{1/2}.
\]
For a smooth cutoff $\eta$ supported away from the
origin, test the regularized Simons inequality with
$v_\varepsilon^{2q}\eta^2$ and the stability inequality
with $v_\varepsilon^{q+1}\eta$.
Integration by parts and Young's inequality give
\[
 \int_{\Sigma_\infty}
v_\varepsilon^{2q+4}\eta^2\,d\mu_{F_0}\leq C\int_{\Sigma_\infty}v_\varepsilon^{2q+2}|\nabla_{F_0}\eta|_{F_0}^2\,d\mu_{F_0}+C\varepsilon^2\int_{\Sigma_\infty}v_\varepsilon^{2q+2}\eta^2\,d\mu_{F_0}.
\]
Indeed, the coefficient needed for absorption is
\[
 2q+1+\gamma-\beta\sigma(q+1)^2,
\]
which is positive for sufficiently small $\delta_*$
because its area value is
\[
 \frac2n-\frac{(n-4)^2}{4}>0,
 \qquad 3\leq n\leq5.
\]
The regularization justifies the calculation also
when $q=-1/2$.

Letting $\varepsilon\downarrow0$, taking
$\eta=\psi^{n/2}$, and applying Young's inequality
once more yields
\[
 \int_{\Sigma_\infty}v^n\psi^n\,d\mu_{F_0}
 \leq C\int_{\Sigma_\infty}
       |\nabla_{F_0}\psi|_{F_0}^n\,d\mu_{F_0}.
\]
Uniform equivalence of the metrics, curvature norms,
and measures therefore gives
\begin{equation}\label{eq:critical-curvature-integral}
 \int_{\Sigma_\infty}
       |A_{\Sigma_\infty}|^n\psi^n\,d\mu_\infty
 \leq C\int_{\Sigma_\infty}
       |\nabla_{\Sigma_\infty}\psi|^n\,d\mu_\infty.
\end{equation}

For $R>2$, choose a radial cutoff $\psi_R$ supported
where $R^{-2}<|x|<R^2$, equal to one where
$R^{-1}\leq|x|\leq R$, and satisfying
\[
 |\nabla_{\Sigma_\infty}\psi_R|
 \leq\frac{C}{|x|\log R}.
\]
The area-growth bound, summed over dyadic annuli,
gives
\[
 \int_{\Sigma_\infty}
       |\nabla_{\Sigma_\infty}\psi_R|^n\,d\mu_\infty
 \leq C(\log R)^{1-n}\longrightarrow0.
\]
Equation~\eqref{eq:critical-curvature-integral} and
Fatou's lemma imply
$A_{\Sigma_\infty}=0$ away from the origin.

Writing $d\mu_j=d\mathcal H^n\res\Sigma_j$,
smooth graphical convergence consequently gives,
for every fixed $\tau\in(0,1)$,
\[
 \int_{\Sigma_j\cap(B_1\setminus B_\tau)}
      |A_{\Sigma_j}|^2\,d\mu_j\longrightarrow0.
\]
On the inner ball, the rescaled curvature and density
bounds instead give
\[
 \int_{\Sigma_j\cap B_\tau}
      |A_{\Sigma_j}|^2\,d\mu_j
 \leq CK^2\tau^{n-2}.
\]
Letting first $j\to\infty$ and then $\tau\downarrow0$,
and using that the initial sequence of radii was
arbitrary, we conclude
\[
 r^{2-n}\int_{\Sigma\cap B_r}
             |A_\Sigma|^2\,d\mu_\Sigma
 \longrightarrow0
 \qquad\text{as }r\downarrow0.
\]

Apply the regularity criterion in
\cite[proof of Theorem~A.1, equation~(A.4)]{DPDRLi2024}
to the constant sequence $\Sigma$.
The mass ratios are uniformly bounded, stationarity
and stability extend across the puncture, and
$\mathcal H^{n-2}(\{0\})=0$.
The preceding limit excludes the origin from the
curvature-concentration set, giving smooth extension.

Stationarity and stability persist on the extension.
The ordering of the graphical sheets and the strong
maximum principle give embeddedness.
\end{proof}

Now we are ready to prove Theorem~\ref{thm:near-area-minmax} and conclude with the section.

\begin{proof}[Proof of Theorem~\ref{thm:near-area-minmax}]
Choose $\delta_0>0$ sufficiently small for the preceding
local regularity results.
Proposition~\ref{prop:pull-tight-almost-minimizing} gives
representatives $\Phi_j\in\Pi_j$ satisfying
\[
 \sup_{x\in X_j}\F(\Phi_j(x))
 \leq \mathbf L(\Pi_j)+j^{-1}
\]
and an $F$-stationary varifold
$V\in\mathbf C(\{\Phi_j\})$ with $\F(V)=L>0$.
By the definition of the critical set, after passing to
a subsequence there are $x_j\in X_j$ such that
$|\Phi_j(x_j)|\rightharpoonup V$.

By Corollary~\ref{cor:regularity-away-finite-points}
and its proof, $V$ is integral and there is a finite
set $P\subset M$ such that
\[
 V\res(M\setminus P)=\theta|\Sigma|,
 \qquad
 \mathcal H^{d-1}(\Sigma)\leq\|V\|(M)<\infty,
\]
where $\Sigma$ is smooth, properly embedded and
$F$-stationary, and $\theta$ is locally constant
with positive integer values.
Moreover, $\Sigma$ is stable in a sufficiently small
punctured ball about each $p\in P$.

Fix $p\in P$.
Since $d\geq3$, this punctured coordinate ball is
simply connected, so proper embeddedness implies
that $\Sigma$ is two-sided there.
In orthonormal normal coordinates, rescaled by
$y\mapsto\exp_p(ry)$, absorb the metric area Jacobian
into the integrand and divide the functional by
$r^{d-1}$. The resulting Euclidean integrand satisfies
\[
 \|F_{p,r}-1\|_{C^{k_0}(G_{d-1}(B_1))}
 \leq C\delta_0+o(1)
 \qquad\text{as }r\downarrow0,
\]
uniformly in $p$.
After decreasing $\delta_0$ and then choosing $r$
small, Theorem~\ref{thm:local-removable-singularity}
therefore applies and extends $\Sigma$ smoothly
across $p$, together with its integer multiplicity.

Applying this at every point of $P$ gives a smooth
closed embedded $F$-stationary hypersurface.
Since $V$ is integral, $\|V\|(P)=0$, so the extended
hypersurface with its multiplicities represents
the whole varifold.
Its compactness gives finitely many connected
components, and hence
\[
 V=\sum_{i=1}^s m_i|\Sigma_i|,
 \qquad m_i\in\mathbb N,
 \qquad \F(V)=L,
\]
with the asserted properties.
All smallness requirements on $\delta_0$ come from
the local analytic results and the coordinate estimates,
and are independent of the parameter dimension,
the homotopy classes, and their min-max values.
\end{proof}


\begin{thebibliography}{10}

\bibitem{Allard1974}
{\sc Allard, W.~K.}
\newblock A characterization of the area integrand.
\newblock In {\em Symposia Mathematica, Vol. XIV}. 1974, pp.~429--444.

\bibitem{Allard1983}
{\sc Allard, W.~K.}
\newblock An a priori estimate for the oscillation of the normal to a
  hypersurface whose first and second variation with respect to an elliptic
  integrand is controlled.
\newblock {\em Invent. Math. 73}, 2 (1983), 287--331.

\bibitem{Almgren-homotopy}
{\sc Almgren, Jr., F.~J.}
\newblock The homotopy groups of the integral cycle groups.
\newblock {\em Topology 1\/} (1962), 257--299.

\bibitem{Bellettini2025}
{\sc Bellettini, C.}
\newblock Extensions of {Schoen--Simon--Yau} and {Schoen--Simon} theorems via
  iteration {\`a} la {De Giorgi}.
\newblock {\em Inventiones Mathematicae 240}, 1 (2025), 1--34.

\bibitem{CatinoMastroliaRoncoroni2024}
{\sc Catino, G., Mastrolia, P., and Roncoroni, A.}
\newblock Two rigidity results for stable minimal hypersurfaces.
\newblock {\em Geometric and Functional Analysis 34}, 1 (2024), 1--18.

\bibitem{Chambers-Liokumovich2020}
{\sc Chambers, G.~R., and Liokumovich, Y.}
\newblock Existence of minimal hypersurfaces in complete manifolds of finite
  volume.
\newblock {\em Invent. Math. 219}, 1 (2020), 179--217.

\bibitem{ChambersRotman}
{\sc Chambers, G.~R., and Rotman, R.}
\newblock Monotone homotopies and contracting discs on riemannian surfaces.
\newblock {\em Journal of Topology and Analysis 10}, 2 (2018), 323--354.

\bibitem{ChenLan2026}
{\sc Chen, L., and Lan, J.}
\newblock Multiplicity one for strictly stable anisotropic allen-cahn minimal
  hypersurfaces, 2026.

\bibitem{ChodoshLi2023}
{\sc Chodosh, O., and Li, C.}
\newblock Stable anisotropic minimal hypersurfaces in {$\mathbf{R}^{4}$}.
\newblock {\em Forum Math. Pi 11\/} (2023), e3.

\bibitem{ChodoshLi2024}
{\sc Chodosh, O., and Li, C.}
\newblock Stable minimal hypersurfaces in {$\mathbf{R}^{4}$}.
\newblock {\em Acta Mathematica 233}, 1 (2024), 1--31.

\bibitem{ChodoshLiMinterStryker2026}
{\sc Chodosh, O., Li, C., Minter, P., and Stryker, D.}
\newblock Stable minimal hypersurfaces in {$\mathbf{R}^{5}$}.
\newblock {\em Annals of Mathematics 204}, 1 (2026), 423--453.

\bibitem{Chodosh-Liokumovich-Spolaor2022}
{\sc Chodosh, O., Liokumovich, Y., and Spolaor, L.}
\newblock Singular behavior and generic regularity of min-max minimal
  hypersurfaces.
\newblock {\em Ars Inven. Anal.\/} (2022), Paper No. 2, 27.

\bibitem{ChodoshMantoulidis2020}
{\sc Chodosh, O., and Mantoulidis, C.}
\newblock Minimal surfaces and the {Allen--Cahn} equation on {$3$}-manifolds:
  index, multiplicity, and curvature estimates.
\newblock {\em Annals of Mathematics 191}, 1 (2020), 213--328.

\bibitem{ChodoshMantoulidis2023}
{\sc Chodosh, O., and Mantoulidis, C.}
\newblock The {$p$}-widths of a surface.
\newblock {\em Publications Math{\'e}matiques de l'IH{\'E}S 137\/} (2023),
  245--342.

\bibitem{DKV2023}
{\sc Daneshpajouh, H.~R., Karasev, R., and Volovikov, A.~Y.}
\newblock {Hedetniemi}'s conjecture from the topological viewpoint.
\newblock {\em Journal of Combinatorial Theory, Series A 195\/} (2023), 105721.

\bibitem{DPDR2024}
{\sc De~Philippis, G., and De~Rosa, A.}
\newblock The anisotropic min-max theory: Existence of anisotropic minimal and
  {CMC} surfaces.
\newblock {\em Comm. Pure Appl. Math. 77}, 7 (2024), 3184--3226.

\bibitem{DDH2019}
{\sc De~Philippis, G., De~Rosa, A., and Hirsch, J.}
\newblock The area blow up set for bounded mean curvature submanifolds with
  respect to elliptic surface energy functionals.
\newblock {\em Discrete Contin. Dyn. Syst. 39}, 12 (2019), 7031--7056.

\bibitem{DPDRLi2024}
{\sc De~Philippis, G., De~Rosa, A., and Li, Y.}
\newblock Existence and regularity of min-max anisotropic minimal
  hypersurfaces, 2024.
\newblock Preprint, arXiv:2409.15232v1.

\bibitem{DeR-survey}
{\sc De~Rosa, A.}
\newblock On the theory of anisotropic minimal surfaces.
\newblock {\em Notices Amer. Math. Soc. 71}, 7 (2024), 853--859.

\bibitem{DeRosa2026Lectures}
{\sc De~Rosa, A.}
\newblock Lectures on existence and regularity of anisotropic minimal surfaces.
\newblock Preprint, arXiv:2608.10275.

\bibitem{DRHW2026}
{\sc De~Rosa, A., Halavati, A., and Wang, L.}
\newblock Min-max construction of anisotropic minimal surfaces with genus
  bound, 2026.
\newblock Preprint, arXiv:2606.30098.

\bibitem{DeRosaPigati2025}
{\sc De~Rosa, A., and Pigati, A.}
\newblock Anisotropic min-max via phase transitions.
\newblock Preprint, arXiv:2512.07039.

\bibitem{Dey2025}
{\sc Dey, A.}
\newblock A sub-additive inequality for the volume spectrum.
\newblock {\em Journal of Topology and Analysis 17}, 5 (2025), 1543--1561.

\bibitem{doCarmoPeng1979}
{\sc do~Carmo, M., and Peng, C.}
\newblock Stable complete minimal surfaces in {$\mathbb{R}^{3}$} are planes.
\newblock {\em Bulletin of the American Mathematical Society (N.S.) 1}, 6
  (1979), 903--906.

\bibitem{FischerColbrieSchoen1980}
{\sc Fischer-Colbrie, D., and Schoen, R.}
\newblock The structure of complete stable minimal surfaces in {$3$}-manifolds
  of non-negative scalar curvature.
\newblock {\em Communications on Pure and Applied Mathematics 33}, 2 (1980),
  199--211.

\bibitem{Givental1990}
{\sc Givental, A.~B.}
\newblock Nonlinear generalization of the {Maslov} index.
\newblock In {\em Theory of Singularities and its Applications}, vol.~1 of {\em
  Advances in Soviet Mathematics}. American Mathematical Society, Providence,
  RI, 1990, pp.~71--103.

\bibitem{GKPS2021}
{\sc Granja, G., Karshon, Y., Pabiniak, M., and Sandon, S.}
\newblock {Givental}'s non-linear {Maslov} index on lens spaces.
\newblock {\em International Mathematics Research Notices 2021}, 23 (2021),
  18225--18299.

\bibitem{Gromov1988Dimension}
{\sc Gromov, M.}
\newblock Dimension, non-linear spectra and width.
\newblock In {\em Geometric Aspects of Functional Analysis}, J.~Lindenstrauss
  and V.~D. Milman, Eds., vol.~1317 of {\em Lecture Notes in Mathematics}.
  Springer, Berlin, Heidelberg, 1988, pp.~132--184.

\bibitem{Guth2009Minimax}
{\sc Guth, L.}
\newblock Minimax problems related to cup powers and {Steenrod} squares.
\newblock {\em Geom. Funct. Anal. 18}, 6 (2009), 1917--1987.

\bibitem{Hardt1977}
{\sc Hardt, R.}
\newblock On boundary regularity for integral currents or flat chains modulo
  two minimizing the integral of an elliptic integrand.
\newblock {\em Communications in Partial Differential Equations 2}, 12 (1977),
  1163--1232.

\bibitem{HLW-bernstein-7}
{\sc Hong, H., Li, H., and Wang, G.}
\newblock The stable {Bernstein} theorem in {$\mathbb{R}^{7}$}.
\newblock Preprint, {arXiv}:2609.15720 [math.{DG}] (2026), 2026.

\bibitem{IMN2018}
{\sc Irie, K., Marques, F.~C., and Neves, A.}
\newblock Density of minimal hypersurfaces for generic metrics.
\newblock {\em Ann. of Math. (2) 187}, 3 (2018), 963--972.

\bibitem{LiXia2026}
{\sc Li, J., and Xia, C.}
\newblock Stable anisotropic minimal hypersurfaces in {$\mathbb{R}^{5}$} and
  {$\mathbb{R}^{6}$}.
\newblock {\em J. Geom. Anal. 36\/} (2026), 312.
\newblock Article 312; preprint version arXiv:2505.16595v2.

\bibitem{LiZhou2021}
{\sc Li, M. M.-c., and Zhou, X.}
\newblock Min-max theory for free boundary minimal hypersurfaces {I}:
  Regularity theory.
\newblock {\em J. Differential Geom. 118}, 3 (2021), 487--553.

\bibitem{Li2023}
{\sc Li, Y.}
\newblock Existence of infinitely many minimal hypersurfaces in
  higher-dimensional closed manifolds with generic metrics.
\newblock {\em J. Differential Geom. 124}, 2 (2023), 381--395.
\newblock Cited theorem numbering refers to arXiv:1901.08440v2.

\bibitem{Lin1990}
{\sc Lin, F.-H.}
\newblock Estimates for surfaces which are stationary for an elliptic
  parametric integral.
\newblock {\em J. Partial Differential Equations 3}, 3 (1990), 78--92.

\bibitem{Weyl-LMN}
{\sc Liokumovich, Y., Marques, F.~C., and Neves, A.}
\newblock {Weyl} law for the volume spectrum.
\newblock {\em Ann. of Math. (2) 187}, 3 (2018), 933--961.
\newblock Cited theorem numbering refers to arXiv:1607.08721v2.

\bibitem{MN2017PositiveRicci}
{\sc Marques, F.~C., and Neves, A.}
\newblock Existence of infinitely many minimal hypersurfaces in positive
  {Ricci} curvature.
\newblock {\em Invent. Math. 209}, 2 (2017), 577--616.

\bibitem{MarquesNeves2021}
{\sc Marques, F.~C., and Neves, A.}
\newblock Morse index of multiplicity one min-max minimal hypersurfaces.
\newblock {\em Advances in Mathematics 378\/} (2021), 107527.

\bibitem{MNS}
{\sc Marques, F.~C., Neves, A., and Song, A.}
\newblock Equidistribution of minimal hypersurfaces for generic metrics.
\newblock {\em Invent. Math. 216}, 2 (2019), 421--443.
\newblock Cited theorem numbering refers to arXiv:1712.06238v2.

\bibitem{Mazet2024}
{\sc Mazet, L.}
\newblock Stable minimal hypersurfaces in {$\mathbb{R}^{6}$}, 2024.

\bibitem{Pitts1981}
{\sc Pitts, J.~T.}
\newblock {\em Existence and Regularity of Minimal Surfaces on {Riemannian}
  Manifolds}, vol.~27 of {\em Mathematical Notes}.
\newblock Princeton University Press, Princeton, NJ, 1981.

\bibitem{Pogorelov1981}
{\sc Pogorelov, A.~V.}
\newblock On the stability of minimal surfaces.
\newblock {\em Doklady Akademii Nauk SSSR 260}, 2 (1981), 293--295.
\newblock English translation: Soviet Mathematics Doklady 24 (1981), 274--276.

\bibitem{Reshetnyak1968}
{\sc Reshetnyak, Y.~G.}
\newblock Weak convergence of completely additive vector functions on a set.
\newblock {\em Siberian Mathematical Journal 9}, 6 (1968), 1039--1045.

\bibitem{SS1981}
{\sc Schoen, R., and Simon, L.}
\newblock Regularity of stable minimal hypersurfaces.
\newblock {\em Comm. Pure Appl. Math. 34}, 6 (1981), 741--797.

\bibitem{SSA1977}
{\sc Schoen, R., Simon, L., and Almgren, Jr., F.~J.}
\newblock Regularity and singularity estimates on hypersurfaces minimizing
  parametric elliptic variational integrals.
\newblock {\em Acta Math. 139\/} (1977), 217--265.

\bibitem{SSY1975}
{\sc Schoen, R., Simon, L., and Yau, S.-T.}
\newblock Curvature estimates for minimal hypersurfaces.
\newblock {\em Acta Mathematica 134}, 3-4 (1975), 275--288.

\bibitem{Smale}
{\sc Smale, S.}
\newblock An infinite dimensional version of {Sard}'s theorem.
\newblock {\em Amer. J. Math. 87}, 4 (1965), 861--866.

\bibitem{Song2023}
{\sc Song, A.}
\newblock Existence of infinitely many minimal hypersurfaces in closed
  manifolds.
\newblock {\em Ann. of Math. (2) 197}, 3 (2023), 859--895.

\bibitem{SongZhou2021}
{\sc Song, A., and Zhou, X.}
\newblock Generic scarring for minimal hypersurfaces along stable
  hypersurfaces.
\newblock {\em Geom. Funct. Anal. 31}, 4 (2021), 948--980.

\bibitem{White1991}
{\sc White, B.}
\newblock The space of minimal submanifolds for varying {Riemannian} metrics.
\newblock {\em Indiana Univ. Math. J. 40}, 1 (1991), 161--200.

\bibitem{WhiteBumpy2017}
{\sc White, B.}
\newblock On the bumpy metrics theorem for minimal submanifolds.
\newblock {\em American Journal of Mathematics 139}, 4 (2017), 1149--1155.

\bibitem{Winklmann2005}
{\sc Winklmann, S.}
\newblock Pointwise curvature estimates for {$F$}-stable hypersurfaces.
\newblock {\em Ann. Inst. H. Poincar{\'e} Anal. Non Lin{\'e}aire 22}, 5 (2005),
  543--555.

\bibitem{Yau1982}
{\sc Yau, S.-T.}
\newblock Problem section.
\newblock In {\em Seminar on Differential Geometry}, S.-T. Yau, Ed., vol.~102
  of {\em Annals of Mathematics Studies}. Princeton University Press,
  Princeton, NJ, 1982, pp.~669--706.

\bibitem{Zhou}
{\sc Zhou, X.}
\newblock On the multiplicity one conjecture in min-max theory.
\newblock {\em Ann. of Math. (2) 192}, 3 (2020), 767--820.
\newblock Cited theorem numbering refers to arXiv:1901.01173v2.

\bibitem{ZhouZhu2019}
{\sc Zhou, X., and Zhu, J.~J.}
\newblock Min-max theory for constant mean curvature hypersurfaces.
\newblock {\em Invent. Math. 218}, 2 (2019), 441--490.

\bibitem{HalavatiCubeToSlab}
\newblock From cubes to a long rectangle.
\newblock Interactive visualization.
\newblock \url{https://ariahala.github.io/cube-to-slab.html}.

\end{thebibliography}
\end{document}